\documentclass[12pt, letterpaper]{article}
\usepackage{geometry}
\usepackage{graphicx}
\usepackage{amssymb}
\usepackage{epstopdf}
\usepackage{amsmath}
\usepackage{amsmath, ams cd}
\usepackage{pictexwd,dcpic}
\usepackage{rcs}
\usepackage{slashbox}

\usepackage{amsthm, amssymb, amscd}
\pdfoutput=1

\newtheorem{prop}{Proposition}[section]
\newtheorem{theorem}[prop]{Theorem}
\newtheorem{lemma}[prop]{Lemma}
\newtheorem{claim}[prop]{Claim}

\newtheorem{remark}[prop]{Remark}

 \def\G{{\Gamma}}
 \def\d{{\delta}}
 
 \def\e{{\epsilon}}
 \def\l{{\lambda}}
 \def\L{{\Lambda}}

 \def\Th{{\Theta}}
 \def\s{{\sigma}}
 
 \def\a{{\alpha}}
 \def\b{{\beta}}
 \def\p{{\partial}}
 
 \def\ra{{\rightarrow}}
 
 \def\g{{\gamma}}

 \def\Ra{{\Rightarrow}}
 
 \def\c{{\mathbb C}}
 
 \def\z{{\mathbb Z}}
 \def\2{{\mathbb Z_2}}
 
 \def\t{{\tau}}

 \def\Th{{\Theta}}

 \def\sl2{{SL(2,\mathbb C)}}

 \def\ni{{\noindent}}

 \def\sl{{{\mbox{\tiny $\L$}}}}

\usepackage{pictexwd,dcpic}
\usepackage{rcs}

\def\ni{\noindent}

\def\d{\delta}

\def\a{\alpha}
\def\b{\beta}
\def\t{\tau}
\def\L{\Lambda}
\def\e{\epsilon}

\def\l{\lambda}
\def\g{\gamma}

\def\s{\sigma}

\begin{document}

\title{Two classes of virtually fibered Montesinos links of type $\widetilde{SL_2}$}

\author{Xiao Guo\thanks{E--mail: xiaoguo@buffalo.edu} \\University at Buffalo, SUNY}

\date{ } 

\maketitle

\begin{abstract}
We find  two new classes of virtually fibered classic Montesinos links
of type $\widetilde{SL_2}$.
\end{abstract}

\section{Introduction}

A $3$-manifold is called  {\it fibered} if it can be given the structure of a surface bundle over the circle.  
If a $3$-manifold can be finitely covered by a fibered $3$-manifold, we call it {\it virtually fibered}. The {\it virtually fibered conjecture} states that every complete hyperbolic 3-manifold with finite volume is virtually fibered. This conjecture was proposed by Thurston in 1982 as a question in \cite{thu}. A {\it link} in $S^3$ is called  {\it virtually fibered} if its exterior is virtually fibered.

\begin{figure}
\begin{center}
\includegraphics{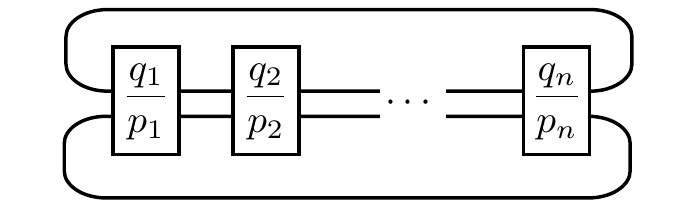}
\end{center}
\caption{\label{cml} Classic Montesinos link}
\end{figure}

 A {\it classic Montesinos link}  has a projection as shown in Figure  \ref{cml}, where a small rectangle with $q_i/p_i$ stands for a rational tangle, $1\leq i\leq n$.
 $n$, $q_i$ and $p_i$ are integers, and we may assume that $n\geq 1$, $p_i\geq 2$, and $q_i$ and $p_i$ are relatively prime, $1\leq i\leq n$. 
By assumption, the absolute value of $q_i/p_i$ may be greater than one, $1\leq i\leq n$.
 $(q_1/p_1,q_2/p_2,...,q_n/p_n)$ is called a {\it cyclic rational tangle decomposition} of the classic Montesinos link which has a projection as shown in Figure \ref{cml}.
Note that a classic Montesinos link may have different cyclic rational tangle decompositions.
Montesinos link $K$ in $S^3$ has a Seifert fibered 2-fold branched covering space, $W_K$. The branch set $\widetilde{K}$ is the preimage of $K$.
$(W_K, \widetilde{K})\ra (S^3, K)$ is induced by the homomorphism $\pi_1(S^3-\overset{\circ}{N}(K))\ra \z_2$, where the image of meridianal generators in $\pi_1(S^3-\overset{\circ}{N}(K))$ is $\bar{1}$. 
$K$ is a Montesinos link of {\it type} $\widetilde{SL_2}$ if $W_K$ has  $\widetilde{SL_2}$ geometric structure. 
 
Recently, the virtually fibered conjecture has been solved  for classic Montesinos links which are not of type $\widetilde{SL_2}$ due to the works of Walsh \cite{wa} and Agol-Boyer-Zhang \cite{abz}. 
Agol-Boyer-Zhang also gave an infinite family of virtually fibered classic Montesinos links of type $\widetilde{SL_2}$, 
in Sec. 6 of \cite{abz}. 
Those links have cyclic rational tangle decompositions $\{({q_1}/{p}, {q_2}/{p}, \ldots, {q_n}/{p}): 
p\geq 3 \text{ odd and } n \text{ is a multiple of } p\}$. Later, this result is extended by removing the condition that $n$ is a multiple of $p$ in \cite{gz}. 
Note that a classic Montesinos link in the set $\{({q_1}/{p}, {q_2}/{p}, \ldots, {q_n}/{p}): p\geq 3 \text{ odd }\}$
 is of type $\widetilde{SL_2}$ when $q_1+q_2+\cdots+q_n\neq0$, and $n=3, q\geq5$, or $n>3$. 
In this paper, we give another two families of virtually fibered classic Montesinos links of type $\widetilde{SL_2}$ by extending the techniques used in Sec. 6 of \cite{abz}.

\begin{theorem}\label{main}
If $K$ is a classic Montesinos link with a cyclic rational tangle decomposition of one of the following forms:

I. $(\displaystyle{\frac{q_1}{p},
\frac{q_2}{p}, \ldots, \frac{q_k}{p},\frac{q_{k+1}}{pr}, \frac{q_{k+2}}{pr}, \ldots, \frac{q_{n}}{pr}})$ where  $p,r\geq 3$ odd, $k$ is a multiple of $p$, and $(n-k)$ is a multiple of $p$ and $r$;

II. $(\displaystyle{\frac{q_1}{p},
\frac{q_2}{p}, \ldots, \frac{q_n}{p}})$ where  $n\geq 4$ even, $p=2m$, and $m$ is odd,\\
then $K$ is virtually fibered.
\end{theorem}
Note that $K$ is of type $\widetilde{SL_2}$ if $r(q_1+\cdots+q_k)+q_{k+1}+\cdots +q_n\neq 0$ in Case I, and $q_1+\cdots+q_n\neq 0$ and $(p,n)\neq (2,4)$ in Case II.

In Case I, the rational tangles in the cyclic rational tangle decomposition of $K$ have different denominators.
In Case II, the denominators are same but they are even numbers.
There are some new issues raised by these two new situations.
We shall mention them later.

For convenience, we follow the notations used in Sec. 6 of \cite{abz}. 
Let $K=K(q_1/p_1,q_2/p_2,\ldots, q_n/p_n)$ be the classic Montesinos link with cyclic rational tangle decomposition $(q_1/p_1,q_2/p_2,\ldots, q_n/p_n)$. 
Let $\mathcal{B}_K$ be the base orbifold of the Seifert fibered space $W_K$.
$W_K$ has the $\widetilde{SL_2}$ geometry structure if
\begin{equation}\label{ew}
e(W_K)=-\sum_{i=1}^{n}\frac{q_i}{p_i}\neq 0;\ \ \ 
\chi(\mathcal{B}_K)=2-n+\sum_{i=1}^{n}\frac{1}{p_i}<0,
\end{equation} 
where $e(W_k)$ is the Euler number of $W_K$, and $\chi({\mathcal{B}_K})$ is the Euler characteristic of $\mathcal{B}_K$. 
$\widetilde{K}$ is geodesic in $W_K$ and orthogonal to Seifert fibers of $W_K$ by Lemma 2.1 of \cite{abz}. 
Let $f:W_K\rightarrow \mathcal{B}_K$
be the Seifert quotient map. 
By the construction of $W_K$, ${\cal B}_K$ is a $2$-sphere with $n$ cone points $\{c_1,c_2,\ldots, c_n\}$, and the order of $c_i$ is $p_i$, $1\leq i\leq n$. 
Let $K^*=f(\widetilde{K})$. Then $K^*$ is a geodesic equator of ${\cal B}_K$ containing all the $n$ cone points.
The number of components of $K$ is also decided by the cyclic rational tangle decomposition.
\begin{equation}\label{cn}
|K|=
\begin{cases}
 1     & \text{if each $p_i$ is odd and $(q_1+\cdots+q_n)$ is odd }, \\
 2    & \text{if each $p_i$ is odd and $(q_1+\cdots+q_n)$ is even},\\
  \#\{i:p_i \text{ is even}\}   & \text{otherwise}.
\end{cases}
\end{equation}

We proof Case I of Theorem \ref{main} in section 2, and  Case II in section 3.

\section{Proof of Theorem \ref{main} in Case I.}
In this section we prove Theorem 1.1 when 
$K=K(\displaystyle{\frac{q_1}{p},
\frac{q_2}{p}, \ldots, \frac{q_k}{p},\frac{q_{k+1}}{pr}, \frac{q_{k+2}}{pr}, \ldots, \frac{q_{n}}{pr}})$ with  $p,r\geq 3$ odd, where $k$ is a multiple of $p$ and $(n-k)$ is a multiple of $p$ and $r$. 
We consider $K$ is of type $\widetilde{SL_2}$, and $k\neq 0$, $n\neq k$. 
Otherwise, $K$ is virtually fibered by \cite{abz}. 
By (\ref{ew}), $K$ is of type $\widetilde{SL_2}$ when $e(W_K)=-(r(q_1+\cdots+q_k)+q_{k+1}+\cdots+q_n)/pr\neq 0$.
From (\ref{cn}), $K$ has one or two components. We first prove Theorem 1.1 when $K$ is a knot in Sec. 2.1. 
As in Sec. 6.2 of \cite{abz}, the proof can be extended to the case that $K$ has two components with several adjustments in Sec. 2.2.

\subsection{$K$ is a knot.}

We first give the outline of the proof. 

At first, we construct a finite cover of $W_K$, say $Y$, so that $Y$ is a locally-trivial circle bundle. 
Let $L$ be the preimage of $\widetilde{K}$ in $Y.$ 
We prove that the exterior of some components of $L$ in $Y$, denoted by $M$, is a surface semi-bundle. 
Next, we construct $\breve{M}$, which is a 2-fold cover of $M$, such that $\breve{M}$ is a fibered manifold. 
This covering can be extended to $Y$. 
Suppose that the corresponding 2-fold cover of $Y$ and $L$ are $\breve{Y}$ and $\breve{L}$ respectively.
 We can isotope $\breve{L}\cap\breve{M}$ and perform Dehn twists to the surface fibers of $\breve{M}$ such that  the components of $\breve{L}\cap\breve{M}$ are transverse to the new surface fibers of $\breve{M}$. 
Then the exterior of $\breve{L}$ in $\breve{Y}$ has a structure of surface bundle over the circle.
In addition, the exterior of $\breve{L}$ in $\breve{Y}$ is a finite cover of the exterior of $\widetilde{K}$ in $W_K$, so it is also a finite cover of the exterior of $K$ in $S^3$. 
Therefore, $K$ is virtually fibered.

Different from \cite{abz}, the tangles in the cyclic rational tangle decomposition of $K$ have different denominators.
 We need compose two finite covers to build $Y$, and 
 perform additional Dehn twist operations to the surface fibers of $\breve{M}$.
 
Recall that $f:W_K\ra{\cal B}_K$ is the Seifert quotient map and $K^*=f(\widetilde{K})$. 
$\mathcal{B}_K$ is a 2-sphere with $n$ cone points.
$K^*$ goes through all cone points of $\mathcal{B}_K:\{c_1,c_2,\ldots,c_n\}$ successively. 
$c_i$ has order $p$ if $1\leq i\leq k$, and has order $pr$ if $k+1\leq i\leq n$.
 $f|:\widetilde{K}\rightarrow K^*$ is a 2-fold cyclic cover because the order of $c_i$ is odd, $1\leq i\leq n$.

At first, we construct a proper covering map, $\psi: F\rightarrow \mathcal{B}_K$, such that $F$ is a smooth orientable closed surface. 
Since ${\cal B}_K$ contains cone points of different orders, we construct $\psi$ by composing two cyclic covers, $\psi_1$ and $\psi_2$.

Let  $\Gamma_1$ be the fundamental group of $\mathcal{B}_K$. 
Then
\begin{equation*}
\Gamma_1=\{x_1,x_2,\ldots,x_n: x_1^{p}=\cdots=x_k^{p}=x_{k+1}^{pr}=\cdots=x_n^{pr}=x_1x_2\cdots x_n=1\}.
\end{equation*}  
 $x_i$ is represented by a circle centered at $c_i$ and in a small regular neighborhood of $c_i$ in $\mathcal{B}_K$, $1\leq i\leq n$. 
$n-k$ is a multiple of $r$. There exists a homomorphism $h_1: \Gamma_1\rightarrow \z/r$ where
\begin{equation*}h_1(x_i)=
\begin{cases}
   \bar{0}   & \text{if $1\leq i\leq k$ }, \\
    \bar{1}  & \text{otherwise}.
\end{cases}
\end{equation*}

Let $\psi_1:F'\ra{\cal B}_K$ be the $r$-fold cyclic covering map corresponding to $h_1$, where $F'$ is the covering space of $\mathcal{B}_K$. 
$F'$ is a $2$-dimensional orbifold with underlying surface $F'_{0}$, which is a closed surface with genus $g=(n-k-2)(r-1)/2$.
When $1\leq j\leq k$, $\psi^{-1}(c_j)$ is a set of $r$ cone points of order $p$, since $h_1(x_j)=\bar{0}$.
Let $\psi^{-1}(c_j)=\{c'_{1,j},c'_{2,j},\ldots, c'_{r,j}\}$, $1\leq j\leq k$.
When $k+1\leq j\leq n$, $h_1(x_j)=\bar{1}$, and $\bar{1}$ has order $r$ in $\z/r$, so $\psi^{-1}_1(c_j)$ is a cone point of order $p$. Let $\psi^{-1}(c_j)=c'_j$, $k+1\leq j\leq n$. 
In summary, $F'$ has $kr+(n-k)$ cone points of order $p$. 

$\psi^{-1}(K^*)$ is a set of $r$ geodesics in $F'$.
Let $\psi^{-1}(K^*)= \{L^*_1,\ldots, L^*_r\}$. 
Each $L^*_i$ goes through $n$ cone points $c'_{i,1},\ldots, c'_{i,k}, c'_{k+1},\ldots,c'_{n}$ successively, $1\leq i\leq r$. 
Let $\t_1$ be the deck transformation of $\psi_1$ corresponding to $\bar{1}\in\z/r$. 
$Fix(\t_1)=\{c'_{k+1},\ldots,c'_n\}$.
Assume that $L^*_ {i+1}=\t_1(L^*_i)$, so $c'_{i+1,j}=\t_1(c'_{i,j})$, $1\leq i< r,1\leq j\leq k$. 
Orient $L^*_1$ and give $L^*_{i}$ the induced orientation, $1< i\leq r$. 
Fix an orientation on $F'$ such that $\t_1$ is a counterclockwise rotation near $c'_l$ by the angle of $2\pi/rp$, $k+1\leq l\leq n$. 
By this construction, $L^*_{i_1}$ only intersects $L^*_{i_2}$ at $c'_l$, and the angle from $L^*_{i_1}$ to $L^*_{i_2}$ is $(i_2-i_1)2\pi/pr$, (we alway assume that counterclockwise is the positive direction.) $k+1\leq l\leq n, 1\leq i_1, i_2\leq r$.

\begin{figure}
\begin{center}
\includegraphics{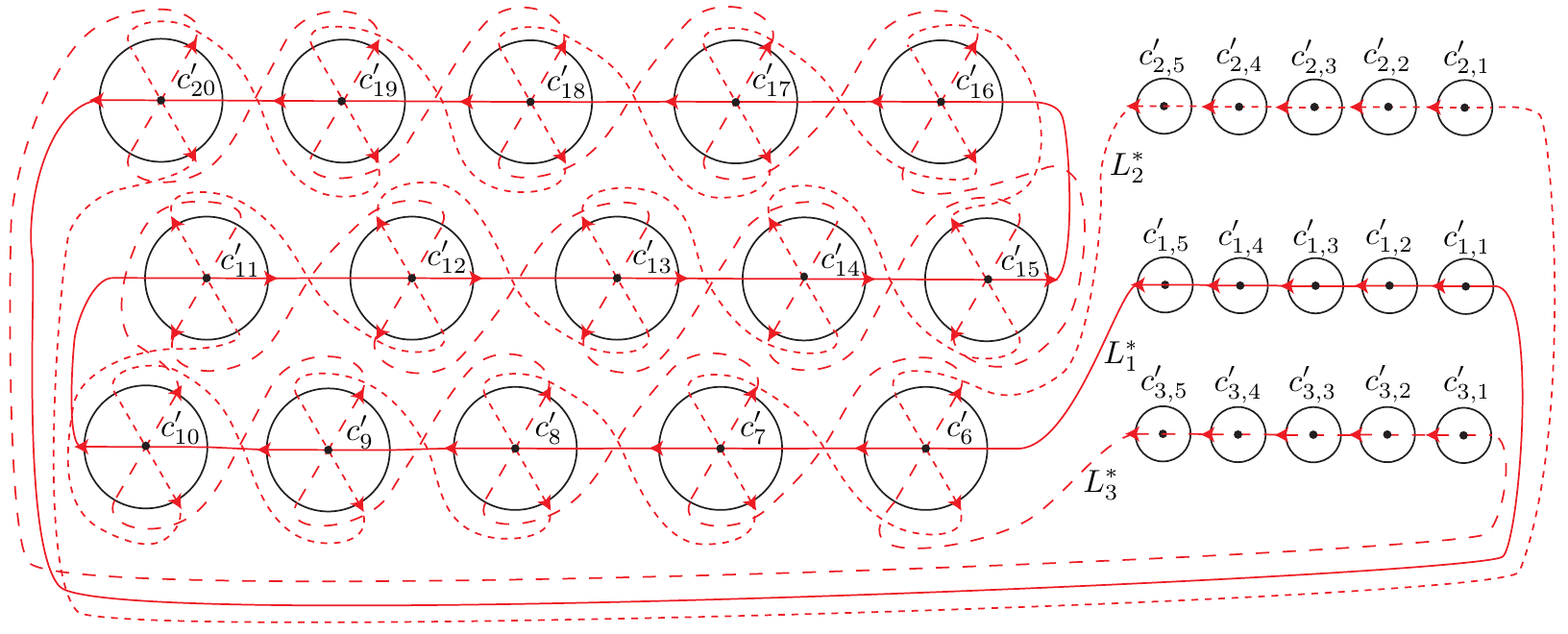}
\end{center}
\caption{\label{f1} $\psi^{-1}(K^*)$}
\end{figure}

\begin{figure}
\begin{center}
\includegraphics[width=14cm]{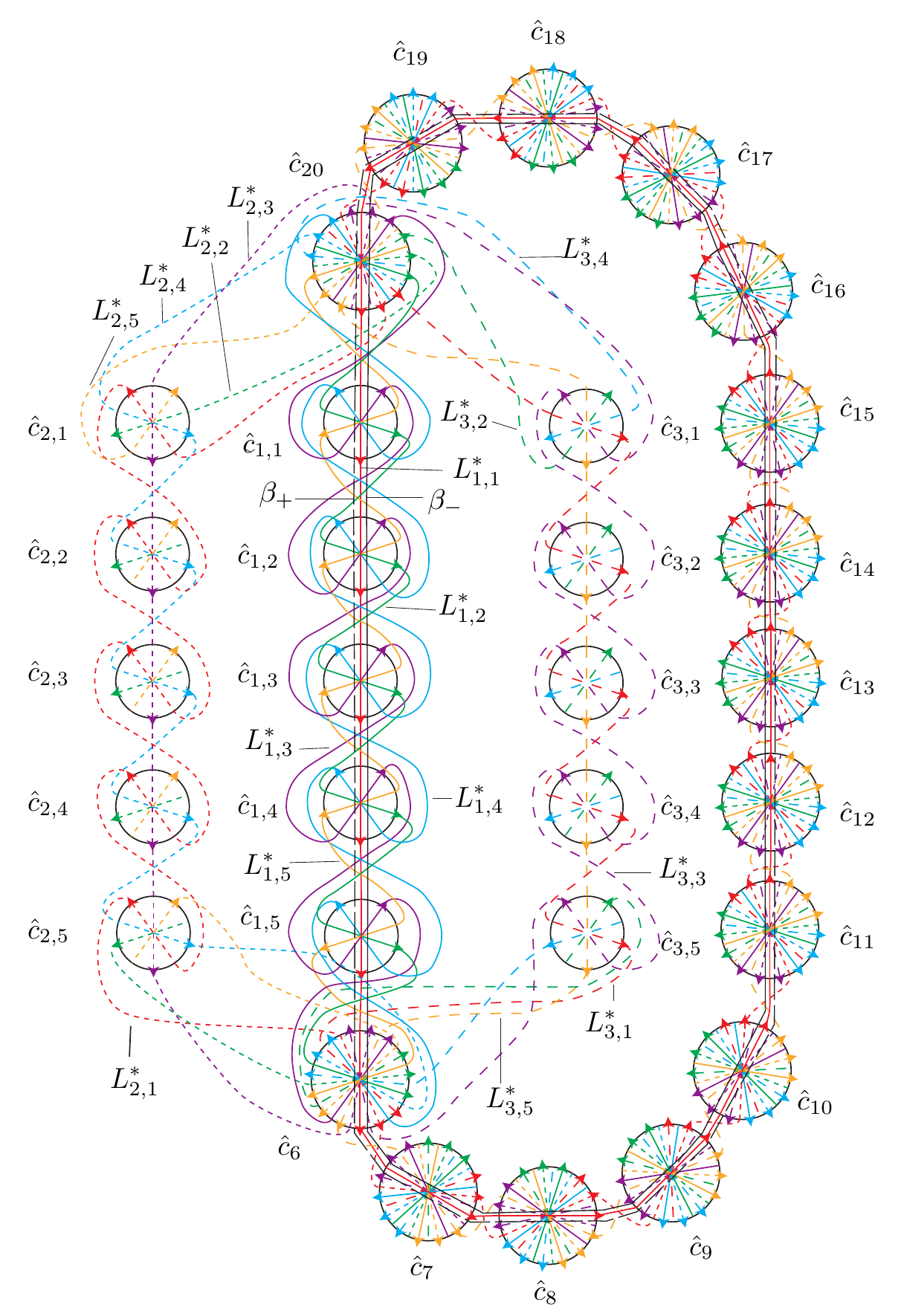}
\end{center}
\caption{\label{df} $\{L_{i,j}^*: 1\leq i\leq r, 1\leq j\leq p\}$}
\end{figure}

Figure  \ref{f1}
shows $\psi^{-1}(K^*)$, when $K=K(\underbrace{2/5,\ldots,2/5}_{5},\underbrace{1/15,\ldots,1/15}_{15})$.
Note that there are some intersections of $L^*_i$'s in Figure \ref{f1} which are not in  $\{ c'_l, k+1\leq l\leq n\}$.
If we draw them on $F'$ which is a surface with genus, those intersections will not appear.
Figure \ref{f1} is a schematic picture, so are Figure \ref{df}, \ref{dfc}, \ref{dfcn}, \ref{25,5f}, \ref{Listar}, \ref{Lijstar6}, \ref{Lijstar6c}, \ref{Lijstar8l} and \ref{610}.

Now we want to define $\psi_2$.
Let $\Gamma_2$ be the fundamental group of $F'$. 
 Then we have:
\begin{equation*}
\begin{split}
\Gamma_2&=\{a_1,b_1,\ldots, a_g,b_g, y_{1,1},y_{1,2},\ldots, y_{1,k},\ldots, y_{r,1}, y_{r,2},\ldots,y_{r,k}, y_{k+1},\ldots,y_{n}: y_{i,j}^{p}=1,\\
&y_{k+1}^{p}=\cdots=y_{n}^{p}=1, \prod_{i=1}^{g}[a_i,b_i]\cdot\prod_{1\leq i\leq r, 1\leq j\leq k}y_{i,j}\cdot y_{k+1}\cdots y_{n}=1,1\leq i\leq r, 1\leq j\leq k\}.
\end{split}
\end{equation*} 
 $y_{i,j}$  and $y_{l}$ are represented by a small circle on $F'$ centered at $c'_{i,j}$ and $c'_l$ respectively, $1\leq i\leq r, 1\leq j\leq k$, $k+1\leq l\leq n$. 
$a_i, b_i$ are the generators of $\pi_1(F'_0)$, $1\leq i\leq g$.

Because $k$ and $n-k$ are both multiple of $p$, there is a homomorphism $h_2: \Gamma_2\rightarrow \z/p$, where
\begin{equation*}
\begin{split}
h_2(a_1)&=\cdots=h_2(a_g)=h_2(b_1)=\cdots=h_2(b_g)=\bar{0},\\ h_2(y_{i,j})&=h_2(y_{k+1})=\cdots=h_2(y_n)=\bar{1}, 1\leq i\leq r, 1\leq j\leq k.
\end{split}
\end{equation*}
Let $\psi_2$ be the covering of $F'$ corresponding to $h_2$, and $F$ the covering space. 
Since the orders of $h_2(y_l)$ ($h_2(y_{i,j})$) equals the order of $c'_l$ ($ c'_{i,j}$) are also $p$, $F$ is a smooth closed orientable surface without cone points, $k+1\leq l\leq n, 1\leq i\leq r,1\leq j\leq k$. 
Denote $\hat{c}_{i,j}=\psi_2^{-1}(c'_{i,j})$ and $\hat{c}_l= \psi_2^{-1}(c'_{l})$, $1\leq i\leq k,1\leq j\leq r, k+1\leq l\leq n$. 

The preimage of $L^*_i$ is a set of $p$ geodesics $\{L^*_{i,1},\cdots,L^*_{i,p}\}$, $1\leq i\leq r$.
Each $L^*_{i,j}$ goes through $n$ points $\hat{c}_{i,1},\ldots,\hat{c}_{i,k}, \hat{c}_{k+1},\ldots, \hat{c}_n$, $1\leq i\leq r,1\leq j\leq p$. 
Let $\t_2$ be the deck transformation of $\psi_2$ corresponding to $\bar{1}\in\z/p$. 
$Fix(\t_2)=\{\hat{c}_{i,1},\ldots,\hat{c}_{i,k}, \hat{c}_{k+1},\ldots, \hat{c}_n: 1\leq i\leq r\}$.
Assume that $L^*_{i,j+1}=\t_2(L^*_{i,j})$, $1\leq i\leq r, 1\leq j< p$. 
Orient $L^*_{i,1}$, such that it passes through $\hat{c}_{i,1},\ldots,\hat{c}_{i,k}, \hat{c}_{k+1},\ldots, \hat{c}_n$  successively, and give $L_{i,j}^*$ the induced orientation, $1\leq i\leq r,1< j\leq p$.  
$F$ admits an orientation such that $\t_2$ is a counterclockwise rotation near the fixed points $\hat{c}_{i,j}$ and $\hat{c}_l$ by $2\pi/p$, $1\leq i\leq r, 1\leq j\leq k, k+1\leq l\leq n$. 
$\{L^*_{i,j},1\leq i\leq r, 1\leq j\leq p\}$ only intersect at those fixed points of $\t_2$, see Figure \ref{df}, for $K=K(\underbrace{2/5,\ldots,2/5}_{5},\underbrace{1/15,\ldots,1/15}_{15})$.
 The following remark gives the details.
 
\begin{remark}\label{angle}
$L^*_{i,j_1}$ intersects $L^*_{i,j_2}$ at $\hat{c}_{i,1},\ldots,\hat{c}_{i,k}, \hat{c}_{k+1},\ldots, \hat{c}_n$ and the angle from $L^*_{i,j_1}$ to $L^*_{i,j_2}$ is $(j_2-j_1)2\pi/p$, $1\leq i\leq r, 1\leq j_1,j_2\leq p$. 
When $i_1\neq i_2$, $L^*_{i_1,j_1}$ intersects $L^*_{i_2,j_2}$ at $\hat{c}_l$, and the angle from $L^*_{i_1,j_1}$ to $L^*_{i_2,j_2}$ is $(i_2-i_1)2\pi/pr+(j_2-j_1)2\pi/p$, $1\leq i_1,i_2\leq r; 1\leq j_1,j_2\leq p, k+1\leq l\leq n$. 
\end{remark}

Let $\psi=\psi_2\circ\psi_1:F\overset{\psi_2} \rightarrow F'\overset{\psi_1}\rightarrow \mathcal{B}_K$. $\psi$ is a $pr$-fold covering map of $\mathcal{B}_K$. Then we have
\begin{gather*}\psi^{-1}(c_j)=
\begin{cases}
\{\hat{c}_{i,j};1\leq i\leq r\}     & \text{if } 1\leq j\leq k, \\
 \hat{c_j}                               & \text{otherwise}.
\end{cases}\\
L^*=\psi^{-1}(K^*)=\{L^*_{i,j}: 1\leq i\leq r; 1\leq j\leq p\}.
\end{gather*}

Next we construct a covering space of $W_K$ from $\psi$ as in Sec. 2 of \cite{abz}.
$W_K$ is a Seifert fibered space with basis $\mathcal{B}_K$. 
There is a covering of $W_K$, say $\Psi$, induced by the Seifert quotient map $f:W_K\rightarrow \mathcal{B}_K$ and associated to $\psi$. 
$Y$ is the covering space of $W_K$ corresponding to $\Psi$. 
$Y$ has a locally-trivial circle bundle Seifert structure since his basis $F$ is a surface. 
Let $\hat{f}:Y\rightarrow F$ be the Seifert quotient map.
$Y$ inherits the $\widetilde{SL_{2}}$ geometry structure from $W_K$. 
Recall $L=\Psi^{-1}(\widetilde{K})$.  
We have the commutative diagram which is analogous to Diagram (5) in Sec. 2 of \cite{abz}.

\begin{figure}
\begin{center}
\includegraphics{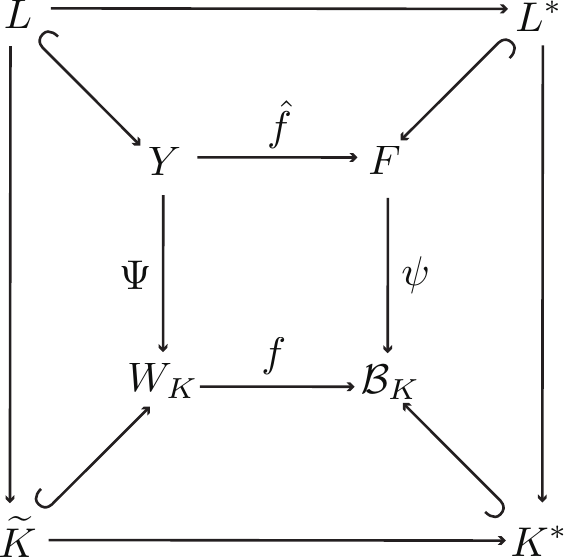}
\end{center}
\caption{\label{diagram} Construction of $Y$}
\end{figure}

From the diagram, we can see that $L=\hat{f}^{-1}(L^*)$. $L$ has exactly $pr$ components by similar discussion as in Sec. 6.1 of \cite{abz}. Let $L=\{L_{i,j}:1\leq i\leq r, 1\leq j\leq p\}$, where $\hat{f}(L_{i,j})=L_{i,j}^*$. 
$\hat{f}|:L_{i,j}\rightarrow L_{i,j}^*$ is a two fold cyclic cover, because $p$ and $pr$ are both odd.

The following proposition can be proved by an analogous  proof of Proposition 6.1 in \cite{abz}.

\begin{prop}\label{semibundle}
The exterior of $L_{1,1}$ in $Y$ is a surface semi-bundle.
\end{prop}

Here we do not repeat the proof. 
We introduce some notations and notes which we need to use later.

Let $T=\hat{f}^{-1}(L^*_{1,1})$, which is a vertical torus overlying $L^*_{1,1}$. 
Because $L_{1,1}^*$ is a geodesic on $F$, $T$ is a totally geodesic torus. $T$ inherits a Euclidean structure from the $\widetilde{SL_{2}}$ structure of $Y$. 
If $\{a,b\}$ is a set of two gedesics which intersects at one point, and $H_1(T)=<a,b>$, then $T$ can be identified to $S^1\times S^1$ where $S^1\times\{*\}$ and $\{*\}\times S^1$ are two geodesics which isotopic to $a$ and $b$ respectively.

Let $F_2=L^*_{1,1}\times[-\e,\e]$, and $\b_-=L_{1,1}^*\times\{-\e\}, \b_+=L_{1,1}^*\times\{\e\}$, where $\e$ is a small positive real number.
($\b_-$ and $\b_+$ are shown in Figure \ref{df}.)
Let 
\begin{equation*}
F_1=F-\overset{\circ}{F}_2, \ Y_i=\hat{f}^{-1}(F_i),\ T_{\pm}=\hat{f}^{-1}(\b_\pm),\; i=1,2.
\end{equation*}
Note that $F_1$ is connected by the construction of $F$.
$Y_2$ is a regular $\e$-neighborhood of $T$ in $Y$.
Define $Y_0=Y_1\cup_{T_-}Y_2$.
The restriction of the Seifert fibration of $Y$ to each of $Y_0, Y_1$ and $Y_2$ is a trivial circle bundle. 
Give the circle fibers of $Y_0$ a consistent orientation. 
Fix two circle fibers of $Y$, $\phi_-,\phi_+$  on $T_-$ and $T_+$ respectively.

$F_0$ is the surface obtained by cutting $F$ open along $\b_+$. 
$Y_0$ is a trivial circle bundle over $F_0$.
Choose a horizontal section $B_0$ of this structure such that $B_0\cap T$ is a geodesic. 
Define $B_i$ to be a subset of $B_0$ in $Y_i$, $i=1,2$. 
Orient $F_0$ and let $F_1$, $F_2$ and their boundaries have the induced orientation. 
Also equip $B_0$, $B_1$, $B_2$ and their boundaries the induced orientation.

Recall that $Y_2=T\times[-\e,\e]$. $T$ is also fibered by geodesics isotopic to $L_{1,1}$. This gives us a {\it new} fibration  of $Y_2$ with base space $\bar{F}_2$. 
We call the fiber of $Y_2$ in this new fibration structure {\it new fiber}, and denoted $\bar{\phi}$. 
We call the fiber from the original fibration of $Y$ the {\it original fiber}, denoted $\phi$.
 Let $\bar{B}_2$ be one horizontal section of $Y_2 \rightarrow \bar{F}_2$ such that $\bar{B}_2\cap T$ is a geodesic. 
Let $N$ be a small regular neighborhood of $L_{1,1}$ in $Y_2$, which is disjoint from other components of $L$ and consists of new fibers.
 Let $M_2=Y_2-\overset{\circ}{N}$, and $M=Y_1\cup_{T_-} M_2$.
 Then $M$ is the exterior of $L_{1,1}$ in $Y$. 
 $\p M_2=T_{2,-}\cup T_{2,+}\cup T_{2,0}$, where $T_{2,0}=\p N$, $T_{2,-}=T_-, T_{2,+}=T_+$.
 $M$ is a graph manifold with boundary $T_{2,0}$ and characteristic tori $T_-$ and $T_+$. Let $\bar{B}_2^0=\bar{B}_2\cap M_2$, which is an annulus with one puncture. Orient $\bar{B}_2^0$ and give $\p\bar{B}_2^0$ the induced orietation. Equip the new fibers of $M_2$ a fixed orientation.

By Prop.6.1 in \cite{abz}, we can construct   essential horizontal surfaces $H_1\subset Y_1$ and $H_2\subset M_2$. 

Let $e$ be the Euler number of the oriented circle bundle of $Y\rightarrow F$. 
Since $Y$ is a $pr$-fold cover of $W_K$,
\begin{equation*}
    \begin{split}
     e&=pr\cdot e(W_K)=pr(-\sum_{i=1}^n \frac{q_i}{p_i})=pr(-\frac{1}{p}(q_1+\cdots +q_k)-\frac{1}{pr}(q_{k+1}+\cdots+q_{n}))\\
     &=-(r(q_1+\cdots q_k)+q_{k+1}+\cdots+q_{n})
\end{split}
\end{equation*}
by (\ref{ew}).
Since $r$ is odd and $K$ only has one component (c.f. (\ref{cn})), $e$ is odd which is the same as in the proof of Prop. 6.1 of \cite{abz}.
Then by choosing the horizontal sections $B_0$ and $\bar{B}_2$ properly,  $H_1$ and $H_2$ will have the same boundary slope as shown in Table 3 in Sec. 6.1 of \cite{abz}.
 We describe the construction of $H_1$ and $H_2$ in the following for later use. 
(c.f. the discussion after the proof of Prop. 6.1 of \cite{abz}.)

\begin{figure}
\begin{center}
\includegraphics{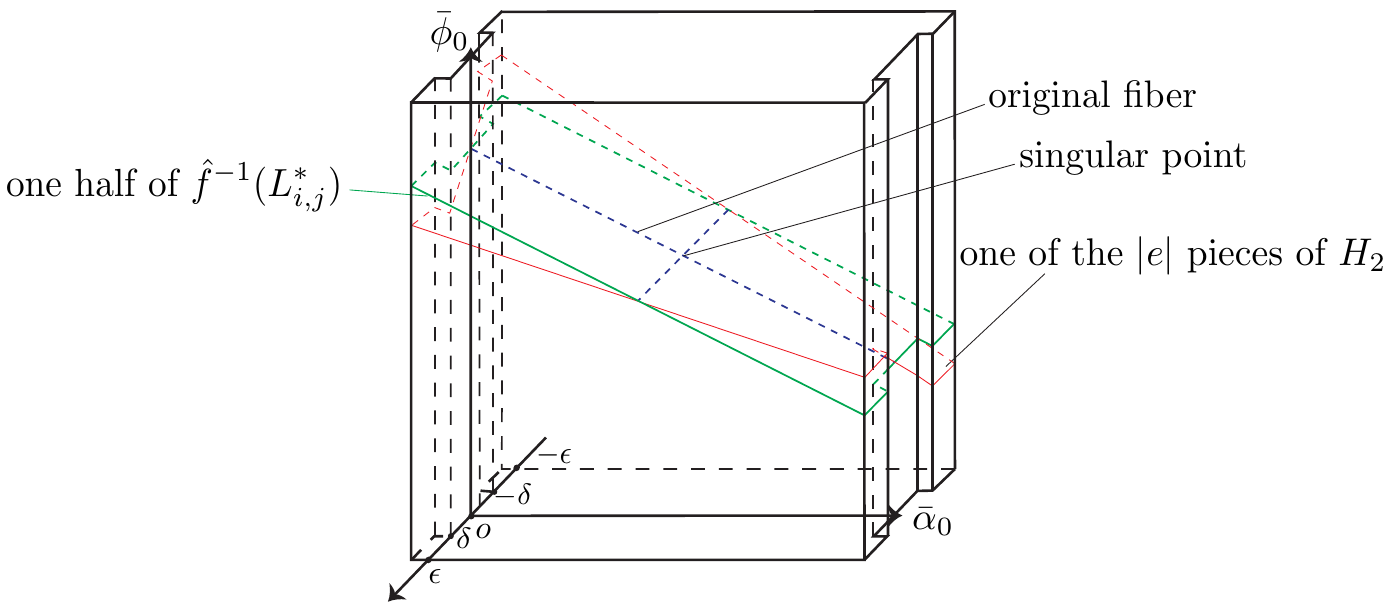}
\end{center}
\caption{\label{fiber} $H_2$}
\end{figure}

\begin{figure}
\begin{center}
\includegraphics[width=15cm]{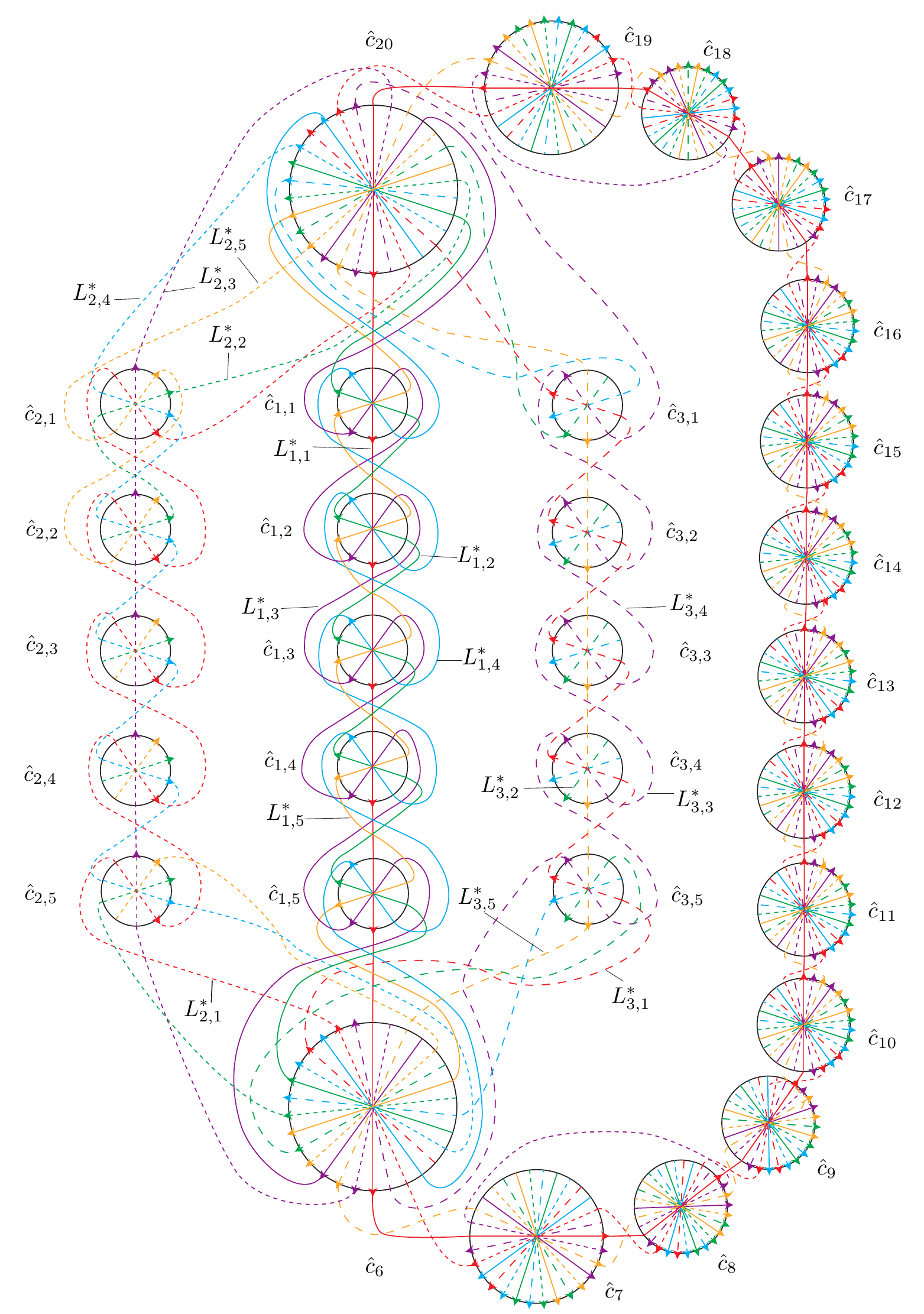}
\end{center}
\caption{\label{dfc} Reorientation of $L^*$}
\end{figure}

\begin{figure}
\begin{center}
\includegraphics[width=15cm]{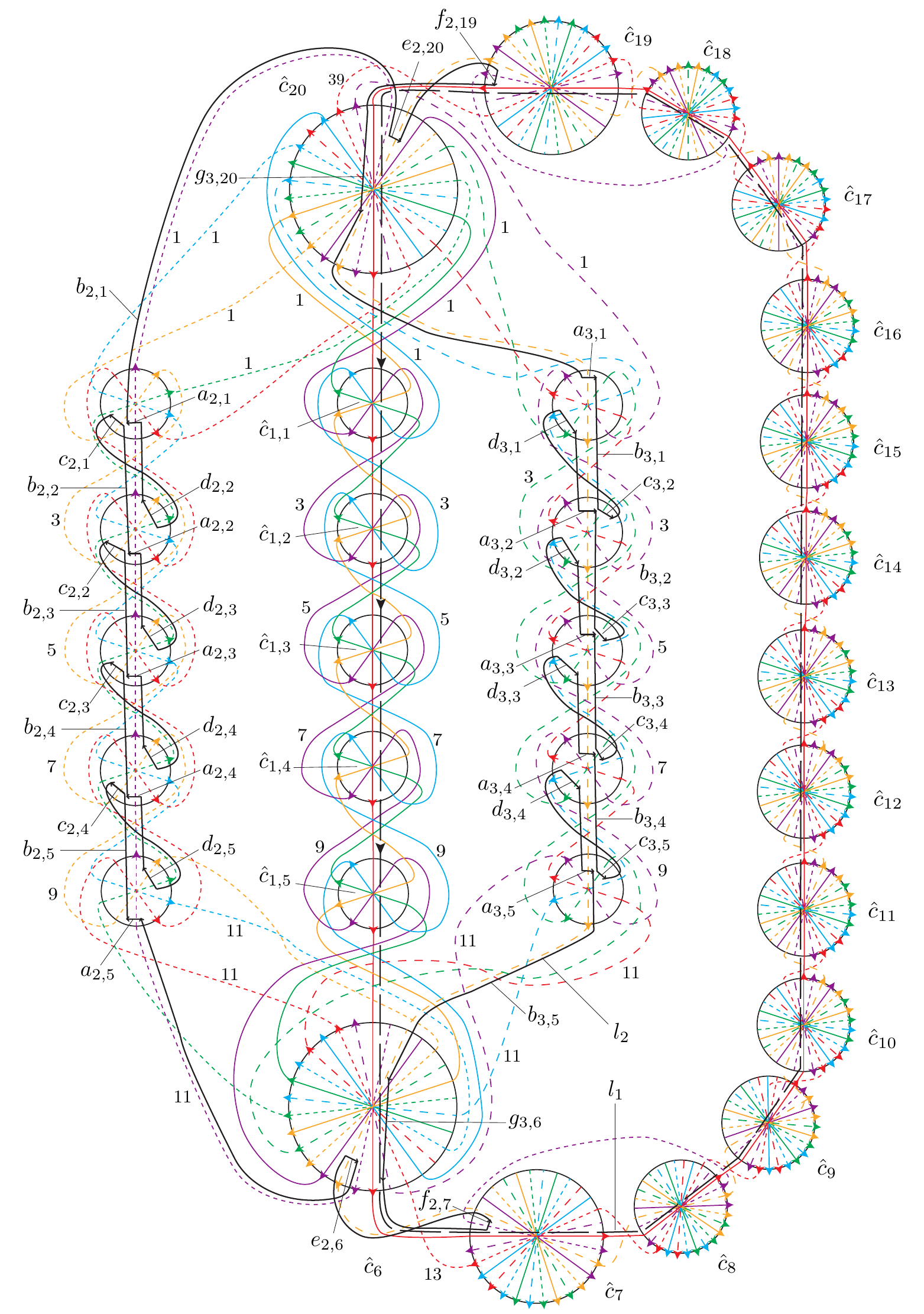}
\end{center}
\caption{\label{dfcn} $\{l_i:1\leq i\leq (r+1)/2\}$, when $p=5, r=3, k=5, n=20$.}
\end{figure}

Set $\bar{\a}_0=\bar{B}_2\cap T$ which is a geodesic. 
Take a fixed new fiber in $T$, denoted $\bar{\phi}_0$.
Then $T=\bar{\a}_0\times\bar{\phi}_0$ and $Y_2=\bar{\a}_0\times\bar{\phi}_0\times[-\e,\e]$. 
$M_2$ can be expressed as following
\begin{equation*}
M_2=(\bar{\a}_0 \times\bar{\phi}_0\times[-\e,\e])-(\overset{\circ}{I}\times \bar{\phi}_0\times(-\d,\d)),
\end{equation*}
where $I=\bar{\a}_0\cap N$, and $\d$ is a positive real number which is less than $\e$. 
Let $\g_1$ be a fixed simple closed geodesic in $T\times\{-\e\}$ of slope $\displaystyle{\frac{1}{2e}-\frac{1}{2}}$ with respect to the basis $\{\bar{\a}_0\times\{-\e\}, \bar{\phi}_0\times\{-\e\}\}$, $\g_2$ a fixed simple closed geodesic in $T\times\{\d\}$ of slope $\displaystyle{-\frac{1}{2e}-\frac{1}{2}}$ with respect to the basis $\{\bar{\a}_0\times\{\d\}, \bar{\phi}_0\times\{\d\}\}$.
As in \cite{abz}, $H_2=\Th_{-}\cup\Th_{0}\cup\Th_{+}$.
$\Th_{-}=\g_1\times[-\e, -\d]$, and $\Th_{+}=\g_2\times[\d, \e]$.
$\Th_0$ is a surface in $(\bar{\a}_0-\overset{\circ}{I})\times\bar{\phi}_0\times[-\d,\d]$ such that $\Th_0\cap (T\times \{t\})$ is a union of $|e|$ geodesic arcs of slope $\displaystyle{-\frac{1}{2}-\frac{t}{2\d e}}$ where $t\in(-\d,\d)$. 
These $|e|$ geodesic arcs are evenly distributed in $(T\times\{t\})\cap M_2$.
$H_2$ is transverse to all the new fibers, and also transverse to the original fibers except when $t=0$. 
Figure \ref{fiber}
illustrates one of the $|e|$ pieces of $H_2$ (c. f. Figure 6 in \cite{abz}). 
The surface fibration ${\cal F}_2$ in $M_2$ is generated by isotoping $H_2$ around the new fibers.

Now see the construction of $H_1$, which is the same as in \cite{abz}.
It's not hard to see that there exists a properly embedded arc, say $\s$, in $B_1$ connecting $T_{-}$ and $T_{+}$.  
Let $\s\times[-1, 1]$ be a regular neighborhood of $\s$ in $B_1$. 
Suppose we pass $(\s\cap T_{-})\times\{-1\}$ to $(\s\cap T_{-})\times\{1\}$ along the orientation of $B_{1}\cap T_{-}$.  
Wrap $\s\times[-1, 1]$ in $B_1$ around the
$\phi$-direction $\displaystyle{(\frac{1-e}{2})}$ times
as we pass from $-1$ to $1$ in $\s\times[-1, 1]\times\phi$.
The resulting surface is $H_1$.  
Let ${\cal F}_1$ be the corresponding surface fibration of $Y_1$
 with $H_1$ as a horizontal surface.

 As in \cite{abz}, we may suppose that $\partial H_1=\partial H_2$. Let $H=H_1\cup H_2$. $H$ is a non-oriented horizontal surface in $M$. 
Then ${\cal F}_1\cup {\cal F}_2$ forms a semi-surface bundle ${\cal F}$ in $M$, as described in Proposition \ref{semibundle}.

Next we want to reorient $\{L_{i,j}:1\leq i\leq r, 1\leq j\leq p, (i,j)\neq(1,1)\}$ such that they travel  from $T_{-}$ to $T_{+}$ in $M_2$ along their new orientation. 
Then we can use Prop. 6.2 of \cite{abz} to show that $\{L_{i,j}:1\leq i\leq r, 1\leq j\leq p, (i,j)\neq(1,1)\}$ are always transverse to ${\cal F}_2$ in $M_2$.

As shown in Figure \ref{df}, $\b_-$ is on the left side of $L_{1,1}^*$ and $\b_+$ is on the right side of $L_{1,1}^*$. 
To get the proper reorientation,  we need that the angle from $L_{1,1}^*$ to $L_{i,j}^*$ (the angle passed when $L_{1,1}^*$ rotates counterclockwise to $L_{i,j}^*$) should between $\pi$ and $2\pi$. 

Define an {\it order} on $\{(i,j):1\leq i\leq r, 1\leq j\leq p\}$.  $(i_1,j_1)\prec(i_2,j_2)$ if $j_1<j_2$ or $j_1=j_2, i_1<i_2$. 
For example $(3,1)\prec(2,3)$.
By Remark \ref{angle}, the angle from $L_{1,1}^*$ to $L_{i_1,j_1}^*$ is less than the angle from $L_{1,1,}^*$ to $L_{i_2,j_2}^*$ at $\hat{c}_l$ if $(i_1, j_1)\prec(i_2,j_2)$, where $1\leq i_1,i_2\leq r, 1\leq j_1,j_2\leq p, k+1<l\leq n$.
At $\hat{c}_{1,s}$, only $L_{1,j}^*$ intersects $L_{1,1}^*$, $1\leq s\leq k$ and $1< j\leq p$.
The angle from $L_{1,1}^*$ to $L_{1,j_1}^*$ is less than the angle from $L_{1,1,}^*$ to $L_{1,j_2}^*$ at $\hat{c}_{1,s}$ if $(1, j_1)\prec(1,j_2)$ i.e. $j_1<j_2$, where $ 1< j_1,j_2\leq p$.

From Remark \ref{angle}, we need to change the direction of $\{L_{i,j}^*: (1,1)\prec(i,j)\preceq ((r+1)/2, (p+1)/2)\}$.
We also change the orientation of the corresponding $L_{i,j}$, $1\leq i\leq r, 1\leq j\leq p$. Then we have the following remark.

\begin{remark}\label{reorientation}
After the reorientation, the angle from $L_{1,1}^*$ to $L_{i,j}^*$ belongs to $(\pi, 2\pi)$.
In addition, the angle from $L^*_{i, (p+1)/2}$ to $L^*_{i,s}$ at $\hat{c}_{i,l}$ belongs to $(\pi, 2\pi)$ when $1<i\leq (r+1)/2$, and the angle from $L^*_{i,p}$ to $L^*_{i,s}$ at $\hat{c}_{i,l}$ belongs to $(\pi, 2\pi)$ when $(r+1)/2<i\leq r$, $1\leq s\leq p$, $1\leq l\leq k$. 
(c.f.Figure \ref{dfc})
\end{remark}

We need introduce some notations before state Prop.\ref{M2}.

By Remark, \ref{angle} 
\begin{multline}\label{intersection}
L^*_{1,1} \text{ intersects } \{L^*_{1,j}:1<j\leq p\} \text{ at }  n \text{ points } \{\hat{c}_{1,s}, \hat{c}_{t}: 1\leq s\leq k, k< t\leq n\}, \\\text{ and intersects } \{L_{i,j}^*:1<i\leq r, 1\leq j\leq p\} \text{ at } n-k \text{ points, } \{\hat{c}_t, k< t\leq n\}. 
\end{multline}
Recall that $F_2=L^*_{1,1}\times [-\e,\e]$, and $F=F_1\cup F_2$. 
Then $L^*_{1,j}$ is separated into $2n$ arcs by $F_1$ and $F_2$, $1<j\leq p$. 
We denote them successively by $\{(L_{1,j}^l)^*: 1\leq l\leq 2n\}$, where $(L_{1,j}^{2l})^*\subset F_2$ and $(L_{1,j}^{2l-1})^*\subset F_1$, $1<j\leq p, 1\leq l\leq n$. 
In addition, $\hat{c}_{1,l}\in (L_{1,j}^{2l})^*$ when $1\leq l\leq k$, and $\hat{c}_{l}\in (L_{1,j}^{2l})^*$ when $k< l\leq n$. 

Different from $L_{1,j}^*$, $L^*_{i,j}$ is separated into $2(n-k)$ arcs by $ F_1$ and $F_2$, when $1<i\leq r, 1\leq j\leq p$.
Note that $L^*_{1,j}$'s do not intersect each other in $F_1$, but $L_{i,j}^*$'s do, $1<i\leq r, 1\leq j\leq p$.  
We divide $L^*_{i,j}\cap F_1$ into $k+1$ parts by the intersection points $\{\hat{c}_{i,1},\ldots, \hat{c}_{i,k}\}$, $1<i\leq r, 1\leq j\leq p$.
We first define $(L_{i,j}^{2l})^*$ to be the segment of $L_{i,j}^*$ in $F_2$ which contains $\hat{c}_l$, $1<i\leq r, 1\leq j\leq p, k<l\leq n$. Let $(L_{i,j}^{2l+1})^*$ to be the segment of $L_{i,j}^*$ between $(L_{i,j}^{2l})^*$ and $(L_{i,j}^{2l+2})^*$, where  $1<i\leq r, 1\leq j\leq p, k<l< n$. 
Let $(L_{i,j}^{2l+1})^*$ be the segment of $L_{i,j}^*$ between $\hat{c}_{i,l}$ and $\hat{c}_{i,l+1}$, $1\leq l< k$, $(L_{i,j}^{1})^*$ the segment of $L_{i,j}^*$ between $(L_{i,j}^{2n})^*$ and $\hat{c}_{i,1}$, and $(L_{i,j}^{2k+1})^*$ the segment of $L_{i,j}^*$ between $\hat{c}_{i,k}$ and $(L_{i,j}^{2(k+1)})^*$, $ 1<i\leq r, 1\leq j\leq p$. 
Now $L_{i,j}^*$ has $k+1+2(n-k)=2n-k+1$ segment components, $1<i\leq r,1\leq j\leq p$.
In summary,
\begin{equation*}
L_{i,j}^*=
\begin{cases}
  \cup_{l=1}^{2n}(L_{i,j}^l)^*    & \text{ if } i=1, 1<j\leq p\\
  (L_{i,j}^1)^*\cup(L_{i,j}^3)^*\cup\cdots\cup(L_{i,j}^{2k+1})^*\cup(\cup_{l=2k+2}^{2n}(L_{i,j}^l)^*)    &\text{ if } {1<i\leq r,1\leq j\leq p}.
  \end{cases}
 \end{equation*}
 See Figure \ref{dfcn} and \ref{25,5f}. The small number next to the arc is the index $l$.
Additional, $(L_{i,j}^{2l})^*\subset { F}_2$, $(i,j,l)\in A$, where 
\begin{equation*}
A=\{i,j,l:  i=1, 1< j\leq p, 1\leq l\leq n, \text{ or } 1< i \leq r , 1\leq j\leq p, k+1\leq l\leq n\}.
\end{equation*}
\begin{figure}
\begin{center}
\includegraphics[width=15cm]{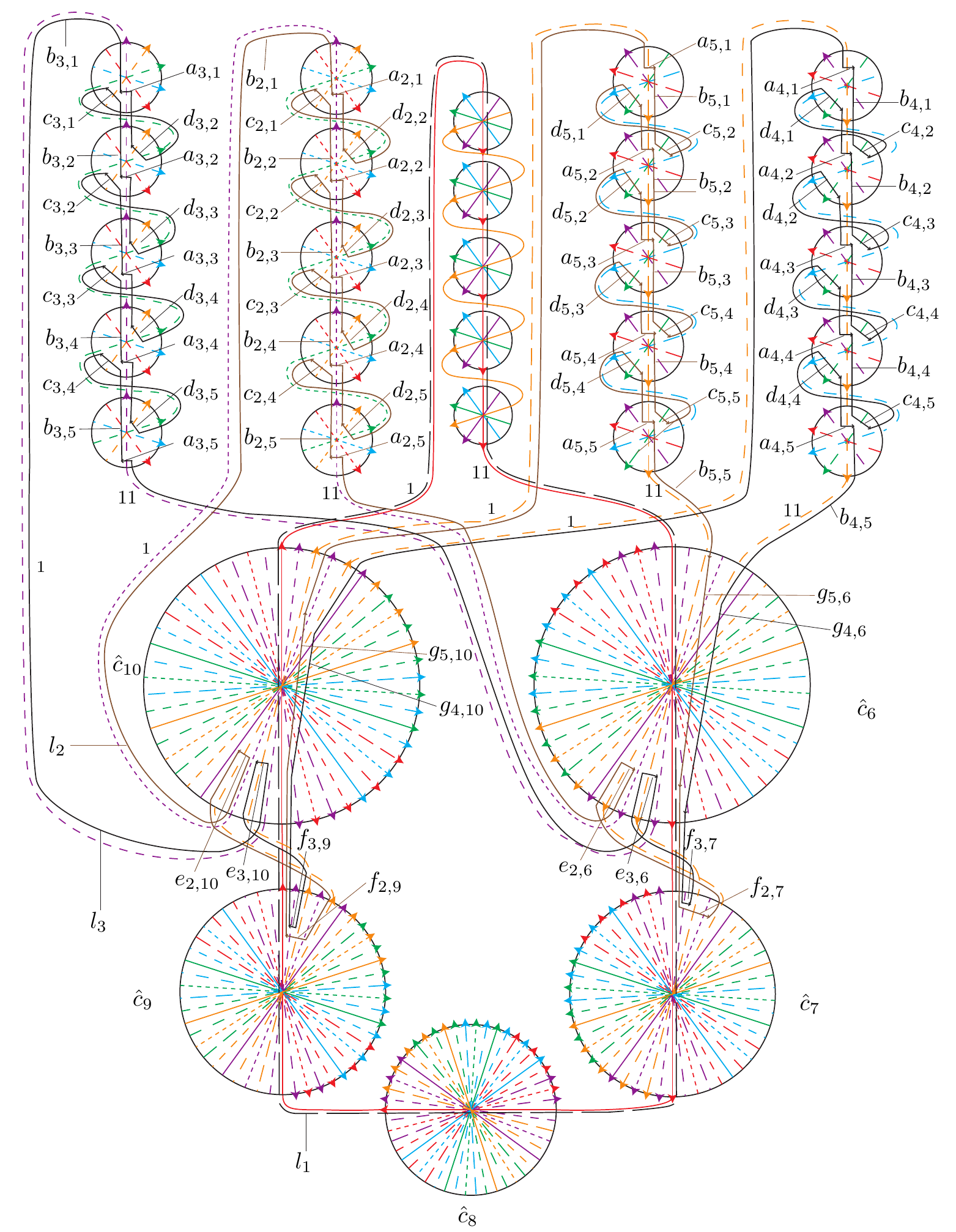}
\end{center}
\caption{\label{25,5f} $\{l_i:1\leq i\leq (r+1)/2\}$, when $p=5, r=5, k=5, n=10$.}
\end{figure}
Let $\phi'_{i,s}=\hat{f}^{-1}(\hat{c}_{i,s})$, the original fiber through $\hat{c}_{i,s}$, and $\phi'_t=\hat{f}^{-1}(\hat{c}_t)$, the original fiber through $\hat{c}_{t}$, where $1\leq i\leq r, 1\leq s\leq k, k+1\leq t\leq n$.
Since $\hat{f}|\  L_{i,j}\rightarrow L_{i,j}^*$
is a 2-fold cyclic cover, $L_{i,j}$ can be separated as the following.

$L_{1,j}$ is separated into $4n$ segments by $T_-$ and $T_+$, $1<j\leq p$.
$L_{i,j}$ is separated into $2(2n-k+1)$ segments by $T_-$, $T_+$ and $\phi'_{i,s}$, $1<i\leq r, 1\leq j\leq p, 1\leq s\leq k$.
We may suppose that
\begin{equation*}
 \hat{f}(L_{i,j}^{l}\cup L_{i,j}^{l+2n})=(L_{i,j}^l)^*,(i,j,l)\in B
 \end{equation*}
where 
\begin{align*}
B=\{i,j,l: 1<j\leq p,1\leq l\leq 2n \text{ when } i=1, \\
\text{ or } 1<i\leq r, 1\leq j\leq p, 1\leq l\leq 2n \text{ and } l\neq 2s, 1\leq s\leq k\}.
\end{align*} 
Then 
\begin{equation*}
L_{i,j}=
\begin{cases}
  \cup_{l=1}^{4n}L_{i,j}^l   & \text{ if } i=1, 1<j\leq p \\
 \underset{s=0,2n}{\cup}( \cup_{l=0}^{k}L_{i,j}^{2l+1+s}\cup(\cup_{l=2k+2}^{2n}L_{i,j}^{l+s}))    & \text{ if } 1<i\leq r, 1\leq j\leq p.
  \end{cases}
\end{equation*}
For each $(L_{i,j}^l)^*$, there exists a vertical annulus 
\begin{equation*}
U_{i,j}^l=\hat{f}^{-1}((L_{i,j}^l)^*), (i,j,l)\in B.
\end{equation*}
Then we have that $ L_{i,j}^{l}\cup L_{i,j}^{l+2n}\subset U_{i,j}^l$, $(i,j,l)\in B$. Suppose 
\begin{equation*}
U_{i,j}^{l,0}=U_{i,j}^l \cap M, (i,j,l)\in B.
\end{equation*}
Note that $ \cup_{l=1}^{2n}U_{1,j}^{l,0}=\hat{f}^{-1}(L^*_{1,j})\cap M$ is a $2n$-punctured torus and $\cup_{l=0}^{k}U_{i,j}^{2l+1,0}\cup(\cup_{l=2k+2}^{2n}U_{i,j}^{l,0})=\hat{f}^{-1}(L^*_{i,j})\cap M$ is a $2(n-k)$-punctured torus, by (\ref{intersection}), $1<i\leq r, 1\leq j\leq p$.
By the construction of $U_{i,j}^{l,0}$, we have 
\begin{equation*}
U_{i,j}^{l,0}\subset
\begin{cases}
M_2      & \text{ if $l$ is even} , \\
 Y_1     & \text{ if $l$ is odd},
\end{cases}
(i,j,l)\in B.
\end{equation*}
\begin{equation}\label{U}
\begin{split}
U_{1,j_1}^{2l,0}\cap U_{1,j_2}^{2l,0}=
\phi'_{1,l}\cap M, &\ \   1<j_1,j_2\leq p, j_1\neq j_2, 1\leq l\leq k;\\
U_{i_1,j_1}^{2l,0}\cap U_{i_2,j_2}^{2l,0}=\phi'_{l}\cap M,  &\ \ 1\leq i_1,i_2\leq r, 1\leq j_1,j_2\leq p, (i_1,j_1),(i_2,j_2)\neq(1,1), \\&\ \ (i_1,j_1)\neq(i_2,j_2), k<l\leq n.
\end{split}
\end{equation}

Analogous to Prop. 6.2 in \cite{abz}, we have the following proposition.

\begin{prop}\label{M2}
Give ${\cal F}_2$ a fixed transverse orientation.
We can isotope $L_{i, j}$ along the original fibers in $U_{i,j}^{2l,0}$, such that $L_{i,j}$ travels from the negative to the positive side of ${\cal{F}}_2$'s leaves.
This isotopy fixes outside a small regular neighborhood of $U_{i,j}^{2l,0}$.
$(i,j,l)\in A$.
\end{prop}

\begin{figure}
\begin{center}
\includegraphics{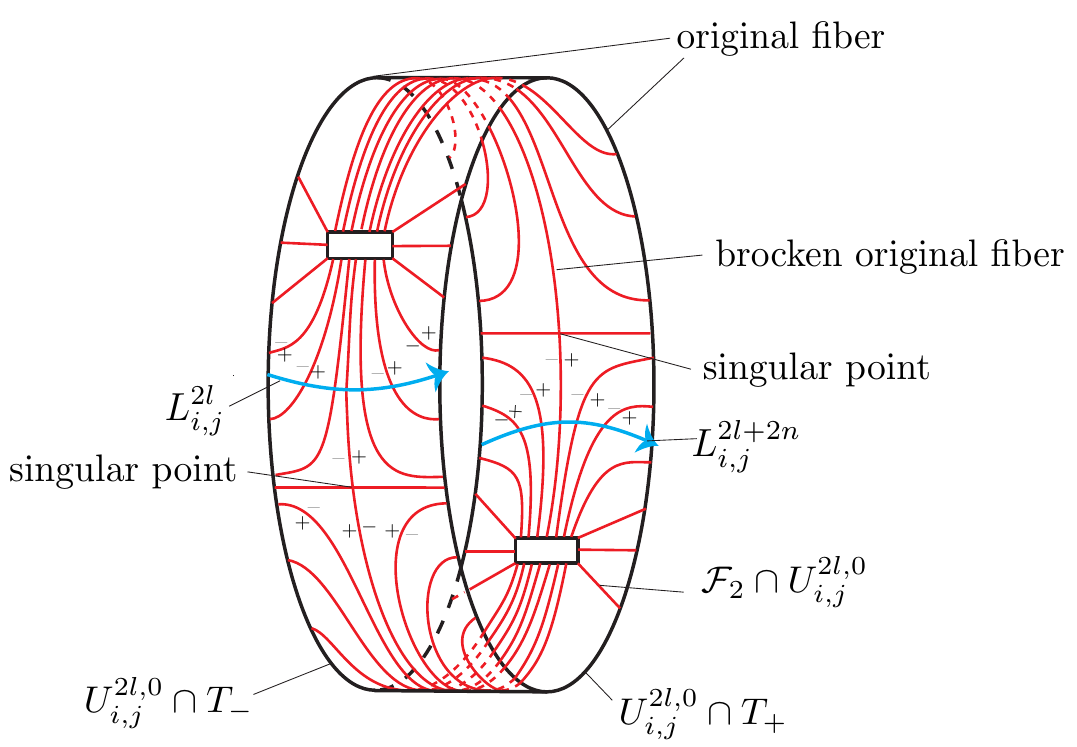}
\end{center}
\caption{\label{foliation} Foliation of $U^{2l,0}_{i,j}$}
\end{figure}

By the construction of $H_2$, there are two
 singular points in the foliation of $U_{i,j}^{2l,0}$\ given by
 $U_{i,j}^{2l,0}\cap{\cal F}_2$, $(i,j,l)\in A$.
 Figure \ref{foliation} (c.f. Figure 7 in Sec. 6.1 of \cite{abz}) illustrates the foliation of $U_{i,j}^{2l,0}$, and the position of $L_{i,j}^{2l}$ and $L_{i,j}^{2l+2n}$ after the isotopy, $(i,j,l)\in A$.
 After the reorientation $L_{i, j}^{2l}$ always travels from $T_{-}$ to $T_{+}$ in $M_2$, so we can isotopy them all above or under the singular points without blocking each other, $(i,j,l)\in A$. 

\begin{figure}
\begin{center}
\includegraphics{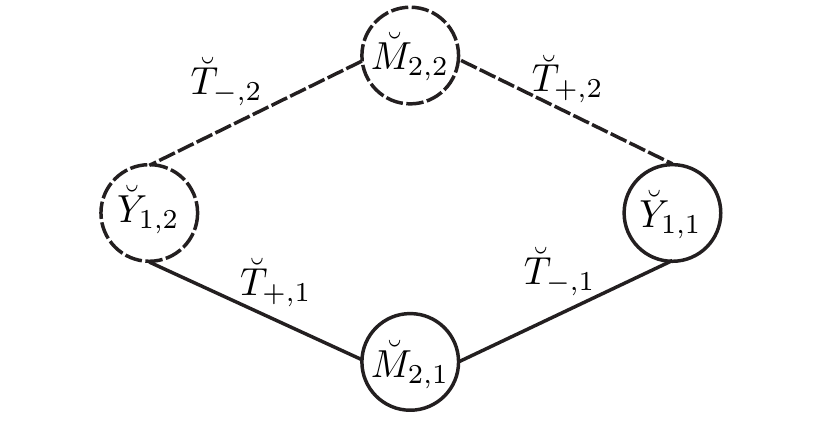}
\end{center}
\caption{\label{double} The graph decomposition of  $\breve{M}$}
\end{figure}

Next we construct a double cover of $M$, say $\breve{M}$. 
Denote $p_2$ to be the covering map. 
The construction of $\breve{M}$ is shown in Figure \ref{double} (c.f. Figure 8 in Sec. 6.1 of \cite{abz}), where
\begin{align*}
\breve{Y}_{1, 1}\cup \breve{Y}_{1, 2}&=p_2^{-1}(Y_1),&
\breve{M}_{2, 1}\cup \breve{M}_{2, 2}&=p_2^{-1}(M_2)\\
\breve{T}_{\pm, 1}\cup \breve{T}_{\pm, 2}&=p_2^{-1}(T_\pm),
\end{align*}
 $\breve{M}$ is also a graph manifold and Figure \ref{double}  shows its JSJ-decomposition.
By construction, $\breve{M}$ is fibered. Denote the surface bundle of $\breve{M}$, $\breve{\cal F}$.
\begin{equation*}
\breve{\cal F}=\breve{\cal F}_{1, 1}\cup \breve{\cal F}_{1, 2},
\end{equation*}
where $\breve{\cal F}_{1, s}=p_2^{-1}({\cal F}_1)\cap \breve{Y}_{1, s}, \breve{\cal F}_{2, s}=p_2^{-1}({\cal F}_2)\cap \breve{M}_{2, s}, s=1, 2$.
Fix a transverse orientation for $\breve{\cal F}$ and let $\breve{\cal F}_{i, s}$ have the induced orientation, $i,s=1, 2$.
$p_2$ can be extended to a double cover of $Y$, denoted $\breve{Y}$. We have
\begin{equation*}
\breve{Y}=\breve{Y}_{1, 1}\cup \breve{Y}_{1, 2}\cup \breve{Y}_{2, 1}\cup\breve{Y}_{2,2}
\end{equation*}
where $\breve{Y}_{2, 1}\cup\breve{Y}_{2,2}=p_2^{-1}(Y_2)$ and $\breve{M}_{2, s}\subset\breve{Y}_{2,s}$, $s=1, 2$.
$\breve{Y}$ has a Seifert fibered structure with fiber $\breve{\phi}=p_2^{-1}(\phi)$. Let $\breve{\phi}$ have the inherited orientation from $\phi$.

By construction, $p_2^{-1}(L_{1,1})$ has two components, denoted $\breve{L}_{1,1,1}$, and $\breve{L}_{1,1,2}$. $\breve{L}_{1,1,s}\subset{\breve{Y}_{2,s}}$, $s=1,2$. Then 
 \begin{equation*}
 \breve{M}=\breve{Y}\setminus (\overset{\circ}{N}(\breve{L}_{1,1,1})\cup \overset{\circ}{N}(\breve{L}_{1,1,2})).
\end{equation*}
$p_2^{-1}(L_{i,j})$ also has two components, $1\leq i\leq r, 1\leq j\leq p, (i,j)\neq(1,1)$. Let $\breve{L}=p_2^{-1}(L)$, then $\breve{L}$ has $2pr$ components.

Proposition \ref{M2} also holds for $\breve{M}_{2,1}$ and $\breve{M}_{2,2}$, so we can isotope $\{\breve{L}_{i,j,s}: 1\leq i\leq r, 1\leq j\leq p, (i,j)\neq(1,1),s=1,2\}$ smoothly along $\breve{\phi}$ such that they always travel from the $``-"$ side of $\breve{\cal F}_{2,s}$ to the $``+"$ side when traveling along their orientations in $\breve{M}_{2,s}$, $s=1,2$.

Next, we perform a sufficiently large number of times of Dehn twist operations on $\{\breve{\cal F}_{1, s} : s=1, 2\}$ along a set of vertical tori $\Gamma$, such that the new surface fiber is transverse to $\{\breve{L}_{i, j, s}, 1\leq i\leq r, 1\leq j\leq p,(i,j)\neq (1,1), s=1,2\}$. 
Since the rational tangle decomposition of $K$, $(\displaystyle{{q_1}/{p},
{q_2}/{p}, \ldots, {q_k}/{p},{q_{k+1}}/{pr},{q_{k+2}}/{pr}, \ldots, {q_{n}}/{pr}} )$, has different denominators, we can't find a set of components of  $L^*$ that is mutually disjoint and goes through all the intersection points. 
$L_{i,j}^*$'s intersect each other in $\breve{M}_{1.s}$, $s=1,2$, $1\leq i\leq r, 1\leq j\leq p, (i,j)\neq (1,1)$. 
The construction of $\Gamma$ is different from \cite{abz} and \cite{gz}.
We need the following discussion before construct $\Gamma$.

Fix $s=1,2$.

Denote $p_2^{-1}(U_{i, j}^{l,0})=\breve{U}_{i, j, 1}^{l, 0}\cup \breve{U}_{i, j, 2}^{l, 0}, (i,j,l)\in B$, where $\breve{U}_{i,j,s}^{2l,0}\subset \breve{M}_{2,s}$ when $(i,j,l)\in A$, and $\breve{U}_{i,j,s}^{2l-1,0}\subset \breve{Y}_{1,s}$ when $(i,j,l)\in C$, where 
\begin{equation*}
C=\{(i,j,l):1\leq i\leq r, 1\leq j\leq p, (i,j)\neq (1,1), 1\leq l\leq n\}.
\end{equation*}
Similarly we denote $\breve{L}_{i,j,1}^{l+t}\cup\breve{L}_{i,j,2}^{l+t}$ to be the lift of $L_{i,j}^{l+t}$, and equip $\breve{L}_{i, j, s}^{l+t}$ the inherited orientation, $(i,j,l)\in B, t=0,2n$. 

$\breve{U}_{i,j,s}^{2l-1,0}$ is a $\breve{\phi}$ vertical annulus in $Y_{1,s}$, $(i,j,l)\in C$. 
$\breve{L}_{i,j,s}^{2l-1}\cup\breve{L}_{i,j,s}^{2l-1+2n}\subset\breve{U}_{i,j,s}^{2l-1,0}$ transverse to $\breve{\phi}$ by Lemma 2.1 in \cite{abz}, $(i,j,l)\in C$.

Different from \cite{abz} and \cite{gz}, the boundaries of some of $\breve{U}_{i,j,s}^{2l-1,0}$'s are not contained in $\breve{T}_{-,s}\cup\breve{T}_{+,s}$, so are the end points of some $\breve{L}_{i,j,s}^{2l-1+t}$, $(i,j,l)\in C, t=0,2n$.
There are two cases. (c.f. Figure \ref{dfcn} and \ref{25,5f}.)

{\bf Case 1}

The head of $\breve{L}_{i,j,s}^{2l-1+t}\in \breve{T}_{-,s}$, where $t=0,2n$, and
\begin{equation} \label{head-1}
\begin{cases}
i=1, 1< j\leq p,  1\leq l\leq n, \text{ or }\\
1<i\leq\frac{r+1}{2},
\begin{cases}
1\leq j\leq \frac{p+1}{2}, l=1 \text{ or } k+2\leq l\leq n, \text{or}\\
\frac{p+1}{2}<j \leq p, k+1\leq l\leq n, \text{ or }
\end{cases}\\
\frac{r+1}{2}<i\leq r,
\begin{cases}
1\leq j< \frac{p+1}{2}, l=1 \text{ or } k+2\leq l\leq n, \text{ or }\\
\frac{p+1}{2}\leq j \leq p, k+1\leq l\leq n.
\end{cases}
 \end{cases}
 \end{equation} 
 The tail of  $\breve{L}_{i,j,s}^{2l-1+t}\in \breve{T}_{+,s}$, where $t=0,2n$, and
\begin{equation}  \label{tail-1}
\begin{cases}
i=1, 1< j\leq p,  1\leq l\leq n, \text{ or }\\
1<i\leq\frac{r+1}{2},
\begin{cases}
1\leq j\leq \frac{p+1}{2}, k+1\leq l\leq n, \text{ or }\\
\frac{p+1}{2}<j \leq p, l=1 \text{ or } k+2\leq l\leq n, \text{ or }
\end{cases}\\
\frac{r+1}{2}<i\leq r,
\begin{cases}
1\leq j< \frac{p+1}{2}, k+1\leq l\leq n, \text{ or }\\
\frac{p+1}{2}\leq j \leq p, l=1 \text{ or } k+2\leq l\leq n.
\end{cases}
 \end{cases}
\end{equation}
In this case, the corresponding boundary of $\breve{U}_{i,j,s}^{2l-1,0}$ is contained in $\breve{T}_{-,s}$ or $\breve{T}_{+,s}$.
Other end points of $\{\breve{L}_{i,j,s}^{2l-1}: 1<i\leq r, 1\leq j\leq p, s=1,2, 1\leq l\leq n\}$ are in $\{\breve{\phi}_{i,j,s}:1<i\leq r, 1\leq j\leq k, s=1,2\}$, where $\breve{\phi}_{i,j,s}$ is the lift of $\phi'_{i,j}$ in $Y_{1,s}$.
Details are given in Case 2.

{\bf Case 2}

The head of $\breve{L}_{i,j,s}^{2l-1+t}, t=0,2n,$ is in
\begin{equation}\label{head-2}
\begin{cases}
\breve{\phi}_{i,l-1,s}&\text{if } 2\leq l\leq k+1, \begin{cases}
1<i\leq \frac{r+1}{2}, 1\leq j\leq \frac{p+1}{2}, \text{or}\\
 \frac{r+1}{2}<i\leq r, 1\leq j< \frac{p+1}{2}, \text{ or }
 \end{cases}\\
 \breve{\phi}_{i,l,s} &\text{if } 1\leq l\leq k, 
 \begin{cases}
1<i\leq \frac{r+1}{2}, \frac{p+1}{2}< j\leq p, \text{or}\\
 \frac{r+1}{2}<i\leq r, \frac{p+1}{2}\leq j< p.
 \end{cases}
 \end{cases}
 \end{equation}
The tail of $\breve{L}_{i,j,s}^{2l-1+t}, t=0,2n,$ is in
\begin{equation} \label{tail-2}
\begin{cases}
\breve{\phi}_{i,l,s} &\text{if } 1\leq l\leq k, 
\begin{cases}
1<i\leq \frac{r+1}{2}, 1\leq j\leq \frac{p+1}{2}, \text{or}\\
 \frac{r+1}{2}<i\leq r, 1\leq j< \frac{p+1}{2}, \text{ or }
 \end{cases}\\
 \breve{\phi}_{i,l-1,s} &\text{if } 2\leq l\leq k+1, \begin{cases}
1<i\leq \frac{r+1}{2}, \frac{p+1}{2}< j\leq p, \text{or}\\
\frac{r+1}{2}<i\leq r, \frac{p+1}{2}\leq j\leq p.
 \end{cases}
 \end{cases}
 \end{equation}
In this case, the corresponding boundary of $\breve{U}_{i,j,s}^{2l-1,0}$ is $\breve{\phi}_{i,l,s}$ or $\breve{\phi}_{i,l-1,s}$.
 
By the construction of ${\cal F}_1$, the surface bundle $\breve{\cal F}_{1,s}$ is transverse to $\breve{\phi}$ in $\breve{Y}_{1,s}$. Then $\breve{\cal F}_{1,s}\cap \breve{U}_{i,j,s}^{2l-1,0}$ gives an interval foliation of $ \breve{U}_{i,j,s}^{2l-1,0}$ which is transverse to the original fiber $\breve{\phi}$. 

In Case 1, we may assume that the oriented arcs $\breve{L}_{i,j,s}^{2l-1}\cup\breve{L}_{i,j,s}^{2l-1+2n}\subset\breve{U}_{i,j,s}^{2l-1,0}$ transverse to the interval foliation $\breve{\cal F}_{1,s}\cap \breve{U}_{i,j,s}^{2l-1,0}$ near $\p\breve{U}_{i,j,s}^{2l-1,0}\subset(\breve{T}_{-,s}\cup\breve{T}_{+,s})$ by Prop \ref{M2}. 
Further more, we may assume that $\breve{L}_{i,j,s}^{2l-1}\cup\breve{L}_{i,j,s}^{2l-1+2n}$ travel from the $``-"$ side to the $``+"$ side of every leaf of $\breve{\cal F}_{1,s}$ near $\p \breve{U}_{i,j,s}^{2l-1,0}$.

\begin{figure}
\begin{center}
\includegraphics{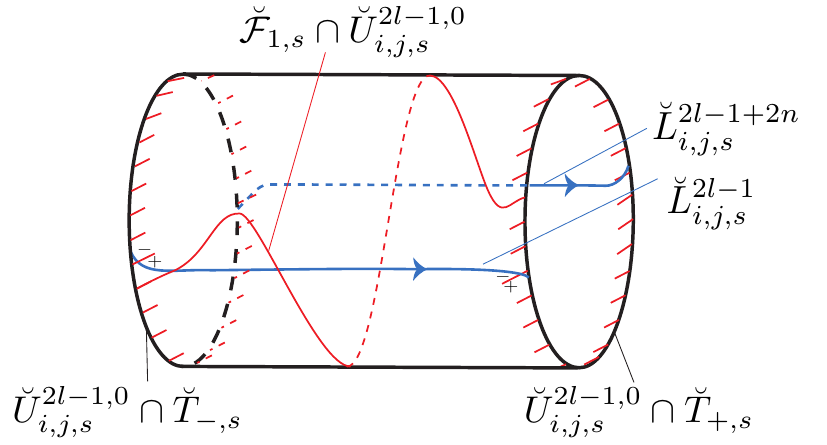}
\end{center}
\caption{\label{odd} $\breve{\cal F}_{1,s}\cap\breve{U}_{i,j,s}^{2l-1,0}$}
\end{figure}

In Case 2, let $N(\breve{\phi}_{i,l,s})$ be a small regular neighborhood of $\breve{\phi}_{i,l,s}$. 
We may isotope $\breve{\cal F}_{1,s}$ along $\breve{\phi}$, such that 
$\breve{L}_{i,j,s}^{2l-1}\cup\breve{L}_{i,j,s}^{2l-1+2n}$ transverse to $\breve{\cal F}_{1,s}\cap \breve{U}_{i,j,s}^{2l-1,0}$ and travel from the $``-"$ to the $``+"$ side of every leaf of $\breve{\cal F}_{1,s}$ near $\p \breve{U}_{i,j,s}^{2l-1,0}$.
Note that this isotopy won't change $\breve{L}\cap \breve{M}_{2,s}$.
Figure \ref{odd} shows $\breve{U}_{i,j,s}^{2l-1}$ in both Case 1 and Case 2.

\begin{figure}
\begin{center}
\includegraphics{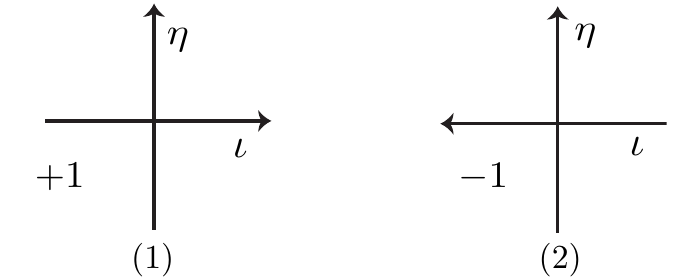}
\end{center}
\caption{\label{sign} Signs of intersection}
\end{figure}

Next, we introduce the algebraic intersection number for two oriented arcs.
Let $\eta$ and $\iota$ be two oriented arcs with only one intersection point.
We say that the algebraic intersection number of $\eta$ and $\iota$, $i(\eta, \iota)= +1$ if the intersection point is shown in Figure \ref{sign}-1, and $i(\eta, \iota)= -1$ if the intersection point is shown in Figure \ref{sign}-2.
 We call $\iota$  intersects $\eta$  \textit{positively} if $i(\eta, \iota)= +1$, otherwise we call  $\iota$ intersects $\eta$ \textit{negatively}. 
 If $\iota$ intersects $\eta$ more than once, the algebraic intersection number of $\eta$ and $\iota$ is the sum of the algebraic intersection numbers for all the intersection points. 

$\G$ is a set of mutually disjoint $(r+1)$ tori.
Let $\Gamma=\{\breve{V}_s^i: 1\leq i\leq (r+1)/2, s=1,2\}$, where $
\breve{V}_1^i\subset\breve{Y}_{1,1}$ and $\breve{V_2^i}\subset\breve{Y}_{1,2}$ are the preimage of a $\phi$-vertical torus, $V^i$ in $M$, $1\leq i \leq (r+1)/2$.
$\hat{f}(V^i)$ is a simple closed curve in $F_1$, denoted $l_i$, $1\leq i \leq (r+1)/2$.
Let ${\cal C}=\{l_i: 1\leq i\leq (r+1)/2\}$ where $l_i=\hat{f}(V^i): 1\leq i\leq (r+1)/2$. 
${\cal C}$ is a set of mutually disjoint $(r+1)/2$ simple closed curves in $F_1$. 
Give every component of ${\cal C}$ an orientation.
We have the following claim.

\begin{claim}\label{torus}
If $i((L_{i,j}^{2l-1})^*, {\cal C})=\sum_{t=1}^{(r+1)/2}i((L_{i,j}^{2l-1})^*, l_t)$ is negative (or positive) for all $(i,j,l)\in C$, then we can perform Dehn twist operations on $\{\breve{\cal F}_{1,s}: s=1,2\}$ along $\G$, such that $\{\breve{L}_{i, j, s}^{2l-1+t}: (i,j,l)\in C, t=0,2n, s=1,2\}$ travels from the $``-"$ to the $``+"$ side of every leaf of the new surface fiber structure of $\breve{Y}_{1,s}$. 
\end{claim}
Proof:
Consider $1\leq t\leq (r+1)/2, (i,j,l)\in C, s=1, 2$ in this proof.
Without losing generality, we assume that $i((L_{i,j}^{2l-1})^*, {\cal C})$ is negative for all $(i,j,l)\in C$.

Let $N(\breve{V}_s^t)$ be a small regular neighborhood of $\breve{V}_s^t$ in $int(\breve{Y}_{1, s})$.
We may suppose that $\{N(\breve{V}_s^t):1\leq t\leq (r+1)/2\}$ are mutually disjoint if we take the regular neighborhoods small enough. 
Define $\p N(\breve{V}_s^t)=\breve{\p}_{1, s}^t\cup \breve{\p}_{2, s}^t$, where $\breve{\p}_{1, s}^t$ and $\breve{\p}_{2, s}^t$ are two $\breve{\phi}$-vertical tori.
$N(\breve{V}_s^t)\cap(\breve{L}_{i, j, s}^{2l-1}\cup\breve{L}_{i, j, s}^{2l-1+2n})$ consists of some arcs if $l_t\cap(L_{i, j, s}^{2l-1})^*\neq\emptyset$.
Each arc in $N(\breve{V}_s^t)\cap\breve{L}_{i, j, s}^{2l-1}$ is corresponding to one point in $l_t\cap(L_{i, j}^{2l-1})^*$.
Note that $\breve{L}_{i, j, s}^{2l-1}$ and $\breve{L}_{i, j, s}^{2l-1+2n}$ travel in the same direction in each $N(\breve{V}_s^t)$ by their orientations.
Now we can assume that the tails of the arcs (with the induced orientation) in $N(\breve{V}_s^t)\cap(\breve{L}_{i, j, s}^{2l-1}\cup\breve{L}_{i, j, s}^{2l-1+2n})$ are contained in $\breve{\p}_{1, s}^t$ if the corresponding intersections  of $l_t$ and $(L_{i,j}^{2l-1})^*$ are negative.

Because $\breve{\cal F}_{1, s}$ is transverse to $\breve{\phi}$, $\breve{\cal F}_{1, s}\cap N(\breve{V}_s^t)$ gives a foliation of $N(\breve{V}_s^t)$ by annuli.

Every $\breve{U}_{i,j,s}^{2l-1,0}$ has a foliation given by $\breve{\cal F}_{1,s}\cap\breve{U}_{i,j,s}^{2l-1,0}$ since $\breve{\cal F}_{1,s}$ is transverse to $\breve{\phi}$. 
One leaf of this foliation is given in Figure \ref{odd}.

\begin{figure}
\begin{center}
\includegraphics{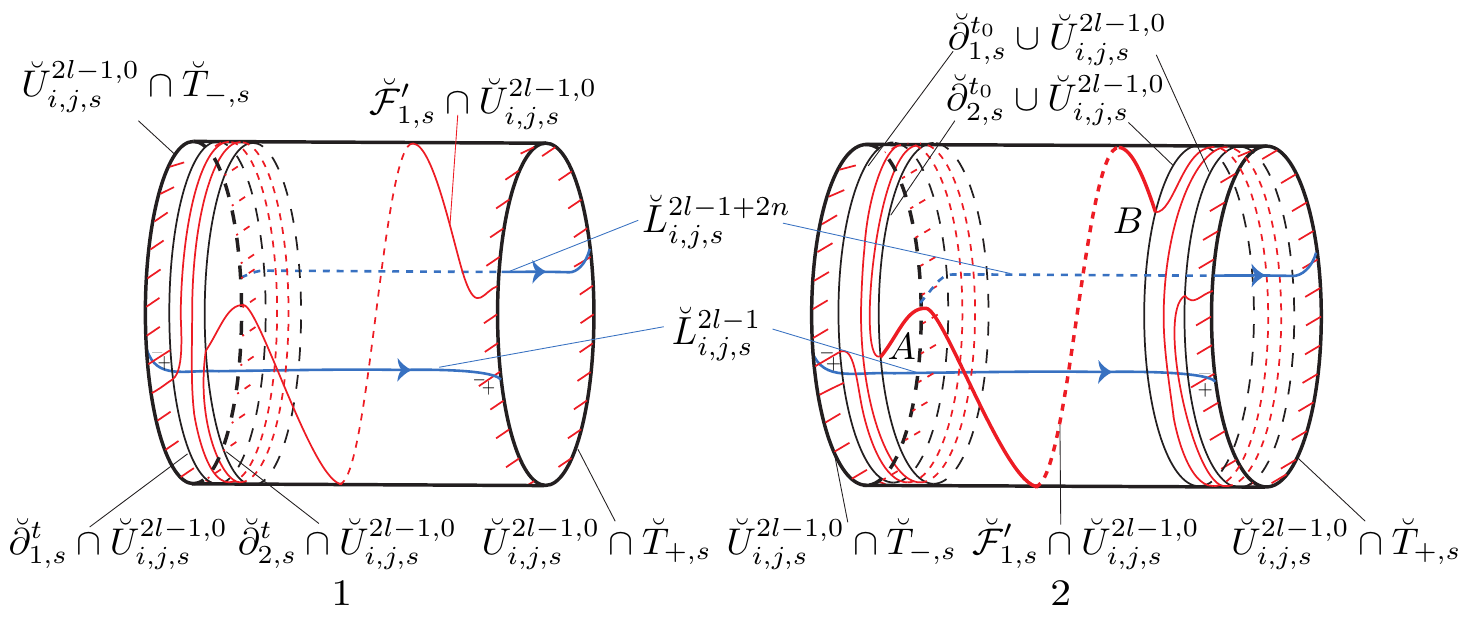}
\end{center}
\caption{\label{dtU} $\breve{\cal F}'_{1,s}\cap\breve{U}_{i,j,s}^{2l-1,0}$.}
\end{figure}

Now we consider two cases.

\textbf{Case A}  The sign of every intersection of $(L_{i,j}^{2l-1})^*$ and ${\cal C}$ is negative.

\begin{figure}
\begin{center}
\includegraphics{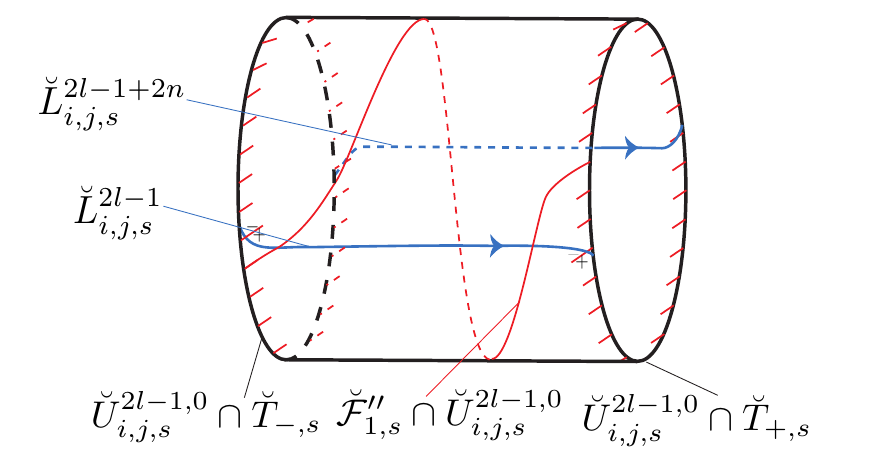}
\end{center}
\caption{\label{dtU4} $\breve{\cal F}''_{1,s}\cap\breve{U}_{i,j,s}^{2l-1,0}$.}
\end{figure}

Let $\breve{\cal F}'_{1,s}$ be the new surface fiber obtained by Dehn twist operations along $\G$ sufficiently large number of times in the direction opposite to the transverse orientation of $\breve{\cal F}_{1, s}$ as we pass from $\breve{\p}_{1, s}^t$ to $\breve{\p}_{2, s}^t$ 
Figure \ref{dtU}-1 shows $\breve{\cal F}'_{1,s}\cap\breve{U}_{i,j,s}^{2l-1,0}$ in this case when $|(L_{i,j}^{2l-1})^*\cap{\cal C}|=1$. 
As in \cite{abz}, we adjust $\breve{\cal F}'_{1, s}$ by isotopy, and 
denote the resulting surface bundle $\breve{\cal F}''_{1, s}$. (c.f. Figure \ref{dtU4}.)
This isotopy only changes $\breve{\cal F}'_{1,s}$ in a small regular neighborhood of $\breve{U}_{i,j,s}^{2l-1,0}$ and fixes $\breve{\cal F}'_{1,s}$ in a small regular neighborhood of $\p\breve{Y}_{1,s}$.
 Then $\breve{L}_{i,j,s}^{2l-1}$ and $\breve{L}_{i,j,s}^{2l-1+2n}$ are travel from the $``-"$ to the $``+"$ side of every leaf of $\breve{\cal F}''_{1,s}$.

If $|(L_{i,j}^{2l-1})^*\cap{\cal C}|>1$, we need wrap $\breve{U}_{i, j, s}^{2l-1,0}\cap\breve{\cal F}_{1,s}$ more times, but in the same direction, since the sign of every intersection of $(L_{i,j}^{2l-1})^*$ and ${\cal C}$ is negative. 
Then we still have that, after isotopy, the new interval foliation in each $\breve{U}_{i, j, s}^{2l-1,0}$ given by $\breve{U}_{i,j,s}^{2l-1,0}\cap \breve{\cal F}''_{1,s}$ becomes transverse to $\breve{L}_{i, j, s}^{2l-1}\cup \breve{L}_{i, j, s}^{2l-1+2n}$. 

\textbf{Case B}  There are two intersection points of $(L_{i,j}^{2l-1})^*$ and ${\cal C}$ with different signs.

Suppose that $(L_{i,j}^{2l-1})^*$ only intersects one component of ${\cal C}$, say $l_{t_0}$, twice in different signs, where $t_0$ is a fixed integer between 1 and $(r+1)/2$. 
$\breve{U}_{i,j,s}^{2l-1,0}\cap\breve{V}_s^{t_0}$ consists of two original circle fibers. 
Each of $\breve{L}_{i,j,s}^{2l-1}\cap N(\breve{V}^{t_0}_s)$ and $\breve{L}_{i,j,s}^{2l-1+2n}\cap N(\breve{V}^{t_0}_s)$ are two arcs.
The tail of one arc of $\breve{L}_{i,j,s}^{2l-1}\cap N(\breve{V}_s^{t_0})$ ($\breve{L}_{i, j,s}^{2l-1+2n}\cap N(\breve{V}_s^{t_0})$) lies in $\breve{\p}_{1,s}^{t_0}$, and the tail of the other one lies in $\breve{\p}_{2,s}^{t_0}$. 
Let $\breve{\cal F}'_{1,s}$ be the surface fiber obtained by Dehn twist operations only along $\breve{V}^{t_0}_s$ in the direction opposite to the transverse orientation of $\breve{\cal F}_{1,s}$ as we pass from $\breve{\p}_{1,s}^{t_0}$ to $\breve{\p}_{2,s}^{t_0}$.
 $\breve{\cal F}'_{1,s}\cap\breve{U}_{i,j,s}^{2l-1,0}$ is as shown in Figure \ref{dtU}-2.

We perform an isotopy of $\breve{\cal F}'_{1,s}$ in a small regular neighborhood of $\breve{U}_{i, j,s}^{2l-1,0}$ in $\breve{Y}_{1,s}$, say $N(\breve{U}_{i, j,s}^{2l-1,0})$, such that $\breve{U}_{i, j,s}^{2l-1,0}\cap\breve{\cal F}'_{1,s}$ changes back to $\breve{U}_{i, j,s}^{2l-1,0}\cap \breve{\cal F}_{1,s}$. 
The isotopy is described as following. 
Push the whole arc $AB$ (as shown in Figure \ref{dtU}-2)
(like a finger move) along $\breve{\phi}$ in $\breve{U}_{i, j,s}^{2l-1,0}$ in the direction opposite to the Dehn twist operation and the same times as the operation, meanwhile fix $\partial N(\breve{U}_{i, j,s}^{2l-1,0}\cap \breve{\cal F }'_{1,s})$ all the time.  
We still call the surface bundle $\breve{\mathcal{F}}_{1, s}^{'}$ after this isotopy.

$(L_{i,j}^{2l-1})^*$ may intersects two different components of ${\cal C}$  in different signs, and there may be more than one pair of different sign intersection points of $(L_{i,j}^{2l-1})^*$ and ${\cal C}$.
By a similar discussion, the new foliation of $\breve{U}_{i,j,s}^{2l-1,0}$ given by the new surface fiberation obtained by the Dehn twist operation corresponding to any pair of different sign intersection points will be the same as $\breve{U}_{i, j,s}^{2l-1,0}\cap \breve{\cal F}_{1,s}$ after isotopy.

Since $i((L_{i,j}^{2l-1})^*,{\cal C})<0$, there are some other intersection points of  $(L_{i,j}^{2l-1})^*$ and ${\cal C}$ with negative signs except for the paired different sign intersection points.
Let $\breve{\cal F}''_{1,s}$ be the surface fiber obtained by Dehn twisting $\breve{\cal F}_{1,s}$ along $\G$ sufficiently large number of times in the direction opposite to the transverse orientation of $\breve{\cal F}_{1, s}$ as we pass from $\breve{\p}_{1, s}^t$ to $\breve{\p}_{2, s}^t$.
By the above discussion and Case A, $\breve{L}_{i,j,s}^{2l-1}$ and $\breve{L}_{i,j,s}^{2l-1+2n}$ are transverse to the new surface fibers, $\breve{\cal F}''_{1,s}$, after isotopy.  $\square$

Now we construct $\{l_i: 1\leq i\leq (r+1)/2\}$.

Let $l_1$ be the simple closed curve in a small regular neighborhood of $\b_-$ in $F_1$, and equip $l_1$ the same orientation as $L_{1,1}^*$. (c.f. Figure \ref{dfcn} and Figure \ref{25,5f}). (The small number next to the segment of $(L_{i,j}^{l})^*$ is the index $l$.) 
We can see that $({L}_{i,j}^{2l-1})^*\cap l_1\neq \emptyset$ if and only if $(i,j,l)$ is as described in (\ref{head-1}), and $i((L_{i,j}^{2l-1})^*, l_1)=-1$ in this case,
by Remark \ref{reorientation}. 

Next, we construct $\{l_i: 2\leq i\leq (r+1)/2\}$.
Let $D_{i,l}$ be a disk in $F$ centered at $\hat{c}_{i,l}$ with radius greater than $\e$, $1<i\leq r, 1\leq l\leq k$. 
Similarly, let $D_l$ be a disk in $F$ centered at $\hat{c}_{l}$ with radius greater than $\e$, $ k+1\leq l\leq n$.
At first, we construct some oriented arcs on $\{D_{i,l}:1\leq i\leq r, 1\leq l\leq k\}\cup\{D_{l}: k+1\leq l\leq n\}$ such that they intersect $L_{i,j}^*$ {\it negatively}, $1\leq i \leq r, 1\leq j\leq p, (i,j)\neq(1,1)$. 
Then we connect them by arcs on $F_1$ and parallel to some component of $L^*$. 
The arcs on $D_{i,l}$ $(D_t)$ are given as the following, $1\leq i\leq r, 1< l\leq k, k+1\leq t\leq n$. (c.f. Figure \ref{dfcn} and  \ref{25,5f})

Let $a_{i,l}$ be an oriented arc on $D_{i,l}\cap F_1$ intersecting the following arcs negatively.
\begin{equation}\label{a}
  \begin{cases}
 (L_{i,(p+1)/2}^{2l+1})^* & 1<i\leq (r+1)/2,\\
 (L_{i,p}^{2l-1})^* & (r+1)/2<i\leq r,
 \end{cases}
 1\leq l\leq k.
 \end{equation}

Let $c_{i,l}$ be an oriented arc on $D_{i,l}\cap F_1$ intersecting the following arcs negatively.
  \begin{equation*}
   \begin{cases}
 (L_{i,p}^{2l-1})^* &1<i\leq (r+1)/2, 1\leq l\leq k-1;\\ 
 (L_{i,(p-1)/2}^{2l+1})^*& (r+1)/2<i\leq r, 2\leq l\leq k.
   \end{cases}
 \end{equation*}

Let $d_{i,l}$ be an oriented arc on $D_{i,l}\cap F_1$ intersecting the following arcs negatively.
 \begin{equation*}
 \begin{cases}
 (L_{i,(p-1)/2}^{2l-1})^* \text{and }(L_{i,p}^{2l+1})^*, &1<i\leq (r+1)/2, 2\leq l\leq k;\\ 
 (L_{i,(p-1)/2}^{2l-1})^* \text{and }(L_{i,p-1}^{2l+1})^*, & (r+1)/2<i\leq r,1\leq l\leq k-1.
  \end{cases}
\end{equation*}

$e_{i,l} \ (l=k+1,n)$ and $f_{i,l}\ (l=k+2,n-1)$ is defined for $1<i\leq (r+1)/2$.

Let $e_{i,l}$ be an oriented arc on $D_{l}\cap F_1$ intersecting the following arcs negatively.
$ \begin{cases}
(L_{i+(r-1)/2,p}^{2k+3})^*     &\text{if } l=k+1, \\
(L_{i+(r-1)/2,p}^{2n-1})^*   &\text{if } l=n,
\end{cases}$
 $1<i\leq (r+1)/2$.
 
Let $f_{i,l}$ be an oriented arc on $D_{l}\cap F_1$ intersecting the following arcs negatively.
 \begin{equation*}
 \begin{split}
 & \begin{cases}
 \{(L_{i,(p+1)/2}^{2k+5})^*:2\leq i\leq(r+1)/2\}\cup\{(L_{i+(r-1)/2,p}^{2k+3})^*: 2<i\leq (r+1)/2\}    & \text{if }l=k+2,\\
  \{(L_{i,(p+1)/2}^{2n-3})^*:2\leq i\leq(r+1)/2\}\cup\{(L_{i+(r-1)/2,p}^{2n-1})^*: 2<i\leq (r+1)/2\}     & \text{if }l=n-1,
\end{cases} \\&1<i\leq (r+1)/2.
\end{split}
\end{equation*}
We arrange $f_{i_1,l}$ closer to $\hat{c}_{l}$ than $f_{i_2,l}$ if $i_1<i_2$. (c.f. Figure \ref{25,5f})

$g_{i,l}\ (l=k+1,n)$ is defined for $(r+1)/2<i\leq r$.
Let $g_{i,l}$ be an oriented arc on $D_{l}\cap{F_1}$ intersecting the following arcs negatively, $(r+1)/2<i\leq r$.
\begin{equation*}
\begin{cases}
(L_{i_1,j_1}^{2k+3})^*&(1,1)\prec(i_1,j_1)\prec((r+1)/2-(r-i),(p+1)/2),\text{ when }l=k+1,\\
(L_{i_1,j_1}^{2k+1})^*&((r+1)/2,(p+1)/2)\prec(i_1,j_1)\prec(i,p), \text{ when }l=k+1,\\
(L_{i_1,j_1}^{2n-1})^*&(1,1)\prec(i_1,j_1)\prec((r+1)/2-(r-i),(p+1)/2),\text{ when }l=n,\\
(L_{i_1,j_1}^{1})^*&((r+1)/2,(p+1)/2)\prec(i_1,j_1)\prec(i,p),\text{ when }l=n.
\end{cases}
\end{equation*}

When $1<i\leq (r+1)/2$, connect the head of $a_{i,1}$ and the tail of $e_{i,n}$ by an arc $b_{i,1}$, and connect the head of $a_{i,l}$ and the tail of $c_{i,l-1}$ by an arc $b_{i,l}$, where $b_{i,1}$ and $b_{i,l}$ are arcs parallel to $L^*_{i, (p+1)/2}$ and in a small regular neighborhood of $L^*_{i,(p+1)/2}$ in $F_1$.
Give $b_{i,l}$ the induce orientation which is the same as $L^*_{i,(p+1)/2}$, $1\leq l\leq k$, $1<i\leq (r+1)/2$.
When $(r+1)/2<i\leq r$, connect the head of $a_{i,l}$ and the tail of $c_{i,l+1}$ by an arc $b_{i,l}$, and connect the head of $a_{i,k}$ and the tail of $g_{i,k+1}$ by an arc $b_{i,k}$,  where $b_{i,l}$ and $b_{i,k}$ are arcs parallel to $L^*_{i,p}$ and in a small regular neighborhood of $L^*_{i,p}$ in $F_1$, $1\leq l\leq k-1$.
Give $b_{i,l}$ the induce orientation which is the same as $L^*_{i,p}$, $1\leq l\leq k$, $(r+1)/2<i\leq r$.
(c.f. Figure \ref{dfcn} and Figure \ref{25,5f}.)

By Remark \ref{reorientation}, $b_{i,l}$ intersects 
\begin{equation}\label{b1}
\begin{cases}
(L_{i,j}^{2l+1})^* \text{ and } i((L_{i,j}^{2l+1})^*,b_{i,l})=-1& 1\leq j< (p-1)/2,\\
(L_{i,j}^{2l-1})^*\text{ and } i((L_{i,j}^{2l-1})^*,b_{i,l})=-1   &  (p+1)/2< j\leq p, 
\end{cases}
\end{equation}
when $1\leq l\leq k, 1<i\leq (r+1)/2$, and
\begin{equation}\label{b2}
\begin{cases}
(L_{i,j}^{2l+1})^* \text{ and } i((L_{i,j}^{2l+1})^*,b_{i,l})=-1& 1\leq j< (p-1)/2,\\
(L_{i,j}^{2l-1})^*\text{ and } i((L_{i,j}^{2l-1})^*,b_{i,l})=-1   &  (p+1)/2\leq j< p, 
\end{cases}
\end{equation}
when $1\leq l\leq k, (r+1)/2<i\leq r$.

Now we connect all these arcs by arcs in a small regular neighborhood of $L^*$ to construct $(r-1)/2$ simple closed curves as following. Note that these connecting arcs are mutually disjoint and do not intersect $L^*$.

When $2\leq i\leq (r+1)/2$, connect the head of $c_{i,l}$ and the tail of $d_{i,l+1}$ by an arc in a small regular neighborhood of $L^*_{i,(p-1)/2}$ and parallel to $L^*_{i,(p-1)/2}$, $1\leq l\leq k-1$;
connect the head of $d_{i,l}$ and the tail of 
$a_{i,l-1}$ by an arc in a small neighborhood of $L^*_{i,(p+1)/2}$ and parallel to $L^*_{i,(p+1)/2}$, $2\leq l\leq k$.

When $(r+1)/2<i\leq p$, connect the head of $c_{i,l}$ and the tail of $d_{i,l-1}$ by an arc in a small regular neighborhood of $L^*_{i,p-1}$ and parallel to $L^*_{i,p-1}$, $2\leq l\leq k$;
connect the head of $d_{i,l}$ and the tail of $a_{i,l+1}$ by an arc in a small regular neighborhood of $L^*_{i,p}$ and parallel to $L^*_{i,p}$, $1\leq l\leq k-1$.

Connect the head of $e_{i,n}$ and the tail of $f_{i,n-1}$ by an arc in a small regular neighborhood of $L^*_{i+(r-1)/2,p}$ and parallel to $L^*_{i+(r-1)/2,p}$, when $1<i\leq (r+1)/2$.

Connect the head of $f_{i,n-1}$ and the tail of $g_{r+2-i,n}$ by an arc in a small regular neighborhood of $L^*_{1,1}$ and parallel to $L^*_{1,1}$, when $1<i\leq (r+1)/2$.

Connect the head of $g_{i,n}$ and the tail of $a_{i,1}$ by an arc in a small regular neighborhood of $L^*_{i,p}$ and parallel to $L^*_{i,p}$, when $(r+1)/2<i\leq r$.

Connect the head of $g_{i,k+1}$ and the tail of $f_{r+2-i,k+2}$ by an arc in a small regular neighborhood of $L^*_{1,1}$ and parallel to $L^*_{1,1}$, when $1<i\leq (r+1)/2$.

Connect the head of $f_{i,k+2}$ and the tail of $e_{i,k+1}$ by an arc in a small regular neighborhood of $L^*_{i+(r-1)/2,p}$ and parallel to $L^*_{i+(r-1)/2,p}$, when $1<i\leq (r+1)/2$.

Connect the head of $e_{i,k+1}$ and the tail of $a_{i,k}$ by an arc in a small regular neighborhood of $L^*_{i,(p+1)/2}$ and parallel to $L^*_{i,(p+1)/2}$, when $1<i\leq (r+1)/2$.

The following diagram shows the construction of $l_i$, $2\leq i\leq (r+1)/2$.

$a_{i,k}\rightarrow b_{i,k}\rightarrow 
c_{i,k-1}\xrightarrow {L_{i,(p-1)/2}^*}d_{i,k}\xrightarrow{L^*_{i,(p+1)/2}} a_{i,k-1}\rightarrow b_{i,k-1}\rightarrow c_{i,k-2}
\xrightarrow{L^*_{i,(p-1)/2}} d_{i,k-1}\cdots a_{i,1}
\rightarrow{b_{i,1}}\rightarrow e_{i,n}
\xrightarrow{L^*_{i+(r-1)/2,p}}f_{i,n-1}\xrightarrow{L^*_{1,1}}
g_{r+2-i,n}\xrightarrow{L_{r+2-i,p}^*}a_{r+2-i,1}
\rightarrow{b_{r+2-i,1}}\rightarrow c_{r+2-i,2}\xrightarrow{L_{r+2-i,p-1}^*} d_{r+2-i,1}\xrightarrow{L^*_{r+2-i,p}}
a_{r+2-i,2}\rightarrow{b_{r+2-i,2}}\rightarrow c_{r+2-i,3}\xrightarrow{L^*_{r+2-i,p-1}}
d_{r+2-i,2}\cdots$ $a_{r+2-i,k}\rightarrow{b_{r+2-i,k}}
{\rightarrow}g_{r+2-i,k+1}\xrightarrow{L_{1,1}^*}f_{i,k+2}\xrightarrow{L^*_{i+(r-1)/2,p}}e_{i,k+1}\xrightarrow{L_{i,(p+1)/2}^*}a_{i,k}$.

The connecting arc is parallel to the component of $L^*$ above the arrow.
(c.f. Figure \ref{dfcn} and Figure \ref{25,5f}.)

Now we finished the construction of $l_{i}$, $1\leq i\leq (r+1)/2$. Note that $\{l_{i}:1\leq i\leq (r+1)/2\}$ are mutually disjoint and $\{L^*_{i,j}: 1\leq i\leq r, 1\leq j\leq p, (i,j)\neq (1,1)\}$ intersects $l_i$ negatively, $1\leq i\leq (r+1)/2$.
(c.f. Figure \ref{dfcn} and Figure \ref{25,5f})

Recall $\breve{L}\cap\breve{Y}_{1,s}=\{\breve{L}_{i,j,s}^{2l-1}\cup\breve{L}_{i,j,s}^{2l+2n-1}, (i,j,l)\in C\}$, $C=\{(i,j,l):1\leq i\leq r, 1\leq j\leq p, (i,j)\neq (1,1), 1\leq l\leq n\}$. 
${\cal C}=\{l_{i}:1\leq i\leq (r+1)/2\}$.
By the construction of ${\cal C}$, the sign of every intersection point of $(L_{i,j}^{2l-1})^*$ and ${\cal C}$ is negative, $(i,j,l)\in C$.
We want to show that $(L_{i,j}^{2l-1})^*\cap{\cal C}\neq \emptyset$, $(i,j,l)\in C$.
The following remark gives the detail.

\begin{remark}\label{intersection2}
When $i=1, 1<j\leq p, 1\leq l\leq n$, $(L_{i,j}^{2l-1})^*\cap l_1\neq\emptyset$ by (\ref{head-1}). 

Consider $1<i\leq (r+1)/2$.
\begin{align*}
 \text{When }&  l=1, 
    \begin{cases}
    (L_{i,j}^{2l-1})^*\cap l_1\neq\emptyset & \text{if }1\leq j\leq (p+1)/2,  \text{ by }(\ref{head-1});\\
    (L_{i,j}^{2l-1})^*\cap b_{i,1}\neq\emptyset& \text{if }(p+1)/2<j\leq p,\text{ by }(\ref{b1}); 
    \end{cases}\\
 \text{When }&   2\leq l\leq k,
     \begin{cases}
      (L_{i,j}^{2l-1})^*\cap b_{i,l-1}\neq\emptyset & \text{if } 1\leq j<(p+1)/2,\text{ by } (\ref{b1});\\
       (L_{i,j}^{2l-1})^*\cap a_{i,l-1}\neq\emptyset &\text{if } j=(p+1)/2,\text{ by } (\ref{a});\\
       (L_{i,j}^{2l-1})^*\cap b_{i,l}\neq\emptyset & \text{if } (p+1)/2<j\leq p,\text{ by } (\ref{b1});\\
     \end{cases} \\ 
 \text{When }&  l=k+1, 
       \begin{cases}
      (L_{i,j}^{2l-1})^*\cap b_{i,k}\neq\emptyset  &\text{if }1\leq j<(p+1)/2, \text{ by } (\ref{b1});\\
    (L_{i,j}^{2l-1})^*\cap a_{i,k}\neq\emptyset  & \text{if }j=(p+1)/2,\text{ by } (\ref{a});\\
     (L_{i,j}^{2l-1})^*\cap l_1\neq\emptyset  & \text{if }(p+1)/2<j\leq p,\text{ by } (\ref{head-1});
     \end{cases}\\
 \text{When }&  k+2\leq l\leq n,(L_{i,j}^{2l-1})^*\cap l_1\neq\emptyset,  \text{ by } (\ref{head-1}).
\end{align*}
Consider $(r+1)/2<i\leq r$.
\begin{align*}
\text{When }&   l=1,
    \begin{cases}
   (L_{i,j}^{2l-1})^*\cap l_1\neq\emptyset  &\text{if }1\leq j< (p+1)/2,\text{ by }(\ref{head-1});\\
    (L_{i,j}^{2l-1})^*\cap b_{i,1}\neq\emptyset  &\text{if }(p+1)/2\leq j< p,\text{ by }(\ref{b2});\\
    (L_{i,j}^{2l-1})^*\cap a_{i,1}\neq\emptyset  &\text{if }j=p,\text{ by } (\ref{a});
     \end{cases}\\
 \text{When }&   2\leq l\leq k,
        \begin{cases}
     (L_{i,j}^{2l-1})^*\cap b_{i,l-1}\neq\emptyset & \text{if } 1\leq j<(p+1)/2,\text{ by } (\ref{b1});\\
     (L_{i,j}^{2l-1})^*\cap b_{i,l}\neq\emptyset & \text{if } (p+1)/2\leq j< p,\text{ by } (\ref{b1});\\ 
       (L_{i,j}^{2l-1})^*\cap a_{i,l}\neq\emptyset &\text{if } j=p,\text{ by } (\ref{a});
     \end{cases}  
 \end{align*}
 \begin{align*}
 \text{When }&  l=k+1, 
       \begin{cases}
   (L_{i,j}^{2l-1})^*\cap b_{i,k}\neq\emptyset   &\text{if }1\leq j<(p+1)/2, \text{ by } (\ref{b2});\\
    (L_{i,j}^{2l-1})^*\cap l_1\neq\emptyset  & \text{if }(p+1)/2\leq j\leq p,\text{ by } (\ref{head-1});
     \end{cases}\\
 \text{When }&  k+2\leq l\leq n,(L_{i,j}^{2l-1})^*\cap l_1\neq\emptyset\text{ by } (\ref{head-1}).
\end{align*}
\end{remark}

By Claim \ref{torus}, we can perform Dehn twist operations to $\breve{F}_{1,s}$ along $\G$, such that the link $\breve{L}-(\breve{L}_{1,1,1}\cup\breve{L}_{1,1,2})$ is transverse to 
the new surface fiber bundle $\breve{\cal F}''_{1,s}$ everywhere.
Then the exterior of $\breve{L}$ in $\breve{Y}$ is a surface bundle. 
Moreover, this exterior is a finite cover of the exterior of $\widetilde{K}$ in $W_k$, which is a free double cover of the exterior of $K$ in $S^3$. 
Thus $K$ is virtually fibered.

Now we finish the proof of Theorem \ref{main} in Case I, when $K$ is a classic Montesinos knot.

\subsection{$K$ has two components}
The proof in this case is generated from the proof when $K$ is a knot, similar to Sec. 6.2 of \cite{abz}. We only give the outline. We still use the same notations introduced in 2.1. 

$K$ has two components, so $f |: \widetilde{K}\rightarrow K^*$ is a trivial 2-fold cover. Then $L=\Psi^{-1}(\widetilde{K})$ is a geodesic link with $2pr$ components. Let $L=\{L_{i,j}^t: 1\leq i\leq r,1\leq j \leq p, t=1,2\}$, where $\hat{f}(L_{i,j}^1\cup L_{i,j}^2)=L^*_{i,j}$.

Let $M$ be the complement  of $L_{1,1}^1$ in $Y$. 
 We have $M=Y_1\cup M_2$, where $M_2=\hat{f}^{-1}(L_{1,1}^*\times[-\epsilon, \epsilon])-\overset{\circ}{N}(L_{1,1}^1)$. 

 $H_2$ and $H_1$ are constructed by the same way as in Sec 6.2 of \cite{abz}. Note that $L_{1, 1}^2$ is a {\it new}  circle fiber in the Seifert fiber structure of $M_2$, so it is transverse to the surface fiber ${\cal F}_2$.

\begin{figure}
\begin{center}
\includegraphics{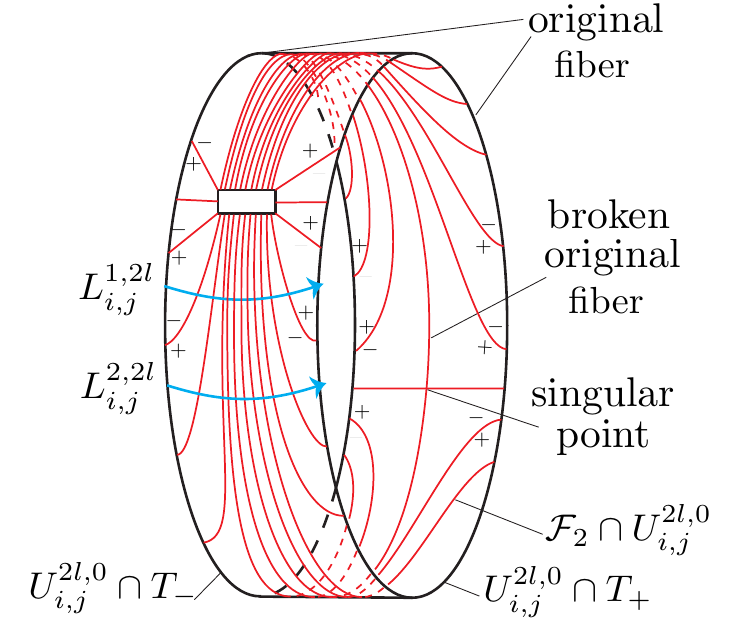}
\end{center}
\caption{\label{foliation2} Foliation of $U_{i,j}^{2l,0}$ when $K$ has two components. }
\end{figure}

 As in Sec. 2.1, $Y_1$ and $M_2$ separate $L_{1,j}^*$ into $2n$ arcs and separate $L_{i,j}^*$ into $2(n-k)$ arcs, $1\leq i\leq r, 1\leq j\leq p, (i,j)\neq (1,1), t=1,2$.
$\hat{f}(L_{i,j}^{1,l}\cup L_{i,j}^{2,l})=(L_{i,j}^l)^*$, $(i,j,l)\in B, t=1,2$.
Now $L_{1, 1}^1$ only intersects the original fiber once, so $U_{i, j}^{2l, 0}=\hat{f}^{-1}((L_{i,j}^{2l})*)\cap M_2$ is a once-punctured annulus, and there is only one singular point on $U_{i, j}^{2l,0}$ as shown in Figure \ref{foliation2} (c.f. Figure 11 in Sec 6.2 of \cite{abz}), $(i,j,l)\in B$.
The following construction is the same as Sec. 2.1.
The proof of Theorem \ref{main} in Case I is finished.

\section{Proof of Theorem \ref{main} in Case II.}

In this section we prove Theorem \ref{main} when 
$K=K(\displaystyle{\frac{q_1}{p},\frac{q_2}{p}, \ldots, \frac{q_{n}}{p}})$ with  $p=2m$, where $m$ is odd, and $n\geq 4$ even. 
We do not consider the case when $p=2$ and $n=4$, since it is already proved in \cite{abz}. 
$K$ is of type $\widetilde{SL_2}$ if $e(W_K)=-(q_1+\cdots+q_n)/p\neq 0$, by (\ref{ew}). 

This proof also follows the outline of Sec. 6.1 of \cite{abz}. 
Recall that $(W_K,\widetilde{K})$ is the 2-fold branched cover of $(S^3,K)$. 
$W_K$ has a Seifert fiber structure: $f:W_K\rightarrow\mathcal{B}_K$.
${\cal B}_K$ is a 2-sphere with $n$ cone points and the order of every cone point is $p$. 
$f(\widetilde{K})=K^*$.
$K^*$ is an equator geodesic which goes through all $n$ cone points $\{c_1,c_2,\ldots, c_n\}$ successively.
The fundamental group of ${\cal B}_K$ is $\G_1$. 
\begin{equation*}
\Gamma_1=\{x_1,x_2,\ldots, x_n: x_1^p=x_2^p=\cdots=x_n^p=x_1x_2\cdots x_n=1\}.
\end{equation*}
$x_i$ presents $\p N(c_i)$, where $N(c_i)$ is a small regular neighborhood of $c_i$ on ${\cal B}_K$, $1\leq i\leq n$.
Let ${\cal B}_K^0={\cal B}_K-\cup_{i=1}^n\overline{N(c_i)}$, an 
$n$-punctured sphere.
$W_K$ is constructed by gluing solid tori along the boundary of ${\cal B}_K^0\times S^1$.
Let $h$ be a simple closed circle in $W_K$ representing the $S^1$ factor. 
The meridian disk of the solid torus is attached to the circle on $\p N(c_i)\times S^1$ of slop $x_i^{p}h^{q_i}$, $1\leq i\leq n$.
Then 
\begin{equation*}
\pi_1(W_K)=\{x_1,x_2,\ldots, x_n, h: x_i^{p}h^{q_i}= [x_i,h]= x_1x_2\cdots x_n=1, 1\leq i\leq n\}.
\end{equation*}

Since $p$ is even, $K$ has $n$ components by (\ref{cn}). 
Then $\widetilde{K}$ also has $n$ components denoted $\{\widetilde{K}_1,\widetilde{K}_2,\ldots, \widetilde{K}_n\}$.  $f(\widetilde{K}_i)$ is the segment of $K^*$ between two cone points, $1\leq i\leq n$.
Suppose that $c_{i}$ and $c_{i+1}$ are two end points of $f(\widetilde{K}_i)$, $1\leq i\leq n$.
(We always consider the index $i$ as mod $n$.) 
By the construction of $W_K$, $\widetilde{K}_i$ corresponds to the element $(x_ix_{i+1})^{p/2}h^{(q_{i}+q_{i+1})/2}=(x_ix_{i+1})^mh^{(q_{i}+q_{i+1})/2}$ in $\Gamma_1$.
($q_i$ and $p$ are relatively prime so $q_i$ is odd. Then $(q_{i}+q_{i+1})/2$ is an integer.) 
If we think $K^*$ as a 1-orbifold and the cone points as mirror points,
 $f|:\widetilde{K}\rightarrow K^*$ is a 2-fold orbifold cover.

At first, we construct a proper covering space of ${\cal B}_K$, $\psi: F\rightarrow {\cal B}_K$, such that $F$ is a smooth oriented closed surface. We also construct this cover by two steps.

\begin{figure}
\begin{center}
\includegraphics{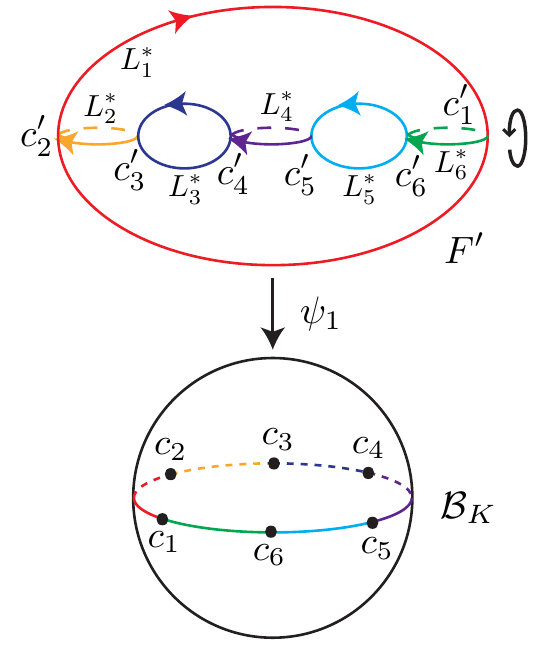}
\end{center}
\caption{\label{2fold} $F'$, when $n=6$.}
\end{figure}

Since $n$ is even, there exists a homomorphism: $h_1: \Gamma_1\rightarrow \z_2$, where 
\begin{equation*}
h_1(x_i)=\bar{1},\ \ \  1\leq i\leq n. 
\end{equation*} 
Let $\psi_1:F'\rightarrow {\cal B}_K$ be the covering space of ${\cal B}_K$ corresponding to $h_1$.  
$F'$ is the two-fold branched covering space of ${\cal B}_K$ with the branch set $\{c'_1,c'_2,\ldots, c'_n\}$,
(c.f. Figure \ref{2fold}),
where  $c'_i=\psi_1^{-1}(c_i)$, $1\leq i\leq n$.
$c'_i$ is a cone point of order $m$, $1\leq i\leq n$. 
Let $\Gamma_2$ be the fundamental group of $F'$. $\Gamma_2$ has the following expression:
\begin{equation*}
\Gamma_2=\{a_1,b_1,\ldots, a_g,b_g, y_{1},\ldots, y_{n}: 
y_{1}^{m}=\cdots=y_{n}^{m}=1, \prod_{i=1}^{g}[a_i,b_i] y_{1}\cdots  y_{n}=1\}
\end{equation*} 
where $y_i$ is represented by a small circle on $F'$ centered at $c'_i$, $g=n/2-1$,  $1\leq i\leq n$. 

\begin{figure}
\begin{center}
\includegraphics{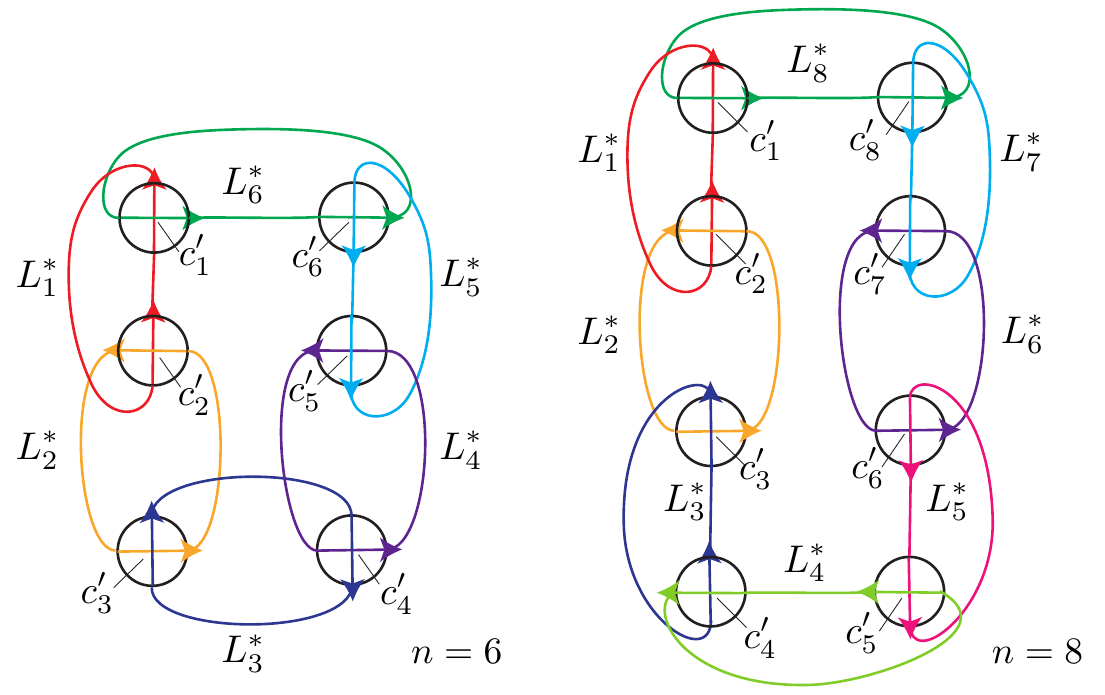}
\end{center}
\caption{\label{Listar} $\{L_1^*,L_2^*,\ldots, L_n^*\}$, when $n=6$ and $n=8$.}
\end{figure}

The covering transformation $\t_1$ is a rotation with angle $\pi$ about the axis which punctures $F'$ at the points $\{c'_i:1\leq i\leq n\}$.  
 $\psi_{1}^{-1}(K^*)$ has $n$ components, $\{L_1^*,L^*_2,\ldots, L_n^*\}$.
$L^*_i$ is a geodesic and goes through $c'_i$ and $c'_{i+1}$, $1\leq i\leq n$ . 
We have the following result:
\begin{equation}\label{A2}
L^*_i\cap L^*_j=
  \begin{cases}
  c'_{j} &\text{ if } i=j-1,\\
  c'_{i} &\text{ if } i=j+1,\\
  \emptyset &\text{otherwise}.
  \end{cases}
\end{equation}
If $L^*_i$ and $L^*_j$ intersect at point $P$, we define {\it the angle between $L_{i}^*$ and $L^*_j$ at a point $P$} to be  $\angle (L^*_{i},L^*_{j})$. ( i.e. If we rotate $L^*_{i}$ to $L^*_j$ around $P$ counterclockwise in a small neighborhood of $P$, the angle that $L^*_{i}$ passes through is $\angle (L^*_{i},L^*_{j})$.)
Note that we require $\angle (L^*_{i}, L^*_{j})\in [0,2\pi)$.
  Then $\angle(L^*_{i},L^*_{j})=2\pi-\angle(L^*_{j},L^*_{i})$.
Orient $F'$ and  $\{L_{i}^*: 1\leq i\leq n\}$, such that $\angle (L^*_{i-1},L^*_i)=2\pi/4m=\pi/2m$ at $c'_{i}$, (note that $c'_i$ is a cone point of order $m$), $1\leq i\leq n$.
Figure \ref{Listar} gives examples when $n=6$ and $8$.

There is a homomorphism $h_2: \Gamma_2\rightarrow \z/m$ which is defined as following.
\begin{equation}\label{h2}
\begin{split}
\text{If } 1\leq i\leq n/2,\ 
h_2(y_i)=&
\begin{cases}
\bar{1} &\text{ if }  i=4s+1 \text{ or }4s+2,\\
-\bar{1} &\text{ if } i=4s \text{ or }4s+3;
\end{cases}\\
\text{if } n/2< i\leq n,\ h_2(y_i)=&-h_2(y_{n/2-(i-n/2-1)})=-h_2(y_{n-i+1}),
\end{split}
\end{equation}
where $s\in\z$.
Let $\psi_2:F\rightarrow F'$ be the covering of $F'$ corresponding to $h_2$ where $F$ is the covering space. 
Since the order of $h_2(y_i)$ is $m$ that is equal to the order of the cone point $c'_i$, $F$ is a smooth closed orientable surface without cone points, $1\leq i\leq n$.
Denote $\hat{c}_{i}=\psi_2^{-1}(c'_{i})$, $1\leq i\leq n$.

\begin{figure}
\begin{center}
\includegraphics{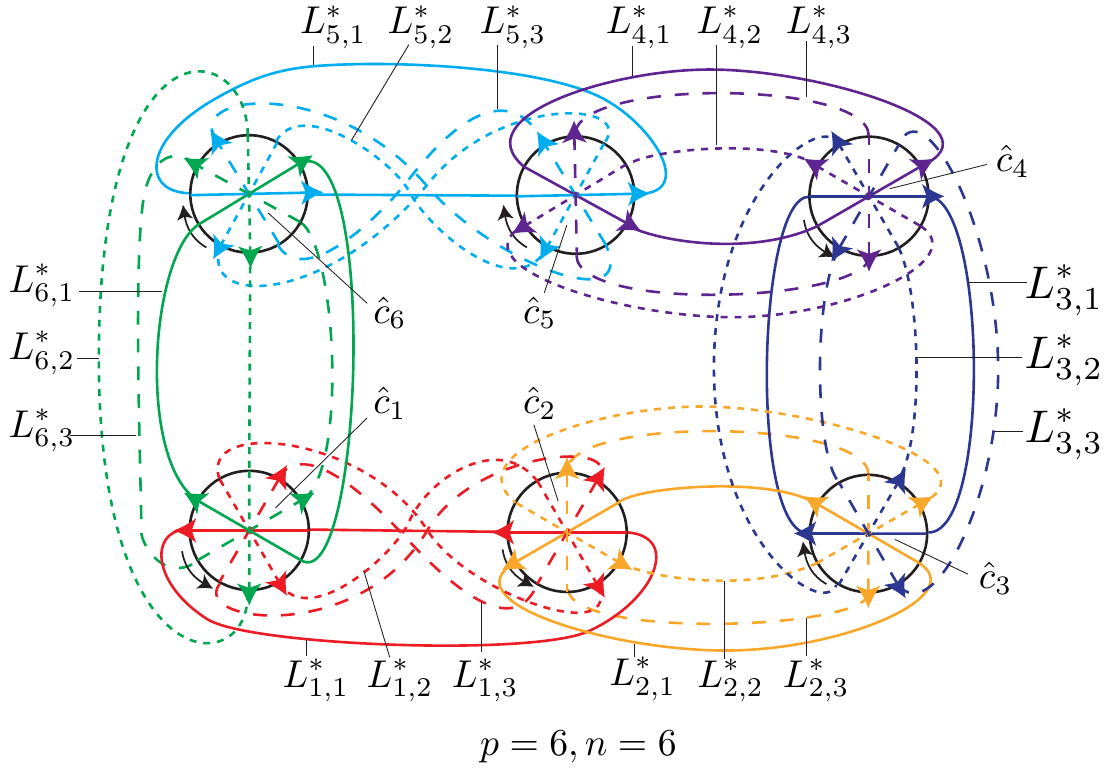}
\end{center}
\caption{\label{Lijstar6} $\{L_{i,j}^*: 1\leq i\leq n, 1\leq j\leq p/2\}$, when $p=6, n=6$.}
\end{figure}

The preimage of $L^*_i$ is a set of $m$ geodesics $\{L^*_{i,1},\cdots,L^*_{i,m}\}$, $1\leq i\leq n$.
Let $L^*=\psi_2^{-1}(K^*)$. Then $L^*$ has $nm$ components, denoted by $L_{i,j}^*$, $1\leq i\leq n, 1\leq j\leq m$.
Each $L^*_{i,j}$ goes through $2$ points $\hat{c}_{i}$ and $\hat{c}_{i+1}$, $1\leq i\leq n, 1< j\leq m$. 
Let $\t_2$ be the deck transformation of $\psi_2$ corresponding to $\bar{1}\in\z/m$. 
We may assume that $L^*_{i,j+1}=\t_2(L^*_{i,j})$, $1\leq i\leq n, 1\leq j\leq m$. 
Orient $L^*_{i,1}$, such that $\angle(L^*_{i-1,1},L^*_{i,1})=\pi/2m$. Then give $L_{i,j}^*$ the induced orientation, $1\leq i\leq n,1< j\leq m$.  
$F$ admits an orientation such that $\t_2$ is a counterclockwise rotation near the fixed points $\hat{c}_{i}$ by $2\pi/m$ if $h(y_i)=\bar{1}$, $1\leq i\leq n$. Then $\t_2$ is a clockwise rotation near the fixed points $\hat{c}_{i}$ by $2\pi/m$ if $h(y_i)=-\bar{1}$, $1\leq i\leq n$.
By the construction of $F$, $\{L_{i,j,}^*: 1\leq i\leq n, 1\leq j\leq m\}$ only intersect at $\{\hat{c_i}: 1\leq i\leq n\}$, and $L^*_{i_1,j_1}\cap L^*_{i_2,j_2}\neq\emptyset $ if and noly if $|i_1-i_2|\leq 1$, $1\leq i_1, i_2\leq n,1\leq j_1,j_2\leq m$.
In addition, when $L^*_{i_1,j_1}\cap L^*_{i_2,j_2}\neq \emptyset$, if   $i_1\neq i_2$, $L^*_{i_1,j_1}\cap L^*_{i_2,j_2}=\hat{c}_i$ where $i=max(i_1,i_2)$;  if  $i_1=i_2=i$, $L^*_{i,j_1}\cap L^*_{i,j_2}=\hat{c}_i\cup\hat{c}_{i+1}$, $1\leq i_1,i_2,i\leq n$, $1\leq j_1,j_2\leq m$.
See Figure \ref{Lijstar6}
for $p=6,n=6$. 
The following remark gives the detail.

\begin{remark}\label{angle2}
We assume that $i, i+1, i-1$ is in $\z/n$. 

 $\displaystyle{\angle(L^*_{i-1,1}, L^*_{i,1})=\frac{\pi}{2m}, \ \angle(L^*_{i,1}, L^*_{i-1,1})=-\frac{\pi}{2m}+2\pi}$.

If $h_2(y_i)=\bar{1}$, 
\begin{equation*}
\begin{split}
\angle(L^*_{i,1}, L^*_{i-1,j})&=\displaystyle{-\frac{\pi}{2m}+\frac{(j-1)2\pi}{m}}, \ \ \angle(L^*_{i,1}, L^*_{i,j})=\displaystyle{\frac{(j-1)2\pi}{m}},\\
  \angle(L^*_{i,1}, L^*_{i+1,j})&=\displaystyle{\frac{\pi}{2m}+\frac{(j-1)2\pi}{m}}, 
1<j\leq m.
\end{split}
\end{equation*}

If $h_2(y_i)=-\bar{1}$, 
\begin{equation*}
\begin{split}
\angle(L^*_{i,1},L^*_{i-1,j})&=\displaystyle{2\pi-\frac{\pi}{2m}-\frac{(j-1)2\pi}{m}},\ \ \angle(L^*_{i,1}, L^*_{i,j})=\displaystyle{2\pi-\frac{(j-1)2\pi}{m}},\\
 \angle(L^*_{i,1},L^*_{i+1,j})&=\displaystyle{2\pi+\frac{\pi}{2m}-\frac{(j-1)2\pi}{m}}, 
 1<j\leq m.
 \end{split}
\end{equation*}
\end{remark}

Next we want to construct a covering space of $W_K$.
There is a nature homomorphism: $h_0: \pi_1(W_K)\rightarrow \pi_1({\cal B}_K)$ by sending $h\in \pi_1(W_K)$ to 1. 
The covering $\Psi_1: Y'\rightarrow W_K$ is induced by the composition $h_1\circ h_0$. 
$F'$ is the base orbifold of $Y'$.

\begin{figure}
\begin{center}
\includegraphics{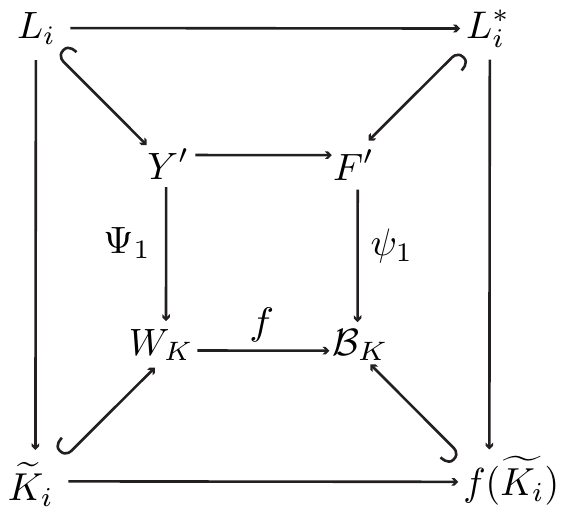}
\end{center}
\caption{\label{diagram2} Construction of $Y'$.}
\end{figure}

We have the diagram (Figure \ref{diagram2}) is commutative, similar to the diagram (5) in Sec. 2 of \cite{abz}.

Recall that $\widetilde{K}$ has $n$ components, and the fundamental class of every component $\widetilde{K}_i$ is $(x_ix_{i+1})^mh^{(q_i+q_{i+1})/2}$. 
\begin{equation*}
h_1(h_0((x_ix_{i+1})^mh^{(q_i+q_{i+1})/2}))=h_1((x_ix_{i+1})^m)=m(\bar{1}+\bar{1})=\bar{0} \text{ in }  \z_2.
\end{equation*}
Then $\Psi_1^{-1}(\widetilde{K}_i)=L_i$ has two components, $1\leq i\leq n$. 
 $\Psi_1^{-1}(\widetilde{K})$ is a geodesic link, and it has exactly $2n$ components. 
Then $L^*_i$ is the Seifert quotient of $L_i$, and $L_i$ is a 2-fold trivial cover of $L^*_i$, $1\leq i\leq n$.

Similarly, let $\Psi_2: Y\rightarrow Y'$ be the $m$-fold covering of $Y'$ induced by $h_2$ and the Seifert structure of $Y'$. 
Since $m$ is odd, by the similar discussion in Sec. 6.1 of \cite{abz}, 
the preimage of every component of $L_i$ has $m$ components.
Define $\Psi=\Psi_1\circ\Psi_2, \psi=\psi_1\circ\psi_2$. 
$\Psi: Y\rightarrow W_K$, $\psi:F\rightarrow {\cal B}_K$.
$Y$ has a locally-trivial circle bundle structure. 
Let $\hat{f}: Y\rightarrow F$ be the Seifert quotient map.
Let $L=\Psi^{-1}(\widetilde{K})$. $L$ has $2nm$ components. Denote
\begin{equation*}
L=\{L_{i,j,k}: 1\leq i\leq n, 1\leq j\leq m, k=1,2\}.
\end{equation*}
$\hat{f}(L)=L^*$. $L^*$ has $nm$ components.
We may suppose that $\hat{f}(L_{i,j,1}\cup L_{i,j,2})=L^*_{i,j}, 1\leq i\leq n, 1\leq j\leq m$.
$L_{i,j,1}\cup L_{i,j,2}$ is a 2-fold trivial cover of $L_{i,j}^*$, $1\leq i\leq n, 1\leq j\leq m$.

We assume that the index $i\in \z/n$, $j\in \z/m$ for $i,j$ in $L_{i,j,k}$ or $L^*_{i,j}$.

Next we want to prove that the exterior of some components of $L$ in $Y$ is a surface semi-bundle.
Since $p$ is even, $K$ has $n$ components. There is no component of $L^*$ which goes through all $\hat{c}_i$'s, but there are $n/2$ mutually disjoint components of $L^*$, $\{L_{2i-1,1}^*:1\leq i\leq n/2\}$, pass all $\hat{c}_i$'s, $1\leq i\leq n$. 
Different from Prop. \ref{semibundle}, we need consider the exterior of more than one components. 
\begin{prop}\label{semi}
The exterior of $\{L_{2i-1,1,1}:1\leq i\leq n/2\}$ in $Y$ is a surface semi-bundle.
\end{prop}

Before we prove this proposition, we introduce some notations and preparations.

We can reorient $\{L^*_{i,j}: 1\leq i\leq n, 1\leq j\leq m\}$ so that 
\begin{equation}\label{neworientation}
\begin{split}
i(L^*_{2i-1,1},L^*_{2i-1,j_1},)&=i(L^*_{2i-1,1}, L^*_{2i-2,j_2})
=sign(h_2(y_{2i-1})) \text{ at }\hat{c}_{2i-1}\\
i(L^*_{2i-1,1}, L^*_{2i-1,j_1})&=i(L^*_{2i-1,1}, L^*_{2i,j_2})
=sign(h_2(y_{2i})) \text{ at } \hat{c}_{2i}
\end{split}
\end{equation}
 $1\leq i\leq n/2, 1< j_1\leq m, 1\leq j_2\leq m$. i.e. 
At $\hat{c}_{2i-1}$, 
\begin{equation*}
\angle(L^*_{2i-1,1}, L^*_{2i-1,j}),\angle(L^*_{2i-1,1}, L^*_{2i-2,j}), \angle(L^*_{2i-1,1},L^*_{2i-2,1})\in
\begin{cases}
(\pi, 2\pi) &\text{if } h_2(y_{2i-1})=\bar{1},\\
(0,\pi) & \text{if }  h_2(y_{2i-1})=-\bar{1}, 
\end{cases}
\end{equation*}
and at $\hat{c}_{2i}$,
\begin{equation*}
\angle(L^*_{2i-1,1}, L^*_{2i-1,j}),\angle(L^*_{2i-1,1}, L^*_{2i,j}), \angle(L^*_{2i-1,1},L^*_{2i,1})\in
\begin{cases}
(\pi, 2\pi) &\text{if } h_2(y_{2i})=\bar{1},\\
(0,\pi) & \text{if }  h_2(y_{2i})=-\bar{1}, 
\end{cases}
\end{equation*}
$1\leq i\leq n/2, 1< j\leq m.$

By Remark \ref{angle2}, we need to change the orientation of the  components of $L^*$ in Table \ref{orientL1}.

\begin{table}
\begin{center}
\renewcommand{\arraystretch}{1.3}
\begin{tabular}[b]{|c|c|c||c|c|c|}
\hline
\multicolumn{3}{|c||}{At $\hat{c}_{2i-1}$} & \multicolumn{3}{|c|}{At $\hat{c}_{2i}$}\\
\hline $h_2(y_{2i-1})$ & $\bar{1}$ & $-\bar{1}$ & $h_2(y_{2i})$ & $\bar{1}$ & $-\bar{1}$\\
\hline\hline $L_{2i-1,j}^*$ & $1<j\leq \frac{m+1}{2}$ & $1<j\leq \frac{m+1}{2}$ & $L_{2i-1,j}^*$ & $1<j\leq \frac{m+1}{2}$ & $1<j\leq \frac{m+1}{2}$\\
\hline $L_{2i-2,j}^*$ & $1<j\leq \frac{m+1}{2}$ & $1\leq j\leq \frac{m+1}{2}$ & $L_{2i,j}^*$ & $1\leq j\leq \frac{m+1}{2}$ & $1< j\leq \frac{m+1}{2}$ \\ \hline
\end{tabular}
\end{center}
\caption{\label{orientL1} The components of $L^*$ whose orientations need to be changed ($1\leq i\leq n/2$)}
\end{table}

At first, we consider the case $4|n$. 
In this case, $h_2(y_{4i-3})=h_2(y_{4i-2})=\bar{1}$ and $h_2(y_{4i-1})=h_2(y_{4i})=-\bar{1}$, by (\ref{h2}), $1\leq i\leq n/4$. 
By Table \ref{orientL1}, we need to change the orientations of the following components of $L^*$.
\begin{equation*}
\begin{split}
\{L_{4i-3,j}^*,1<j\leq \frac{m+1}{2}\}&\cup\{L_{4i-4,j}^*,1<j\leq \frac{m+1}{2}\} \text{ at } \hat{c}_{4i-3};\\
  \{L_{4i-3,j}^*,1<j\leq \frac{m+1}{2}\}&\cup\{L_{4i-2,j}^*,1\leq j\leq \frac{m+1}{2}\} \text{ at }  \hat{c}_{4i-2};\\
\{L_{4i-1,j}^*,1<j\leq \frac{m+1}{2}\}&\cup\{L_{4i-2,j}^*,1\leq j\leq \frac{m+1}{2}\}\text{ at }\hat{c}_{4i-1};\\
 \{L_{4i-1,j}^*,1<j\leq \frac{m+1}{2}\}&\cup\{L_{4i,j}^*,1<j\leq \frac{m+1}{2}\}\text{ at }\hat{c}_{4i}.
 \end{split}
 \end{equation*}
 In summary, when $4|n$, we need to change 
\begin{equation*}
\{L^*_{i,j}: 1\leq i\leq n, 1< j\leq \frac{m+1}{2}\}\cup\{L^*_{4i-2,1}: 1\leq i\leq \frac{n}{4}\}.
\end{equation*}

Next, consider the case: $4\nmid n$.
According to (\ref{h2}),
\begin{align*}
\text{When }1\leq i\leq n/2,
h_2(y_{i})=
\begin{cases}
\bar{1} & \text{if } i\equiv 1 \text{ or } 2 \text{ mod } 4,\\
-\bar{1} & \text{otherwise}.
\end{cases}\\
\text{When } n/2<i\leq n,
h_2(y_{i})=
\begin{cases}
-\bar{1} & \text{if } i\equiv 1 \text{ or } 2 \text{ mod } 4, \\
\bar{1} & \text{otherwise}.
\end{cases}
\end{align*}
According to Table \ref{orientL1}, we need to change the orientations of the following components of $L^*$.
\begin{equation*}
\text{When }1\leq i\leq n/2,
\begin{cases}
\{L_{i,j}^*,1<j\leq \frac{m+1}{2}\}\cup\{L_{i-1,j}^*,1< j\leq \frac{m+1}{2}\}, & \text{if } i\equiv 1 \text{ mod } 4;\\
\{L_{i-1,j}^*,1<j\leq \frac{m+1}{2}\}\cup\{L_{i,j}^*,1\leq j\leq \frac{m+1}{2}\}, & \text{if } i\equiv 2 \text{ mod } 4;\\
\{L_{i,j}^*,1<j\leq \frac{m+1}{2}\}\cup\{L_{i-1,j}^*,1\leq j\leq \frac{m+1}{2}\}, & \text{if } i\equiv 3 \text{ mod } 4;\\
\{L_{i-1,j}^*,1<j\leq \frac{m+1}{2}\}\cup\{L_{i,j}^*,1< j\leq \frac{m+1}{2}\}, & \text{if } i\equiv 0 \text{ mod } 4.
 \end{cases}
\end{equation*}
\begin{equation*}
\text{When }n/2< i\leq n,
\begin{cases}
\{L_{i,j}^*,1<j\leq \frac{m+1}{2}\}\cup\{L_{i-1,j}^*,1\leq j\leq \frac{m+1}{2}\}, & \text{if } i\equiv 1 \text{ mod } 4;\\
\{L_{i-1,j}^*,1<j\leq \frac{m+1}{2}\}\cup\{L_{i,j}^*,1< j\leq \frac{m+1}{2}\}, & \text{if } i\equiv 2 \text{ mod } 4;\\
\{L_{i,j}^*,1<j\leq \frac{m+1}{2}\}\cup\{L_{i-1,j}^*,1< j\leq \frac{m+1}{2}\}, & \text{if } i\equiv 3 \text{ mod } 4;\\
\{L_{i-1,j}^*,1<j\leq \frac{m+1}{2}\}\cup\{L_{i,j}^*,1\leq j\leq \frac{m+1}{2}\}, & \text{if } i\equiv 0 \text{ mod } 4.
 \end{cases}
\end{equation*}
Note that there is no conflict if one component $L^*_{i,j}$ connecting two points $\hat{c}_i$ and $\hat{c}_{i+1}$ when $i,i+1$ both greater than $n/2$ or both less than $n/2$, $1\leq j\leq m$. 
Now we check $L_{n/2,j}^*$ that goes through $\hat{c}_{n/2}$ and $\hat{c}_{n/2+1}$, $1\leq j\leq m$. $4\nmid n$, so $n/2$ is odd. 
From the above discussion we need to change $\{L^*_{n/2,j}:1< j\leq (m+1)/2\}$ at both $\hat{c}_{n/2}$ and $\hat{c}_{n/2+1}$.
Then this reorientation works when $4\nmid n$.

\begin{figure}
\begin{center}
\includegraphics{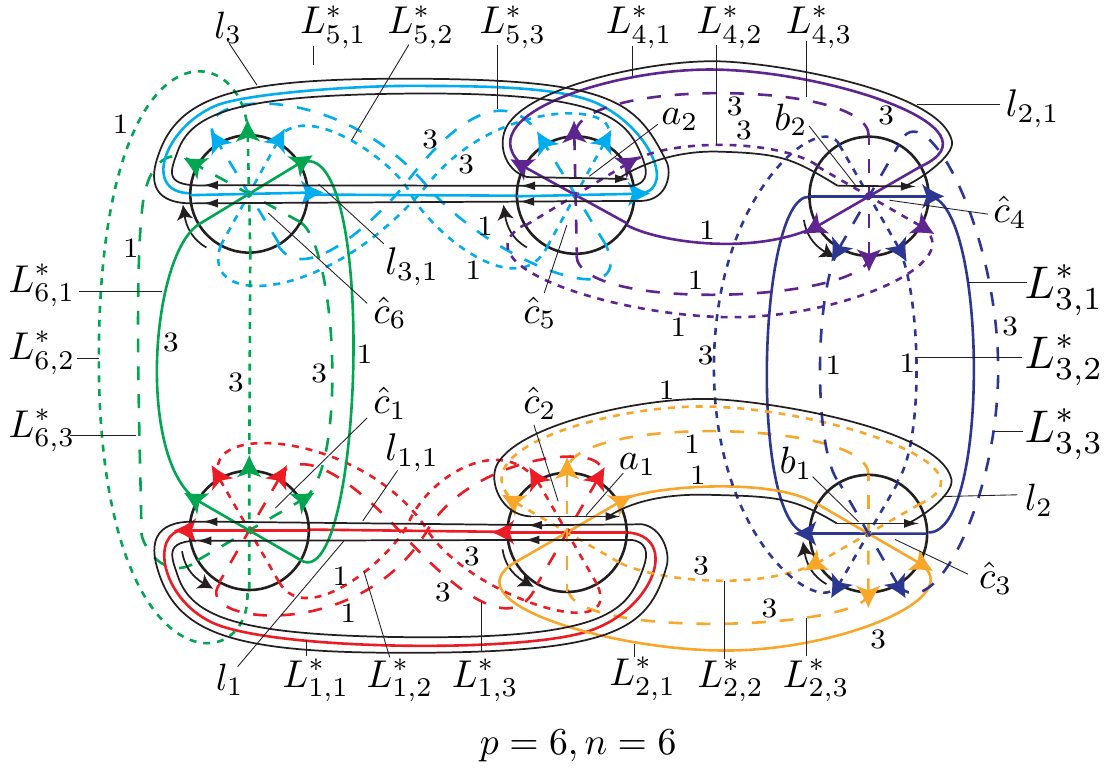}
\end{center}
\caption{\label{Lijstar6c} Reorientation of $\{L_{i,j}^*: 1\leq i\leq n, 1\leq j\leq m\}$.}
\end{figure}

In summary, when $4\nmid n$, we need change 
\begin{equation*}
\begin{split}
\{L^*_{i,j}&: 1\leq i\leq n, 1< j\leq \frac{m+1}{2}\}\\
\cup
\{L^*_{i,1}&: 4|(i-2) \text{ and } 1\leq i\leq \frac{n}{2},\text{or } 4|i \text{ and } \frac{n}{2}<i\leq n\}.
\end{split}
\end{equation*}
In the following discussion, we consider $\{L_{i,j}^*: 1\leq i\leq n, 1\leq j\leq m\}$ with the new orientation. Figure \ref{Lijstar6c} shows the new orientation when $p=6, n=6$.

By the construction of $F$, $\{L^*_{2i-1,1}: 1\leq i\leq n/2\}$ are mutually disjoint. 
Let $F_2^i=L_{2i-1,1}^*\times [-\e,\e]$, where $\e$ is a small positive number such that $F_2^i$'s are mutually disjoint, $1\leq i\leq n/2$.
Denote $\partial F_2^i=\b_-^i\cup\b_+^i$, where $\b_-^i=L_{2i-1,1}^*\times\{-\e\}$, $\b_+^i=L_{2i-1,1}^*\times\{\e\}$, $1\leq i\leq n/2$.
We may suppose that $\b_-^i$ is on the left side of $L_{2i-1,1}^*$.
Define
\begin{equation*}
F_1={F-(\underset {i=1}{\overset{n/2}{\cup}}\overset{\circ}{F_2^i})}.
\end{equation*}
\begin{equation*}
\text{Set }A=\{i,j:1\leq j\leq m \text{ when }i=2s,\text{or }1<j\leq m\text{ when } i=2s-1. 1\leq s\leq n/2\}.
\end{equation*}
$L_{i,j}^*$ is  separated into four parts by $F_1$ and $\{F_2^i: 1\leq i\leq n/2\}$, denoted by $\{(L_{i,j}^k)^*: 1\leq k\leq 4\}$, $(i,j)\in A$.
$(L_{2i-1,j}^2)^*\cup (L_{2i-1,j}^4)^* \subset F_2^i$ and $\hat{c}_{2i-1}\in (L_{2i-1,j}^2)^*$ and $\hat{c}_{2i}\in (L_{2i-1,j}^4)^*$, $1\leq i\leq n/2, 1< j\leq m$.
$\hat{c}_{2i+1}\in(L_{2i,j}^2)^* \subset F_2^{i+1}$ and $\hat{c}_{2i}\in(L_{2i,j}^4)^* \subset F_2^{i}$, $1\leq i\leq n/2, 1\leq j\leq m$.
$(L_{i,j}^1)^*\cup (L_{i,j}^3)^*\in F_1$, $(i,j)\in A$. 
We assume that $(L_{i,j}^1)^*, (L_{i,j}^2)^*, (L_{i,j}^3)^*, (L_{i,j}^4)^*$ are connected successively along the orientation of $L_{i,j}^*$, $(i,j)\in A$.
In Figure \ref{Lijstar6c},  the small number next to the segment $(L_{i,j}^k)^*$ is the index $k$, $(i,j)\in A, 1\leq k\leq 4$.
Then we have the following remark.
\begin{remark}\label{arc}
$1\leq i\leq n/2$.

If $h_2(y_{2i-1})=h_2(y_{2i})=\bar{1}$, $(L_{2i-1,j}^1)^*$ and $(L_{2i-1,j}^3)^*$ are from $\b_+^i$ to $\b_-^i, 1<j\leq m$.

If $h_2(y_{2i-1})=h_2(y_{2i})=-\bar{1}$, $(L_{2i-1,j}^1)^*$ and $(L_{2i-1,j}^3)^*$ are from $\b_-^i$ to $\b_+^i, 1<j\leq m$.

If $h_2(y_{2i-1})=-\bar{1},h_2(y_{2i})=\bar{1}$, $(L_{2i-1,j}^1)^*$ is from $\b_+^i$ to $\b_+^i$,
and $(L_{2i-1,j}^3)^*$ is from $\b_-^i$ to $\b_-^i, 1<j\leq m$.

If $h_2(y_{2i-1})=\bar{1},h_2(y_{2i})=-\bar{1}$, $(L_{2i-1,j}^1)^*$ is from $\b_-^i$ to $\b_-^i$,
and $(L_{2i-1,j}^3)^*$ is from $\b_+^i$ to $\b_+^i, 1<j\leq m$.

When $h_2(y_{2i})=\bar{1}, h_2(y_{2i+1})=-\bar{1}$,
$(L_{2i,j}^1)^*$ is from $\b_+^{i}$ to $\b_+^{i+1}$, $(L_{2i,j}^3)^*$  is from $\b_-^{i+1}$ to $\b_-^{i}$, $1\leq j\leq m$.

If $h_2(y_{2i})=-\bar{1}, h_2(y_{2i+1})=\bar{1}$,
$(L_{2i,j}^1)^*$ is from  $\b_-^{i}$ to $\b_-^{i+1}$, $(L_{2i,j}^3)^*$  is from  $\b_+^{i+1}$ to $\b_+^{i}$, $1\leq j\leq m$.\\
$\{(L_{2i-1,j}^2)^*: 1<j\leq m\} \cup\{(L_{2i-2,j}^2)^*: 1\leq j\leq m\}$ is from 
$\begin{cases}
\b_-^i \text{ to } \b_+^i & \text{ if } h_2(y_{2i-1})=\bar{1};\\
\b_+^i \text{ to } \b_-^i & \text{ if } h_2(y_{2i-1})=-\bar{1}.
\end{cases}$
$\{(L_{2i-1,j}^4)^*: 1<j\leq m\} \cup\{(L_{2i,j}^4)^*: 1\leq j\leq m\}$ is from 
$\begin{cases}
\b_-^i \text{ to } \b_+^i & \text{ if } h_2(y_{2i})=\bar{1};\\
\b_+^i \text{ to } \b_-^i & \text{ if } h_2(y_{2i})=-\bar{1}.
\end{cases}$
\end{remark}

We need the following lemma for Proposition \ref{semi}
\begin{lemma}\label{connected}
$F_1$ is connected.
\end{lemma}

\ni\textbf{Proof}: Suffice to prove that the boundary components of
$F_1$, which is the set $\{\b_-^i\cup\b_+^i, 1\leq i
\leq n/2\}$, can be mutually connected to each other by arcs in
$F_1$.

 From Remark \ref{arc}, $\b_+^i$ and $\b_-^i$ are connected by $(L_{2i-1,j}^1)^*$, $1\leq i\leq n/2$, $1<j\leq m$. 
 $\b_-^i$ and $\b_-^{i+1}$ are connected by $(L_{2i,j}^1)^*$ or $(L_{2i,j}^3)^*$, $1\leq i\leq n/2$, $1\leq j\leq m$.
Similarly, $\b_+^i$ and $\b_+^{i+1}$ are also connected by $(L_{2i,j}^1)^*$ or $(L_{2i,j}^3)^*$, $1\leq i\leq n/2$, $1\leq j\leq m$.
$\square$

Let $T^i=\hat{f}^{-1}(L_{2i-1, 1}^*)\subset Y$, which is a vertical
torus over $L_{2i-1,1}^*$, $1\leq i
\leq n/2$. 
Then $L_{i, j,t}$ is transverse to $T^i$ for all $(i,j)\in A$ and $t=1,2$.
Since $L_{2i-1,1}^*$ is a geodesic, the torus  $T^i$ is a totally geodesic torus which inherits a Euclidean structure from the $\widetilde{SL_2}$ structure on $Y$, $1\leq i\leq n/2$. 

Next we give the proof of Proposition \ref{semi}

\ni\textbf{Proof of Proposition \ref{semi}}: 

Suppose $1\leq i\leq n/2$ in this proof. 

Let $Y_2^i$ be the submanifold of $Y$ lying over $F_2^i$.
$Y_2^i=\hat{f}^{-1}(L_{2i-1, 1}^*\times [-\e, \e])$. 
\begin{equation*}
Y_1={Y-(\underset {i=1}{\overset{n/2}{\cup}}\overset{\circ}{Y_2^i})}
\end{equation*}
is the submanifold of $Y$ lying over $F_1$. 
By Lemma \ref{connected}, $F_1$ is connected so $Y_1$ is connected.
Define
\begin{equation*}
T_-^i=\hat{f}^{-1}(\b_-^i), \ T_+^i=\hat{f}^{-1}(\b_+^i).\end{equation*}

Let $Y_0$ be the union of $Y_1$ and $\{Y_2^i: 1\leq i\leq n.2\}$ where $Y_2^i$ is glued  to $Y_1$ along $T_-^i$. 
By construction, each Seifert fibration of  $Y_0, Y_1$, and  $Y_2^i$ is a trivial circle bundle. 
Let $F_0$ be the base of $Y_0$ in the Seifert fiber construction. 
$F_0$ is obtained by cutting $F$ open along $\{\b_+^i: 1\leq i\leq n/2\}$.  
We give the circle fibers of $Y_0$ a consistent orientation. 
Choose a horizontal section $B_0$ of the bundle $Y_0\rightarrow F_0$ such that $T^i\cap B_0$ is a geodesic. 
Let $B_1= B_0\cap Y_1$, $B_2^i=B_0\cap Y_i$.
Fix an orientation of $B_0$ and let $B_1, B_2^i$  and their boundaries have the induced orientation. 
We denote the tori $T_\pm^i $ by  $T_{1, \pm}^i$ and $T_{2, \pm}^i$ when
we think of them as the boundary components of  $Y_1$ and $Y_2^i$ respectively. 
Let $\phi_{k,\pm}^i$ be a fixed circle fiber  in the torus $T_{k,\pm}^i$ for $k=1, 2$. 
Let $\a_{k,\pm}^i=B_0\cap T_{k,\pm}^i$.
Then $\a_{k, -}^i=-\a_{k, +}^i$, $\phi _{k, -}^i=\phi _{k, +}^i$, and
$\{\a _{k, \pm}^i, \phi _{k, \pm}^i\}$ forms a basis of $H_1(T_{k,\pm}^i)$, $k=1, 2$.

The Seifert manifold $Y$ can be obtained by gluing $Y_1$ and $Y_2^i$ along $T_{-}^i$ and $T_{+}^i$.
Suppose $g_\pm^i: T_{2, \pm}^i\rightarrow T_{1, \pm}^i$ is the gluing map satisfying the following conditions.
\begin{align*}
(g_-^i)_*(\a_{2, -}^i)&=-\a_{1, -}^i, & (g_-^i)_*(\phi _{2, -}^i)&=\phi _{1, -}^i,\\
(g_+^i)_*(\a_{2, +}^i)&=-\a _{1, +}^i+e^i\phi _{1, +}^i, & (g_+^i)_*(\phi _{2, +}^i)&=\phi _{1, +}^i,
\end{align*}
where $\sum_{i=1}^{n/2}e^i=e\in \z$ is the Euler number of the Seifert bundle $Y\rightarrow F$. Since $Y\rightarrow W_K$ is a $p$-fold cover,
\begin{equation*}
e=pe(W_K)=-p(\frac{q_1}{p}+\frac{q_2}{p}+\cdots+\frac{q_n}{p})=-(q_1+q_2+\cdots +q_n).
\end{equation*}
 $q_{j}$ is odd because $p$ is even and $gcd(p,q_j)=1$, $1\leq j\leq n$.
Then $e=-(q_1+q_2+\cdots +q_n)$ is an even number since $n$ is even.
For convenience, we may assume 
\begin{equation*}
e^i=2q, 1\leq i<\frac{n}{2}; \; e^{n/2}=e'=e-2q(n/2-1)=e-q(n-2),
\end{equation*}
where $q$ is a nonzero integer so that $e'\neq0$.
Note that $e'$ is even since $e$ and $n$ are both even.

$T^{i}$ is also fibered by geodesics isotopic to $L_{i,1,1}$. 
There is a {\it new} fibration of $Y_2^i$ corresponding to the new fibration of $T^{i}$.
Suppose $\bar{F}_2^i$ is the base space of the new fibration of $Y_2^i$. 
We call the fiber of $Y_{2}^i$ in this new fibration structure {\it new} fiber, and denoted $\bar{\phi}^i$.
Fix an orientation for the new circle fibers.
Note that the other component of $\hat{f}^{-1}(L_{i,1}^*)$, $L_{i,1,2}$, is parallel to $L_{i,1,1}$ so $L_{i,1,2}$ is also a new fiber.
We call the circle fiber of $M_2^i$ from the original fibration of $Y$ the {\it original fiber}, denoted $\phi^i$.
 Let $\bar{B}_2^i$ be one horizontal section of $Y_2^i \rightarrow \bar{F}_2^i$ such that $\bar{B}_2^i\cap T^i$ is a geodesic. 
 
 Let $N^i$ be a regular neighborhood of $L_{i,1,1}$ in $Y_2$. We may assume that $N^i$ is disjoint from other components of $L$. 
Note that $N^i$ is consists of new fibers.
 Let 
 \begin{equation*}
 M_2^i=Y_2-\overset{\circ}{N^i},  \ \ \ M=Y-\cup_{i=1}^{n/2}\overset{\circ}{N^i}=Y_1\cup(\cup_{i=1}^{n/2}M_2^i)
 \end{equation*}
 $M$ is the exterior of $\{L_{i,1,1}:1\leq i\leq n/2\}$ in $Y$. 
 $\p M_2^i=T^i_{2,-}\cup T^i_{2,+}\cup T^i_{2,0}$, where $T^i_{2,0}=\p N^i$.
 $M$ is a graph manifold with boundary $\cup_{i=1}^{n/2}T^i_{2,0}$ and characteristic tori $\{T^i_-\cup T^i_+:1\leq i\leq n/2\}$. 
Let $\bar{B}_2^{i,0}=\bar{B}^i_2\cap M_2^i$, which is an annulus with one puncture. Orient $\bar{B}_2^{i,0}$ and give the new fiber of $M_2^i$ a fixed orientation.
 
Define $\bar{\a}^i_{2,s}$ to be the component of $\p\bar{B}_2^{0,i}$ in $T^i_{2,s}$, and $\bar{\phi}^i_{2,s}$ to be a fixed new fiber in $T^i_{2,s}$ with the given orientation, $s=-,+,0$. 
Then $\{\bar{\a}^i_{2,\pm}, \bar{\phi}^i_{2,\pm}\}$ forms another basis for $H_1(T^i_{2,\pm})$, and $\{\bar{\a}^i_{2,0}, \bar{\phi}^i_{2,0}\}$ forms a basis for $H_1(T^i_{2,0})$.

As in the proof of Prop. 6.1 in \cite{abz}, we use the method introduced in \cite{wy} to construct essential horizontal surfaces $H_1$ in $Y_1$ and $H_2^i$ in $M_2^i$. 

Suppose that the old basis
 $\{\a_{2, -}^i, \phi_{2, -}^i\}$ of $H_1(T_{2, -}^i)$ and
 the new one $\{\bar{\a}_{2, -}^i, \bar{\phi}_{2, -}^i\}$
satisfies the following relation.
\begin{align*}
\bar{\a}_{2, -}^i&=a^i\a_{2, -}^i+b^i\phi_{2, -}^i, & \bar{\phi}_{2, -}^i&=c^i\a_{2, -}^i+d^i\phi_{2, -}^i,
\end{align*}
where $a^i, b^i, c^i, d^i$ are integers.
 We may assume that $a^id^i-b^ic^i=1$,
by reversing the orientation of the original fibers if necessary.
For convenience, we can suppose that $a^1=a^2= \cdots =a^{n/2}=a,\ b^1=b^2=\cdots =b^{n/2}=b,\ c^1=c^2=\cdots =c^{n/2}=c$, and $d^1=d^2=\cdots =d^{n/2}=d$.

In $M_2^i$, we have $\bar{\a}_{2, -}^i=-\bar{\a}_{2, +}^i$ and $\bar{\phi}_{2, -}^i=\bar{\phi}_{2, +}^i$. Thus
\begin{align*}
\bar{\a}_{2, +}^i&=a\a_{2, +}^i-b\phi_{2, +}^i, & \bar{\phi}_{2, +}^i&=-c\alpha_{2, +}^i+d\phi_{2, +}^i.
\end{align*}
The gluing maps $g_-^i: T_{2, -}^i\rightarrow T_{1, -}^i$ and $g_+^i: T_{2, +}^i\rightarrow T_{1, +}^i$ can be expressed as
\begin{align*}
(g_-^i)_*(\bar{\a}_{2, -}^i)&=-a\a_{1, -}^i+b\phi_{1, -}^i,& (g_-^i)_*(\bar{\phi}_{2, -}^i)&=-c\a_{1, -}^i+d\phi_{1, -}^i;\\
(g_+^i)_*(\bar{\a}_{2, +}^i)&=-a\a_{1, +}^i+(ae^i-b)\phi_{1, +}^i, & (g_+^i)_*(\bar{\phi}_{2, +}^i)&=c\a_{1, +}^i+(d-ce^i)\phi_{1, +}^i.
\end{align*}
Let $G_-^i$ and $G_+^i$ be the associated matrices corresponding to $g_-^i$ and $g_+^i$ respectively.
\begin{align*}
G_-^i=(g_-^i)_*&=\begin{pmatrix}-a&b\\-c&d\end{pmatrix} , & G_+^i=(g_+^i)_*&=\begin{pmatrix}-a&{2qa-b}\\c&{d-2qc}\end{pmatrix},\\
(G_-^i)^{-1}&=\begin{pmatrix}-d&b\\-c&a\end{pmatrix}, & 
(G_+^i)^{-1}&=\begin{pmatrix}{2qc-d}&{2qa-b}\\c&a\end{pmatrix},
\end{align*}
when $1\leq i< n/2$.
\begin{align*}
G_-^{n/2}=(g_-^{n/2})_*&=\begin{pmatrix}-a&b\\-c&d\end{pmatrix} , & G_+^{n/2}=(g_+^{n/2})_*&=\begin{pmatrix}-a&{ae'-b}\\c&{d-ce'}\end{pmatrix},\\
(G_-^{n/2})^{-1}&=\begin{pmatrix}-d&b\\-c&a\end{pmatrix}, & 
(G_+^{n/2})^{-1}&=\begin{pmatrix}{ce'-d}&{ae'-b}\\c&a\end{pmatrix}.
\end{align*}

$M$ is a graph manifold. The JSJ-decomposition of $M$ consists of
 $n/2+1$ vertices corresponding to
 $Y_1, M_2^1, M_2^2, \cdots, M_2^{n/2}$, and $n$ edges corresponding to $T_-^i$'s and $T_+^i$'s.
By \cite{wy}, every non-zero solution of the following equation gives a horizontal surface of $M$.

\begin{equation*}
(Y-Z)\left(\begin{array}{c}\l \\ \bar{\l}_1 \\  \vdots \\ \bar{\l}_{n/2}\end{array}\right)=\left(\begin{array}{c}0\\{\ast}\\ \vdots \\{\ast} \end{array}\right),\ \l, \bar{\l}_1, \cdots, \bar{\l}_{n/2}, \ast \in \mathbb{Z},
\end{equation*}
where  $Y$ and $Z$ are $(n/2+1)\times(n/2+1)$ matrices defined on page 450 of \cite{wy}.
The entries of $Y$, and $Z$ are decided by the gluing matrix.
 Let $Y=\{y_{j_1,j_2}\}, 1\leq j_1,j_2\leq n/2+1$.
According to \cite{wy}, we have $y_{1, j}=y_{j,
1}=\displaystyle{\frac{2}{c}}, 2\leq j\leq {n/2+1}$, and
other entries of $Y$ are all zeroes.

$Z=diag(z_1,z_2,\cdots,z_{n/2+1})$ is a diagonal matrix with
\begin{align*}
z_1&=\Bigl[\sum_{j=1}^{n/2-1}(\frac{d}{c}+\frac{2qc-d}{c})\Bigr]+(\frac{d}{c}+\frac{ce'-d}{c})\\
&=\Bigl[\sum_{j=1}^{n/2-1}\frac{2qc}{c}\Bigr]+e'=2q(\frac{n}{2}-1)+e'=q(n-2)+e'=e,\\
z_j&=\displaystyle{\frac{a}{c}-\frac{a}{c}}=0, 2\leq j\leq
{\frac{n}{2}+1}.
\end{align*}
According to the equation (1.6) of \cite{wy}, we have
\begin{equation*}
\begin{pmatrix}{-e}&{\displaystyle{\frac{2}{c}}}&{\displaystyle{\frac{2}{c}}}&{\cdots}&{\displaystyle{\frac{2}{c}}}\\
{\displaystyle{\frac{2}{c}}}&0&0&{\cdots}&0\\[2ex]
{\displaystyle{\frac{2}{c}}}&0&0&{\cdots}&0\\{\vdots}&{\vdots}&{\vdots}&{\vdots}&{\vdots}\\{\displaystyle{\frac{2}{c}}}&0&0&{\cdots}&0\end{pmatrix}_{(n/2+1)\times (n/2+1)}\begin{pmatrix}\l\\\bar{\l}_1\\\overline{\l}_2\\{\vdots}\\\overline{\l}_{n/2}\end{pmatrix}=\begin{pmatrix}0\\{\ast}\\{\ast}\\{\vdots}\\{\ast}\end{pmatrix}\\
\end{equation*}
\begin{align}\label{lambda}
\Ra \begin{cases}{-e\l +\displaystyle{\frac{2}{c}}\overline{\l}_1+\displaystyle{\frac{2}{c}}\overline{\l}_2+{\cdots}+\displaystyle{\frac{2}{c}}\overline{\l}_{n/2}=0}\\{\displaystyle{\frac{2}{c}}\l=\ast}\end{cases}
\end{align}
We may assume $\bar{\l}_j=qc\l$, $1\leq j< n/2$.

From (\ref{lambda}), we have
\begin{align*}
-e\l+\frac{2}{c}[(\frac{n}{2}-1)qc\l+\bar{\l}_{n/2}]=0\\
[-e+(n-2)q]\l+\frac{2}{c}\bar{\l}_{n/2}=0\\
-e'\l+\frac{2}{c}\bar{\l}_{n/2}=0
\end{align*}
Then 
\begin{equation}\label{solution}
\frac{\l}{\bar{\l}_{j}}=\frac{1}{cq}, \frac{\bar{\l}_{j}}{\l}=\frac{cq}{1}, 1\leq j<\frac{n}{2}; \;\;\frac{\l}{\bar{\l}_{n/2}}=\frac{2}{ce'}, \frac{\bar{\l}_{n/2}}{\l}=\frac{ce'}{2}.
\end{equation}

Suppose that 
\begin{equation*}
\p H_1=u_\pm^i\a_{1, \pm}^i+t_\pm^i\phi_{1, \pm}^i, 
\end{equation*}
on $T_\pm^i$ with respect to the basis $\{\a_{1, \pm}^i, \phi_{1, \pm}^i\}$, and 
\begin{equation*}
\p H_2^i= \bar{u}_s^i\bar{\a}_{2, s}^i+\bar{t}_s^i\bar{\phi}_{2, s}^i,
\end{equation*}
on $T_s^i$ with respect to the basis $\{\bar{\a}_{2, s}^i, \bar{\phi}_{2, s}^i\}$, $s=-,+,0$. Then there are $\e_-^i,\e_+^i\in\{\pm1\}$, so that equation (1.2) of \cite{wy} becomes:
\begin{align*}
\frac{t_-^i}{u_-^i}&=\frac{\e_-^i\bar{\l}_i}{-\l c}-\frac{d}{c}=-q\e_-^i-\frac{d}{c}, & 
\frac{\bar{t}_-^i}{\bar{u}_-^i}&=\frac{\e_-^i\l}{-\bar{\l}_ic}+\frac{a}{-c}=\frac{-\e_-^i-acq}{c^2q},\\
\frac{t_+^i}{u_+^i}&=\frac{\e_+^i\bar{\l}_i}{\l c}-\frac{2cq-d}{c}=q(\e_+^i-2)+\frac{d}{c},&
\frac{\bar{t}_+^i}{\bar{u}_+^i}&=\frac{\e_+^i\l}{\bar{\l}_ic}+\frac{a}{c}=\frac{\e_+^i+acq}{c^2q},
\end{align*}
when $1\leq i<n/2$.
\begin{align*}
\frac{t_-^{n/2}}{u_-^{n/2}}&=\frac{\e_-^{n/2}\bar{\l}_{n/2}}{-\l c}-\frac{d}{c}=\frac{-\e_-^{n/2}ce'-2d}{2c},&
 \frac{\bar{t}_-^{n/2}}{\bar{u}_-^{n/2}}&=\frac{\e_-^{n/2}\l}{-\bar{\l}_{n/2}c}+\frac{a}{-c}=\frac{-2\e_-^{n/2}-ace'}{c^2e'},\\
\frac{t_+^{n/2}}{u_+^{n/2}}&=\frac{\e_+^{n/2}\bar{\l}_{n/2}}{\l c}-\frac{ce'-d}{c}=\frac{2d-(2-\e_+^{n/2})ce'}{2c},&
\frac{\bar{t}_+^{n/2}}{\bar{u}_+^{n/2}}&=\frac{\e_+^{n/2}\l}{\bar{\l}_{n/2}c}+\frac{a}{c}=\frac{2\e_+^{n/2}+ace'}{c^2e'}.
\end{align*}
Since $H_1$ is a horizontal surface, $\sum_{i=1}^{n/2}(\displaystyle{\frac{t_-^i}{u_-^i}+\frac{t_+^i}{u_+^i}})=0$. 
Then
\begin{align*}
\sum_{i=1}^{\frac{n}{2}-1}(-q\e^i_--\frac{d}{c}+q(\e^i_+-2)+\frac{d}{c})+\frac{-\e_-^{n/2}ce'-2d+2d-2ce'+\e^{n/2}_+ce'}{2c}&=0\\
\sum_{i=1}^{\frac{n}{2}-1}(-\e^i_-+\e^i_+-2)q+\frac{(-\e_-^{n/2}+\e^{n/2}_+-2)ce'}{2c}&=0
\end{align*}

From the above equation, we may assume that $\e_-^i=-1, \e_+^i=1, 1\leq i\leq n/2$.
 
Then
\begin{equation*}
\frac{t_-^i}{u_-^i}=q-\frac{d}{c}, \frac{t_+^i}{u_+^i}=-q+\frac{d}{c}; \ \ 
\frac{\bar{t}_-^i}{\bar{u}_-^i}=\frac{1-acq}{c^2q},  \frac{\bar{t}_+^i}{\bar{u}_+^i}=\frac{1+acq}{c^2q}, 1\leq i<\frac{n}{2}.
\end{equation*}
\begin{equation*}
\frac{t_-^{n/2}}{u_-^{n/2}}=\frac{ce'-2d}{2c},
\frac{t_+^{n/2}}{u_+^{n/2}}=\frac{2d-ce'}{2c}, \ \ 
\frac{\bar{t}_-^{n/2}}{\bar{u}_-^{n/2}}=\frac{2-ace'}{c^2e'},
\frac{\bar{t}_+^{n/2}}{\bar{u}_+^{n/2}}=\frac{2+ace'}{c^2e'}.
\end{equation*}

$H_2$ is also a horizontal surface, so $\displaystyle{\frac{\bar{t}_-^i}{\bar{u}_-^i}+\frac{\bar{t}_+^i}{\bar{u}_+^i}+\frac{\bar{t}_0^i}{\bar{u}_0^i}}=0$. We have
\begin{equation*}
\begin{split}
\frac{\bar{t}_0^i}{\bar{u}_0^i}&=-(\frac{\bar{t}_-^i}{\bar{u}_-^i}+\frac{\bar{t}_+^i}{\bar{u}_+^i})=-(\frac{1-acq+1+acq}{c^2q})=-\frac{2}{c^2q}, 1\leq i<n/2.\\
\frac{\bar{t}_0^{n/2}}{\bar{u}_0^{n/2}}&=-(\frac{\bar{t}_-^{n/2}}{\bar{u}_-^{n/2}}+\frac{\bar{t}_+^{n/2}}{\bar{u}_+^{n/2}})=-(\frac{2-ace'+2+ace'}{c^2e'})=-\frac{4}{c^2e'}.
\end{split}
\end{equation*}

Recall that $L_{i,1,1}$ is a 1-fold cover of $L_{i, 1}^*$, so $c=1$. 
By (\ref{solution}), $\l=1, \bar{\l}_{i}=q, 1\leq i<{n/2}$. 
$\bar{\l}_{n/2}=\displaystyle{\frac{e'}{2}}$, is a set of solution of (\ref{lambda}). ($\bar{\l}_{n/2}$ is an integer since $e'$ is even.)
We can take $b=0, a=d=1$ by choosing proper horizontal sections. 
The slopes becomes
\begin{align*}
\frac{t_-^i}{u_-^i}&=q-1, &  \frac{t_+^i}{u_+^i}&=-q+1, &
\frac{\bar{t}_-^i}{\bar{u}_-^i}&=\frac{1-q}{q}, & \frac{\bar{t}_+^i}{\bar{u}_+^i}&=\frac{1+q}{q}, & \frac{\bar{t}_0^i}{\bar{u}_0^i}&=\frac{-2}{q},\\
\text{when }&1\leq i<n/2.\\
\frac{t_-^{n/2}}{u_-^{n/2}}&=\frac{e'-2}{2},&
\frac{t_+^{n/2}}{u_+^{n/2}}&=\frac{-e'+2}{2}, & 
\frac{\bar{t}_-^{n/2}}{\bar{u}_-^{n/2}}&=\frac{2-e'}{e'}, & 
\frac{\bar{t}_+^{n/2}}{\bar{u}_+^{n/2}}&=\frac{2+e'}{e'}, &
\frac{\bar{t}_0^{n/2}}{\bar{u}_0^{n/2}}&=-\frac{4}{e'}.
\end{align*}
We can determine $t_\pm^i, u_\pm^i, \bar{t}_\pm^i, \bar{u}_{\pm}^i, \bar{t}_0^i$ 
 and $\bar{u}_0^i$, as in Table \ref{bs}. 

\begin{table}
\begin{center}
\renewcommand{\arraystretch}{1.8}
\begin{tabular}[b]{|c|c||c|c|c|c|}
\hline
\multicolumn{2}{|c||}{\raisebox{1.5ex}{$\p H_1$}} & \multicolumn{4}{|c|}{\raisebox{1.5ex}{$\p H_2^i$}}\\[-1.8ex]
\hline \hline
$t_-^i=q-1$ & $ t_+^i=1-q$ & $\bar{t}_-^i=1-q$ & $ \bar{t}_+^i=1+q$ & $\bar{t}_0^i=-2,  2\nmid q$ & $\bar{t}_0^i=-1, 2| q$ \\
\hline 
$u_-^i=1$ & $u_+^i=1$ & $\bar{u}_-^i=q$ & $\bar{u}_+^i=q$ & 
$\bar{u}_0^i=q,  2\nmid q$ & $\bar{u}_0^i=q/2, 2| q$\\ \hline
 \hline
$t_-^\frac{n}{2}=\displaystyle{\frac{e'}{2}-1}$ & $ t_+^\frac{n}{2}=1-\displaystyle{\frac{e'}{2}}$ & $\bar{t}_-^\frac{n}{2}=1-\displaystyle{\frac{e'}{2}}$ & $ \bar{t}_+^\frac{n}{2}=1+\displaystyle{\frac{e'}{2}}$ & $\bar{t}_0^\frac{n}{2}=-2, 4\nmid e'$ & $\bar{t}_0^\frac{n}{2}=-1, 4| e'$ \\
\hline 
$u_-^\frac{n}{2}=1$ & $u_+^\frac{n}{2}=1$ & $\bar{u}_-^\frac{n}{2}=\displaystyle{\frac{e'}{2}}$ & $\bar{u}_+^\frac{n}{2}=\displaystyle{\frac{e'}{2}}$ & 
$\bar{u}_0^\frac{n}{2}=\displaystyle{\frac{e'}{2}},  4\nmid e'$ & $\bar{u}_0^\frac{n}{2}=\displaystyle{\frac{e'}{4}}, 4| e'$\\ \hline
\end{tabular}
\end{center}
\caption{\label{bs} The boundary slope of $H_1$, $H_2^{i}$ and $H_2^{n/2}$, where $1\leq i<n/2$.}
\end{table}

From \cite{wy}, $\epsilon_1^i=-\epsilon_2^i$ implies that $H=H_1\cup(\underset{i=1}{{\overset{n/2}{\cup}}}H_2^i)$ is non-orientable, so $M$ is a surface semi-bundle.
$\square$

Without further notice, we consider $1\leq i\leq n/2$ in the following discussion.

The construction of $H_2^i$ is the same as the construction of $H_2$ in Sec 6.1 of \cite{abz}. 
$H_2^i$ is a surface which interpolates between the slope $\displaystyle{\frac{1-q}{q}}$ on $T_-^i$ and $\displaystyle{\frac{1+q}{-q}}$ on $T_+^i$, $1\leq i<n/2$.
$H_2^{n/2}$ is a surface which interpolates between the slope $\displaystyle{\frac{2-e'}{e'}}$ on $T_-^{n/2}$ and $\displaystyle{\frac{2+e'}{-e'}}$ on $T_+^{n/2}$. 
$\p H_1\cap T_{1, \pm}^i$ has $\displaystyle{|\frac{\l}{u_\pm^i}|=|\frac{1}{1}|=1}$ component, and $\p H_2^i\cap T_{2, \pm}^i$ has $\displaystyle{|\frac{\bar{\l}_i}{\bar{u}_\pm^i}|=|\frac{q}{q}|=1}$ component, when $1\leq i<n/2$.
$\p H_2^{n/2}\cap T_{2,\pm}^{n/2}$ has $\displaystyle{|\frac{\bar{\l}_{n/2}}{\bar{u}_\pm^{n/2}}|=|\frac{e'/2}{e'/2}|=1}$ component. 
Denote the associated fibering in $M_2^i$ by ${\cal F}_2^i$.  
${\cal F}_2^i$ is transverse to all the new fibers, and in particular to $L_{i, 1,2}$.

The construction of $H_1$ is slightly different.  
Recall that $B_1$ is a horizontal section of $Y_1\rightarrow F_1$, such that $B_1\cap (T^i_{1,-}\cup T^i_{1, +})$ are geodesics. 
Note that $B_1$ is connected sine $F_1$ is, by Lemma \ref{connected}. 
Then there exists a properly embedded arc $\s^i$ in $B_1$ connecting $T^i_{1, -}$ and $T^i_{1, +}$. 
$\s^i\times[-1, 1]$ is a regular neighborhood of $\s^i$ in $B_1$.  
Wrap $\s^i\times[-1, 1]$ around the
$\phi^i$-direction $q-1$ times as we pass from $-1$ to $1$ in $\s^i\times[-1, 1]\times\phi^i$, when $1\leq i<n/2$; ($\displaystyle{1-\frac{e'}{2}}$) times when $i=n/2$.  
The resulting surface is $H_1$.
Let ${\cal F}_1$ be the corresponding surface fibration of $Y_1$
 with $H_1$ as a surface fiber.

We may suppose that $\p H_1=\p H_2^1\cup\p H_2^2 \cdots \cup\p H_2^{n/2}$ and that ${\cal F}_1\cup {\cal F}^1_2\cdots \cup{\cal F}_2^{n/2}$ forms a semi-surface bundle ${\cal F}$ in $M$, as described in Proposition \ref{semi}.

Next, we isotopy $L$ such that $L-\{L_{2i-1,1,1}: 1\leq i\leq n/2\}$ is transverse to $\{H_2^i: 1\leq i\leq n/2\}$. 
Recall that we reorient $\{L_{i,j}^*: 1\leq i\leq n/2, 1\leq j\leq m\}$  before the proof of Proposition \ref{semi}, and 
\begin{equation*} 
L_{i,j}^*=(L_{i,j}^1)^*\cup(L_{i,j}^2)^*\cup(L_{i,j}^3)^*\cup(L_{i,j}^4)^*
\end{equation*}
where $(i,j)\in A$.
\begin{align*}
\hat{c}_{2i-1}=\{(L_{2i-1,j}^2)^*:1<j\leq m\}&\cap
\{(L_{2i-2,j}^2)^*: 1\leq j\leq m\}, 1\leq i\leq n/2;\\
\hat{c}_{2i}=\{(L_{2i-1,j}^4)^*:1<j\leq m\}&\cap\{(L_{2i,j}^4)^*: 1\leq j\leq m\}, 1\leq i\leq n/2;\\
(L_{2i-1,j}^2)^*\cup(L_{2i-1,j}^4)^*&\subset F_2^i,1\leq i\leq n/2, 1< j\leq m;\\
(L_{2i-2,j}^2)^* \cup(L_{2i,j}^4)^*&\subset F_2^i, 1\leq i\leq n/2, 1\leq j\leq m;\\
(L_{i,j}^1)^*\cup (L_{i,j}^3)^*&\subset F_1, (i,j)\in A. 
\end{align*}
Correspondingly, $L_{i,j,1}$ and $L_{i,j,2}$ are also split into four parts, $(i,j)\in A$. 
We denote them $L_{i,j,t}^l$, $(i,j)\in A, t=1,2, l=1,2,3,4$.
$\hat{f}(L_{i,j,1}^l\cup L_{i,j,2}^l)=(L_{i,j}^l)^*$, $(i,j)\in A, l=1,2,3,4.$

$U_{i,j}^l=\hat{f}^{-1}((L_{i,j}^l)^*)$ is a vertical annulus in $Y$, and $L_{i,j,1}^{l}\cup L_{i,j,2}^{l}\subset U_{i,j}^l$, where $(i,j)\in A, l=1,2,3,4.$
 Suppose 
\begin{equation*}
U_{i,j}^{l,0}=U_{i,j}^l \cap M, (i,j)\in A, l=1,2,3,4.
\end{equation*}
Note that $ \cup_{l=1}^{4}(U_{i,j}^{l,0})=\hat{f}^{-1}(L^*_{1,j})\cap M$ is an $2$-punctured torus, $(i,j)\in A$. 
By the construction of $U_{i,j}^{l,0}$, we have 
\begin{equation}\label{Ul}
U_{2i-1,j}^{l,0}\subset
\begin{cases}
M_2^i      & \text{if $l$ is even} , \\
 Y_1     & \text{otherwise},
\end{cases}
 1<j\leq m;\;\;
U_{2i,j}^{l,0}\subset
\begin{cases}
M_2^{i+1}      & \text{if }l=2 , \\
M_2^i            & \text{if }l=4,\\
 Y_1     & \text{otherwise},
\end{cases}
1\leq j\leq m,
\end{equation}
where $1\leq i\leq n/2$.

Different from Case I, $L_{i,1,1}$ only intersects each original fiber once, so $U_{i, j}^{2l,0}$ is a once-punctured annulus, and there is only one singular point on $U_{i, j}^{2l,0}, (i,j)\in A, l=1,2$.
Let $\hat{\phi}_{i}$ be the original fiber at $\hat{c}_i$, $1\leq i\leq n$.
The singular point on $U_{2i-1,j}^{2,0}$ is contained in $\hat{\phi}_{2i-1}$ and the singular point on $U_{2i-1,j}^{4,0}$ is contained in $\hat{\phi}_{2i}$, $1\leq i\leq n/2, 1<j\leq m$. 
The singular point on $U_{2i,j}^{2,0}$ is contained in $\hat{\phi}_{2i+1}$ and the singular point on $U_{2i,j}^{4,0}$ is contained in $\hat{\phi}_{2i}$, $1\leq i\leq n/2, 1\leq j\leq m$. 
$U_{i,j}^{2l,0}$ is similar as in Figure \ref{foliation2}.

Similar to Prop. \ref{M2}, we have the following proposition.

\begin{prop}\label{M2i}
Give ${\cal F}_2^i$ a fixed transverse orientation in $M_2^i$. 
We can isotope $L_{i, j,t}$ along the original fibers in $U_{i,j,t}^{2l,0}$, such that $L_{i,j,t}$ travels from the negative to the positive side of ${\cal{F}}_2^i$'s leaves.
This isotopy fixes outside a small regular neighborhood of $U_{i,j,t}^{2l,0}$.
$(i,j)\in A, t, l=1,2$.
\end{prop}

Different from \cite{abz}, not all arcs of $M_2^i\cap \{L_{i,j,t}: (i,j)\in A, t=1,2\}$ are from $T^i_{2,-}$ to $T^i_{2,+}$. 
We need the following claim.

\begin{claim}\label{isotopy}
If every two segments $(L_{i_1, j_1}^{2l_1})^*$ and $(L_{i_2, j_2}^{2l_2})^*$ with $(L_{i_1, j_1}^{2l_1})^*\cap(L_{i_2, j_2}^{2l_2})^*\neq \emptyset$ intersect $\{L_{2i-1,1}^*: 1\leq i\leq n/2\}$ in the same direction, then Prop \ref{M2i} is true, $(i_1,j_1), (i_2,j_2)\in A, l_1,l_2=1,2$.
\end{claim}
{\bf Proof}: 
Assume $(i_1,j_1), (i_2,j_2)\in A, l_1,l_2=1,2$ in this proof.
Suppose that $(L_{i_1, j_1}^{2l_1})^*\cap(L_{i_2, j_2}^{2l_2})^*=\hat{c}_{i_0}$, $1\leq i_0\leq n/2$.
Then $U_{i_1,j_1}^{2l_1,0}\cap U_{i_2,j_2}^{2l_2,0}=\hat{\phi}_{i_0}$.
 $U_{i_1,j_1}^{2l_1,0}$ ($U_{i_2,j_2}^{2l_2,0}$) is a once-punctured annulus, and there is one singular point contained in $\hat{\phi}_{i_0}\subset U_{i_1,j_1}^{2l_1,0}$ ($U_{i_2,j_2}^{2l_2,0}$).
By the construction of $F$, $(L_{i_1, j_1}^{2l_1})^*$ and $(L_{i_2, j_2}^{2l_2})^*$ intersect  $L^*_{i_0,1}$ when $i_0$ is odd, or $L^*_{i_0-1, 1}$ if $i_0$ is even.
Without losing of generality, we may suppose that $i_0$ is odd.
By assumption, $(L_{i_1, j_1}^{2l_1})^*$ and $(L_{i_2, j_2}^{2l_2})^*$ intersect  $L^*_{i_0,1}$ in the same direction.
We have that the heads of $\{L_{i_1,j_1,t}^{2l_1}: t=1,2\}$ and $\{L_{i_2,j_2,t}^{2l_2}: t=1,2\}$ are all contained in $T_{2,-}^{(i_0+1)/2}$ or $T_{2,+}^{(i_0+1)/2}$.
By the proof of Prop. 6.2 in \cite{abz}, we can isotope $\{L_{i_1,j_1,t}^{2l_1}: t=1,2\}$ and $\{L_{i_2,j_2,t}^{2l_2}: t=1,2\}$ along the original fibers above or below all the singular points in $\hat{\phi}_{i_0}$ such that they are all travel from the $``-"$ to the $``+"$ side of the leaves of ${\cal F}_2^{(i_0+1)/2}$.
This isotopy fixes outside a small regular neighborhood of 
$U_{i_1,j_1}^{2l_1,0}$ and $U_{i_2,j_2}^{2l_2,0}$.

Consider the case $(L_{i_1, j_1}^{2l_1})^*\cap(L_{i_2, j_2}^{2l_2})^*=\emptyset$. $\{L_{i_1,j_1,t}^{2l_1}: t=1,2\}$ and $\{L_{i_2,j_2,t}^{2l_2}: t=1,2\}$ won't block each other during the isotopies, since $U_{i_1,j_1}^{2l_1,0}\cap U_{i_2,j_2}^{2l_2,0}=\emptyset$.
Then Prop. \ref{M2i} is true. $\square$

By the construction of $F$, $(L_{i_1, j_1}^{2l_1})^*\cap(L_{i_2, j_2}^{2l_2})^*\neq\emptyset$ when $l_1=l_2=1$ and $i_1=i_2+1$, or  $l_1=l_2=2$ and $i_1=i_2-1$, $(i_1, j_1), (i_2, j_2)\in A$. 
According to  Remark \ref{arc}, $\{L_{i,j}^*:(i,j)\in A\}$ satisfies the assumption of Claim \ref{isotopy}, then Prop. \ref{M2i} is true.

\begin{figure}
\begin{center}
\includegraphics{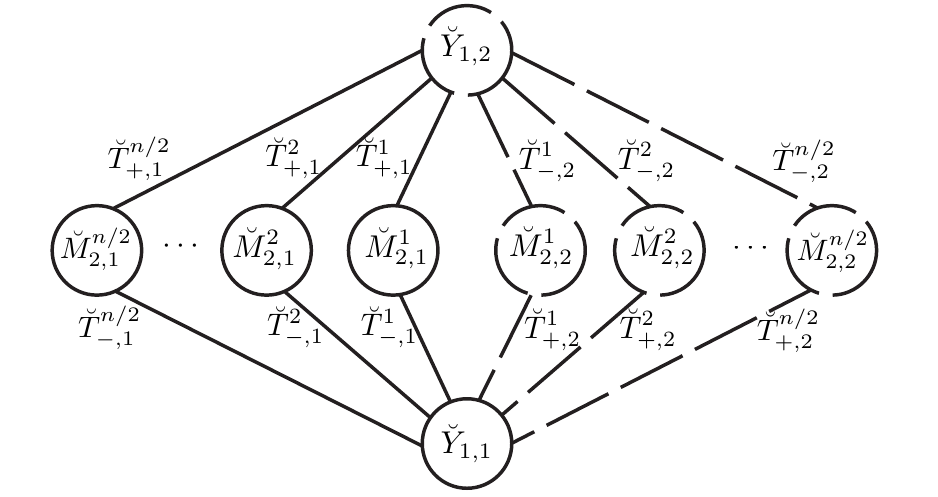}
\end{center}
\caption{\label{double2} Construction of $\breve{M}$.}
\end{figure}

As in Section 2, we construct, a free 2-fold cover of $M$, $\breve{M}$, which is a fibered graph manifold.
The construction is shown in Figure \ref{double2}, (c.f. Figure 8 in Sec. 6.1 of \cite{abz}).
Let $p_2$ be the covering map. 
\begin{align*}
\breve{Y}_{1, 1}\cup \breve{Y}_{1, 2}=p_2^{-1}(Y_1),\; 
\breve{M}_{2, 1}^i\cup \breve{M}_{2, 2}^i=p_2^{-1}(M_2^i),\;
\breve{T}_{\pm, 1}^i\cup \breve{T}_{\pm, 2}^i=p_2^{-1}(T_\pm^i),\; 1\leq i\leq n/2.
\end{align*}
$\breve{M}$ is a graph manifold obtained by gluing $\breve{M}_{2,1}^i$ and $\breve{Y}_{1,2}$ together along $\breve{T}_{+,1}^i$, and gluing $\breve{M}_{2,2}^i$ and $\breve{Y}_{1,1}$ together along  $\breve{T}_{+,2}^i$,  $1\leq i\leq n/2$. 
The gluing map is the same as the one used to glue $M_2^i$ back to $Y_1$ along $T_{+}^i$ to get $M$ for every $1\leq i\leq n/2$. 
Figure \ref{double2} also illustrates the JSJ-decomposition of $\breve{M}$.
Let $H=H_1\cup(\overset{n/2}{\underset{i=1}{\cup}}H_2^i)$ and $\breve{H}=p_2^{-1}(H)$. $\breve{H}$ is a connected orientable surface in $\breve{M}$, so $\breve{M}$ is a fibered graph manifold.

$p_2$ can be extended to $Y$. 
Denote $\breve{Y}$ to be the covering space of $Y$. 
$\breve{Y}$ inherits the Seifert fibered structure from $Y$ with fiber $\breve{\phi}$, where $\breve{\phi}$ is the preimage of the original fiber $\phi$ in $Y$ under $p_2$.
Let $\breve{L}=p_2^{-1}(L)$.
As in \cite{abz}, we use Prop. \ref{M2i} and perform Dehn twist operations to the surface fibers of $\breve{M}$ such that $\breve{L}\cap\breve{M}$ always transverses the new surface fiber.

Let 
\begin{equation*}
\breve{\cal F}=\breve{\cal F}_{1, 1}\cup \breve{\cal F}_{1, 2}\cup (\underset{s=1}{\overset{2}{\cup}}\underset{i=1}{\overset{n/2}{\cup}}\breve{\cal F}_{2, s}^i),
\end{equation*}
 where $\breve{\cal F}_{1, s}=p_2^{-1}({\cal F}_1)\cap \breve{Y}_{1, s}$, $\breve{\cal F}_{2, s}^i=p_2^{-1}({\cal F}_2^i)\cap \breve{M}_{2, s}^i, s=1, 2, 1\leq i\leq n/2$.
Then $\breve{\cal F}=p_2^{-1}({\cal F})$ is a surface bundle in $\breve{M}$, and $\breve{H}$ represents one leaf of $\breve{F}$.  
Fix a transverse orientation for $\breve{\cal F}$ and let $\breve{\cal F}_{1, s}$ and $\breve{\cal F}_{2, s}^i$ have the induced orientation, $s=1, 2$, $1\leq i\leq n/2$.
\begin{equation*}
\breve{Y}=\breve{Y}_{1, 1}\cup \breve{Y}_{1, 2}\cup (\underset{s=1}{\overset{2}{\cup}}\underset{i=1}{\overset{n/2}{\cup}}\breve{Y}_{2, s}^i)
\end{equation*}
 where $\breve{Y}^i_{2, 1}\cup\breve{Y}_{2,2}^i=p_2^{-1}(Y_2^i)$ and $\breve{M}_{2, s}^i\subset\breve{Y}_{2,s}^i$, $s=1, 2$, $1\leq i\leq n/2$. 
By construction, $p_2^{-1}(L_{2i-1,1,1})=\breve{L}_{2i-1,1,1,1}\cup\breve{L}_{2i-1,1,1,2}$ are two copies of $L_{2i-1,1,1}$, where $\breve{L}_{2i-1,1,1,s}\subset{\breve{Y}_{2,s}^i}$, $1\leq i\leq n/2$, $s=1,2$. Equip $\breve{L}_{2i-1,1,1,s}$ the inherited orientation from $L_{2i-1,1,1}$, $1\leq i\leq n/2, s=1,2$.
Then 
 \begin{equation*}
 \breve{M}=\breve{Y}\setminus (\underset{s=1}{\overset{2}{\cup}}\underset{i=1}{\overset{n/2}{\cup}} N(\breve{L}_{2i-1,1,1,s})),
\end{equation*}
where $N(\breve{L}_{2i-1,1,1,s})$ is a regular  small neighborhood of $\breve{L}_{2i-1,1,1,s}$ in $\breve{Y}$, $1\leq i\leq n/2, s=1,2.$

By the construction of $\breve{M}$, $p_2^{-1}(U_{i, j}^{l,0})$ is two copies of $U_{i,j}^{l,0}$, denoted by $\breve{U}_{i, j, 1}^{l,0}\cup \breve{U}_{i, j, 2}^{l,0}$, $1\leq l\leq 4$, $(i,j)\in A$.
By (\ref{Ul})
\begin{align*}\label{Ul}
\breve{U}_{2i-1,j,s}^{l,0}&\subset
\begin{cases}
\breve{M}_{2,s}^i      & \text{if $l$ is even} , \\
 \breve{Y}_{1,s}     & \text{otherwise},
\end{cases}
 1<j\leq m; \;\;
\breve{U}_{2i,j,s}^{l,0}&\subset
\begin{cases}
\breve{M}_{2,s}^{i+1}      & \text{if }l=2 , \\
\breve{M}_{2,s}^i            & \text{if }l=4,\\
 \breve{Y}_{1,s}     & \text{otherwise},
\end{cases}
1\leq j\leq m,
\end{align*}
where $1\leq i\leq n/2, s=1,2$.
Let $\breve{L}_{i,j,t,1}^l\cup\breve{L}_{i,j,t,2}^l$ be the lift of $L_{i,j,t}^l\subset U_{i,j}^{l,0}$, where $\breve{L}_{i,j,t,s}^l\subset\breve{U}_{i,j,s}^{l,0}$, and equip $\breve{L}_{i,j,t,s}^{l}$ the inherited orientation, $(i,j)\in A, l=1,2,3,4, s,t=1,2$. 

By Proposition \ref{M2i}, we can isotope $\breve{L}$ along $\breve{\phi}$ in $\{\breve{U}_{i,j,s}^{2l,0}: (i,j) \in A, l,s=1,2\}$, such that $\{\breve{L}_{2i-1,j,t,s}^{2l}: 1<j\leq m; l,s,t=1,2\}\cup\{\breve{L}_{2i-2l,j,t,s}^{4-2l}: 1\leq j\leq m; s,t=1,2; l=0,1\}$ travel from the negative to the positive side of $\breve{\cal{F}}_{2,s}^i$'s leaves along the orientation on them, $1\leq i\leq n/2$.
This isotopy fixes outside a small regular neighborhood of $\{\breve{U}_{i,j,s}^{2l,0}: (i,j) \in A, l,s=1,2\}$. 
Note that $\breve{L}_{2i-1,1,2,s}=p_2^{-1}(L_{2i-1,1,2})\cap \breve{M}_{2,s}$ is parallel to $\breve{\phi}$, so it is transverse to $\breve{\cal{F}}_{2,s}^i$. $1\leq i\leq n/2$.

Next, we consider the arcs $\{\breve{L}_{i,j,t,s}^{2l-1}: (i,j)\in A, s,t,l=1,2\}=\breve{L}\cap (\breve{Y}_{1,1}\cup\breve{Y}_{1,2})$.
As in Section 2, we can perform Dehn twist operations a sufficiently large number of times on $\{\breve{\cal F}_{1, s} : s=1, 2\}$ along a set of mutually disjoint $\breve{\phi}$-vertical tori $\G\subset(\breve{Y}_{1,1}\cup\breve{Y}_{1,2})$, such that the new surface fibers are transverse to $\{\breve{L}_{i,j,t,s}^{2l-1}: (i,j)\in A, s,t,l=1,2\}$ everywhere. 
$p_2(\G)$ is a set of $\phi$-vertical tori in $Y_1$.
Suppose ${\cal C}$ is the Seifert quotient of $p_2(\G)$.
Then ${\cal C}$ is a set of mutually disjoint simple closed curves in $F_1$.
By Claim \ref{torus}, we need to show that $i((L_{i,j}^{2l-1})^*,{\cal C})$ is negative (or positive) for all   $(i,j)\in A, l=1,2$.

${\cal C}$ is constructed in two different ways when $4\mid n$ and $4\nmid n$.

\begin{figure}
\begin{center}
\includegraphics[width=15cm]{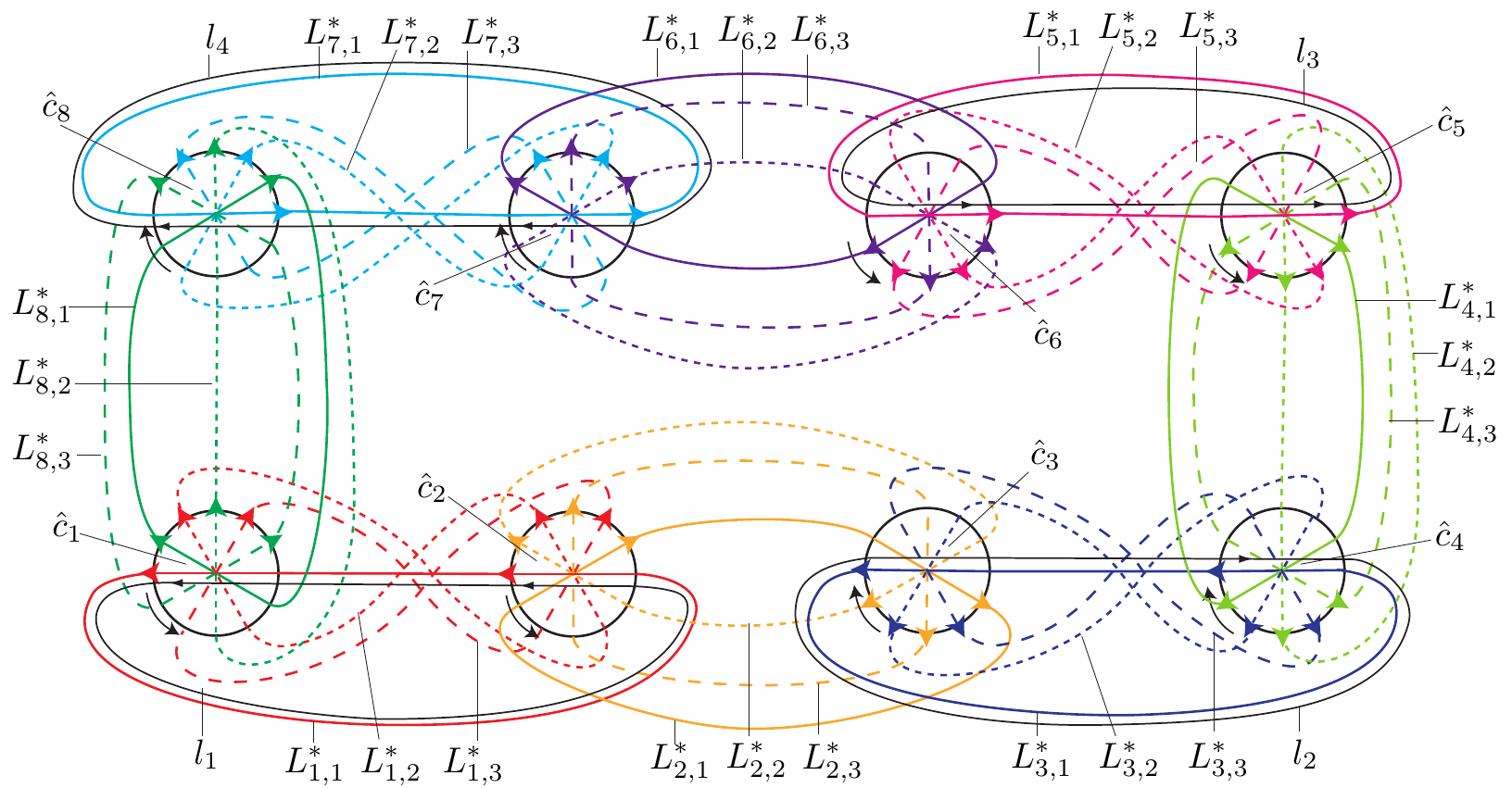}
\end{center}
\caption{\label{Lijstar8l} $\{l_i:1\leq i\leq n/2\}$ when $p=6$, $n=8$.}
\end{figure}

At first we consider $4\mid n$. 
In this case ${\cal C}$ consists $n/2$ simple closed curves in $F_1$.
${\cal C}=\{l_i:1\leq i\leq n/2\}$.
\begin{equation*}
l_i \text{ is in a small regular neighborhood of }  
\begin{cases}
\b_-^i \text{ in } F_1 \text{ and parallel to } \b_-^i & \text{ if}\ i \text{ is odd},\\
\b_+^i \text{ in } F_1 \text{ and parallel to } \b_+^i & \text{ if}\ i \text{ is even}.\\
\end{cases}
\end{equation*}
Then $l_i$ is parallel to $L^*_{2i-1,1}$ in $F$, as shown in Figure \ref{Lijstar8l}, when $p=6,n=8$, $1\leq i\leq n/2$.
Give $l_i$ the same orientation as  $L^*_{2i-1,1}$ if $i$ is odd, the opposite orientation if $i$ is even, $1\leq i\leq n/2$. 

Since $4\mid n$, $h_2(y_{4i-3})=h_2(y_{4i-2})=\bar{1}$ and $h_2(y_{4i-1})=h_2(y_{4i})=-\bar{1}$ by (\ref{h2}), $1\leq i\leq n/4$.

At first, consider that $i$ is odd.
In this case, $2i-2\equiv0$ (mod 4), $2i-1\equiv1$ (mod 4), $2i\equiv 2$ (mod 4), and $2i+1\equiv 3$ (mod 4), $1\leq i\leq n/2$. 
Then $h_2(y_{2i-1})=h_2(y_{2i})=\bar{1}$, $h_2(y_{2i-2})=h_2(y_{2i+1})=-\bar{1}$.
By Remark \ref{arc}, the arcs intersect $\b_-^i$ are shown as the following.
\begin{equation*}
\begin{cases}
(L^1_{2i-1,j})^* \text{ and } (L^3_{2i-1,j})^* \text{ are from } \b_+^i \text{ to }\b_-^i,1<j\leq m;\\ 
(L^1_{2i-2,j})^* \text{ is from } \b_-^{i-1} \text{ to } \b_-^i, 1\leq j\leq m;\\ 
(L^3_{2i,j})^* \text{ is from } \b_-^{i+1} \text{ to } \b_-^i, 1\leq j\leq m.
\end{cases}
\end{equation*}
When $i$ is odd, $l_i$ is a copy of $\b_{-}^i$. Then $l_i$ intersects each of $(L^1_{2i-1,j_1})^*$, $(L^3_{2i-1,j_1})^*$, $(L^1_{2i-2,j_2})^*$ and $(L^3_{2i,j_2})^*$ once, $1<j_1\leq m, 1\leq j_2\leq m$.
Because $l_i$ has the same orientation as $L^*_{2i-1,1}$, we have the following result by (\ref{neworientation}). 
\begin{align*}
i((L^{1}_{2i-1,j_1})^*,l_i)&=-i( L^*_{2i-1,1},L_{2i-1,j_1}^*)=-sign(h_2(y_{2i-1}))=-1,\\
i((L^{3}_{2i-1,j_1})^*,l_i)&=-i( L^*_{2i-1,1},L_{2i-1,j_1}^*)=-sign(h_2(y_{2i}))=-1,\\
i(L^{1}_{2i-2,j_2})^*,l_i)&=-i( L^*_{2i-1,1},L_{2i-2,j_2}^*)=-sign(h_2(y_{2i-1}))=-1,\\
i((L^{3}_{2i,j_2})^*,l_i)&=-i( L^*_{2i-1,1},L_{2i,j_2}^*)=-sign(h_2(y_{2i}))=-1,
\end{align*}
where $1\leq i\leq n/2$, and $i$ is odd, $1<j_1\leq m, 1\leq j_2\leq m$.

When $i$ is even,
 $2i-2\equiv 2$ (mod 4), $2i-1\equiv 3$ (mod 4),  $2i\equiv 0$ (mod 4), and $2i+1\equiv 1$ (mod 4), $1\leq i\leq n/2$. 
Then $h_2(y_{2i-1})=h_2(y_{2i})=-\bar{1}$, $h_2(y_{2i-2})=h_2(y_{2i+1})=\bar{1}$.
By Remark \ref{arc}, the arcs intersect $\b_+^i$ are shown as the following.
\begin{equation*}
\begin{cases}
(L^1_{2i-1,j})^* \text{ and } (L^3_{2i-1,j})^* \text{ are from } \b_-^i \text{ to }\b_+^i,1<j\leq m;\\ 
(L^1_{2i-2,j})^* \text{ is from } \b_+^{i-1} \text{ to } \b_+^i, 1\leq j\leq m;\\ 
(L^3_{2i,j})^* \text{ is from } \b_+^{i+1} \text{ to } \b_+^i, 1\leq j\leq m.
\end{cases}
\end{equation*}
When $i$ is evev, $l_i$ is a copy of $\b_{+}^i$. 
Then $l_i$ intersects each of $(L^1_{2i-1,j_1})^*$, $(L^3_{2i-1,j_1})^*$, $(L^1_{2i-2,j_2})^*$ and $(L^3_{2i,j_2})^*$ once, $1<j_1\leq m, 1\leq j_2\leq m$.
Because $l_i$ has the different orientation from $L^*_{2i-1,1}$, we have the following result by (\ref{neworientation}). 
\begin{align*}
i((L^{1}_{2i-1,j_1})^*,l_i)&=i( L^*_{2i-1,1},L_{2i-1,j_1}^*)=sign(h_2(y_{2i-1}))=-1,\\
i((L^{3}_{2i-1,j_1})^*,l_i)&=i( L^*_{2i-1,1},L_{2i-1,j_1}^*)=sign(h_2(y_{2i}))=-1,\\
i((L^{1}_{2i-2,j_2})^*,l_i)&=i( L^*_{2i-1,1},L_{2i-2,j_2}^*)=sign(h_2(y_{2i-1}))=-1,\\
i((L^{3}_{2i,j_2})^*,l_i)&=i( L^*_{2i-1,1},L_{2i,j_2}^*)=sign(h_2(y_{2i}))=-1,
\end{align*}
where $1\leq i\leq n/2$, and $i$ is even, $1<j_1\leq m, 1\leq j_2\leq m$.

By the above discussion, when $4\mid n$, $i((L^{2l-1}_{2i-1,j_1})^*,{\cal C})=i((L^{2l-1}_{2i,j_2})^*, {\cal C})=-1, 1<j_1\leq m, 1\leq j_2\leq m, 1\leq i\leq n/2$, $l=1,2$.

Now we consider the case $4\nmid n$. $n$ is even, so $4\mid (n+2)$. 
Let $\displaystyle{I=\frac{n+2}{4}}$.
In this case, ${\cal C}$ consists $n/2+3$ mutually disjoint simple closed curves in $F_1$. Let
\begin{equation*}
{\cal C}=\{l_i: 1\leq i\leq n/2\}\cup\{l_{j,1}: j=I-1,I,I+1\}.
\end{equation*}
Since $4\nmid n$, we have the following by (\ref{h2}).
\begin{equation}\label{hn4}
\begin{split}
\text{ When } 1\leq i\leq n/2,
h_2(y_{i})=
\begin{cases}
\bar{1} &\text{ if } i\equiv 1 \text{ or } 2 \text{ (mod) } 4,\\
-\bar{1} & \text{ otherwise}.
\end{cases}\\
\text{ When } n/2< i\leq n,
h_2(y_{i})=
\begin{cases}
-\bar{1} &\text{ if } i\equiv 1 \text{ or } 2 \text{ (mod) } 4,\\
\bar{1} & \text{ otherwise}.
\end{cases}
\end{split}
\end{equation}

\begin{figure}
\begin{center}
\includegraphics[width=15cm]{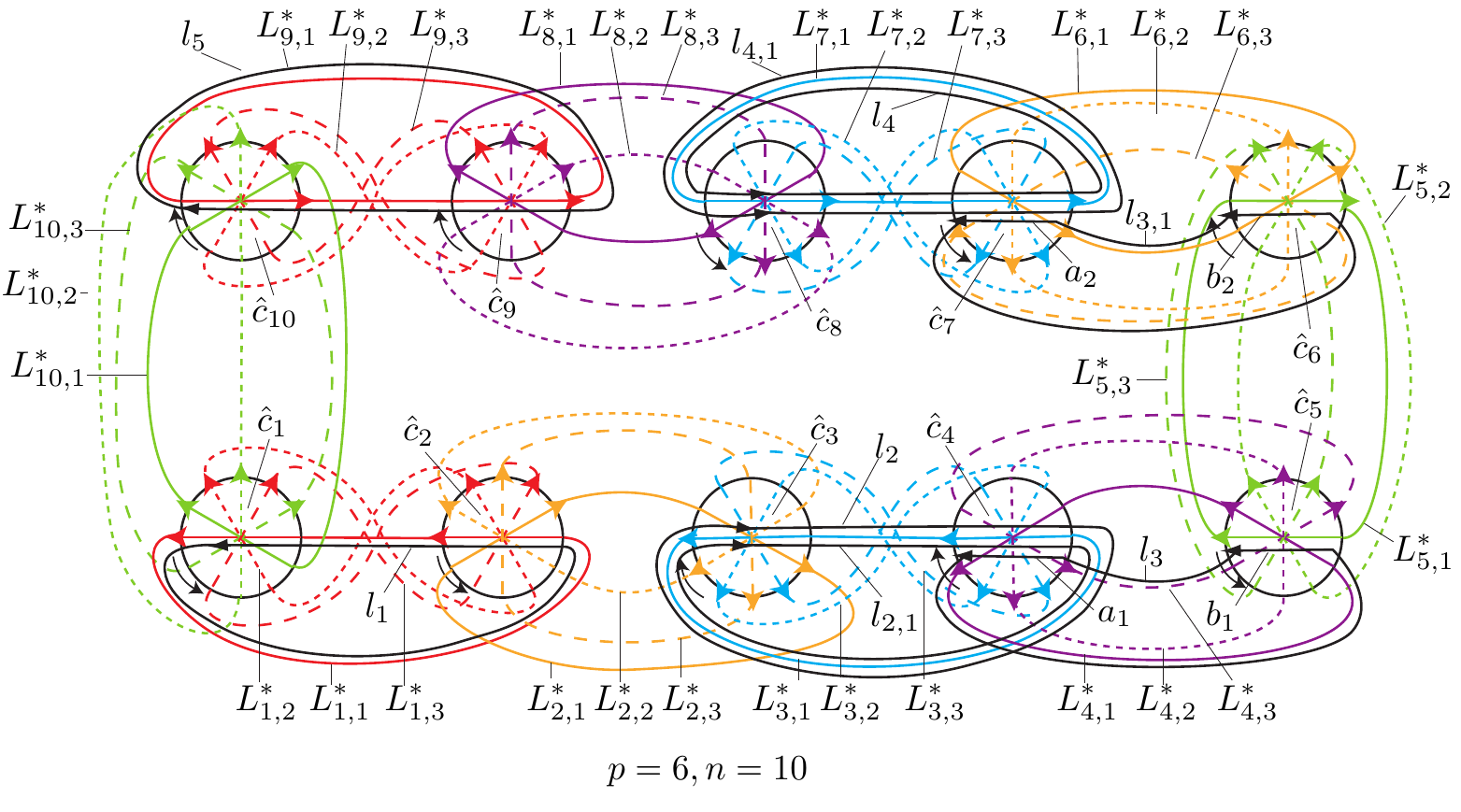}
\end{center}
\caption{\label{610} $\{l_i:1\leq i\leq n/2\}\cup\{l_{i,1}:i=I-1,I,I+1\}$ when $p=6$, $n=10$.}
\end{figure}

We construct ${\cal C}$ in the following cases. (c.f. Figure \ref{Lijstar6c} and Figure \ref{610}).

Case 1: Construct $\{l_i:1\leq i<I$ and $i$ is odd, or $I<i\leq n/2$ and $i$ is even$\}$.

Let $l_i$ be a simple closed curve parallel to $\b_-^i$ and in a small regular neighborhood of $\b_-^i$ in $F_1$.  
Give $l_i$ the same orientation as $L_{2i-1,1}^*$.
By (\ref{hn4}),
$h_2(y_{2i-1})=h_2(y_{2i})=\bar{1}$, $h_2(y_{2i-2})=h_2(y_{2i+1})=-\bar{1}$.
By a similar discussion as in the case $4\mid n$ and $i$ is odd,  $l_i$ intersects each of $(L^1_{2i-1,j_1})^*$, $(L^3_{2i-1,j_1})^*$, $(L^1_{2i-2,j_2})^*$ and $(L^3_{2i,j_2})^*$ once,  and 
\begin{equation}\label{c1}
i((L_{2i-1,j_1}^1)^*,l_i)=i((L^3_{2i-1,j_1})^*,l_i)=i((L^1_{2i-2,j_2})^*,l_i)=i((L^3_{2i,j_2})^*,l_i)=-1,
\end{equation}
$1<j_1\leq m, 1\leq j_2\leq m$.

Case 2: Construct $\{l_i:1\leq i<I$ and $i$ is even, or $I<i\leq n/2$ and $i$ is odd$\}$.

Let $l_i$ be a simple closed curve parallel to $\b_+^i$ and in a small regular neighborhood of $\b_+^i$ in $F_1$. 
Give $l_i$ the opposite orientation of $L_{2i-1,1}^*$.
By (\ref{hn4}),
$h_2(y_{2i-1})=h_2(y_{2i})=-\bar{1}$,  $h_2(y_{2i-2})=h_2(y_{2i+1})=\bar{1}$.
By a similar discussion as in the case $4\mid n$ and $i$ is even,  $l_i$ intersects each of $(L^1_{2i-1,j_1})^*$, $(L^3_{2i-1,j_1})^*$, $(L^1_{2i-2,j_2})^*$ and $(L^3_{2i,j_2})^*$ once,  and 
\begin{equation}\label{c2}
i((L_{2i-1,j_1}^1)^*,l_i)=i((L^3_{2i-1,j_1})^*,l_i)=i((L^1_{2i-2,j_2})^*,l_i)=i((L^3_{2i,j_2})^*,l_i)=-1,
\end{equation}
$1<j_1\leq m, 1\leq j_2\leq m$.

Case 3: Construct $\{l_{j,1}$: $j=I-1$ when $I$ is even, or $j=I+1$ when $I$ is odd$\}$. 

Let $l_{j,1}$ be a simple closed curve parallel to $\b_+^j$ and in a small regular neighborhood of $\b_+^j$ in $F_1$. 
Equip $l_{j,1}$ the same orientation as $L_{2j-1,1}^*$.
By (\ref{hn4}),
$h_2(y_{2j-1})=h_2(y_{2j})=\bar{1}$,  $h_2(y_{2j-2})=h_2(y_{2j+1})=-\bar{1}$.
By Remark \ref{arc}, 
\begin{equation*}
\begin{cases}
(L^1_{2j-1,j_1})^* \text{ and } (L^3_{2j-1,j_1})^* \text{ are from } \b_+^{j} \text{ to }\b_-^{j},1<j_1\leq m;\\ 
(L^3_{2j-2,j_2})^* \text{ is from } \b_+^{j} \text{ to } \b_+^{j-1}, 1\leq j_2\leq m;\\ 
(L^1_{2j,j_2})^* \text{ is from } \b_+^{j} \text{ to } \b_+^{j+1}, 1\leq j_2\leq m.
\end{cases}
\end{equation*}
Then $l_{j,1}$ intersects each of $(L^1_{2j-1,j_1})^*$, $(L^3_{2j-1,j_1})^*$, $(L^3_{2j-2,j_2})^*$ and $(L^1_{2j,j_2})^*$ once, $1<j_1\leq m, 1\leq j_2\leq m$.
By  (\ref{neworientation}),
\begin{equation}\label{c3}
\begin{split}
i((L^{1}_{2j-1,j_1})^*,l_{j,1})&=-i( L^*_{2j-1,1},L_{2j-1,j_1}^*)=-sign(h_2(y_{2j}))=-1,\\
i((L^{3}_{2j-1,j_1})^*,l_{j,1})&=-i( L^*_{2j-1,1},L_{2j-1,j_1}^*)=-sign(h_2(y_{2j-1}))=-1,\\
i((L^{3}_{2j-2,j_2})^*, l_{j,1})&=-i( L^*_{2j-1,1},L_{2j-2,j_2}^*)=-sign(h_2(y_{2j-1}))=-1,\\
i((L^{1}_{2j,j_2})^*,l_{j,1})&=-i( L^*_{2j-1,1},L_{2j,j_2}^*)=-sign(h_2(y_{2j}))=-1,
\end{split}
\end{equation}
$1<j_1\leq m, 1\leq j_2\leq m$.

Case 4: Construct $\{l_{j,1}$ :$j=I-1$ when $I$ is odd, or $j=I+1$ when $I$ is even$\}$. 

Let $l_{j,1}$ be a simple closed curve parallel to $\b_-^j$ and in a small regular neighborhood of $\b_-^j$ in $F_1$. 
Equip $l_{j,1}$ the orientation opposite to $L_{2j-1,1}^*$.
By (\ref{hn4}),
$h_2(y_{2j-1})=h_2(y_{2j})=-\bar{1}$,  $h_2(y_{2j-2})=h_2(y_{2j+1})=\bar{1}$.
By Remark \ref{arc}, 
\begin{equation*}
\begin{cases}
(L^1_{2j-1,j_1})^* \text{ and } (L^3_{2j-1,j_1})^* \text{ are from } \b_-^{j} \text{ to }\b_+^{j},1<j_1\leq m;\\ 
(L^3_{2j-2,j_2})^* \text{ is from } \b_-^{j} \text{ to } \b_-^{j-1}, 1\leq j_2\leq m;\\ 
(L^1_{2j,j_2})^* \text{ is from } \b_-^{j} \text{ to } \b_-^{j+1}, 1\leq j_2\leq m.
\end{cases}
\end{equation*}
Then $l_{j,1}$ intersects each of $(L^1_{2j-1,j_1})^*$, $(L^3_{2j-1,j_1})^*$, $(L^3_{2j-2,j_2})^*$ and $(L^1_{2j,j_2})^*$ once, $1<j_1\leq m, 1\leq j_2\leq m$.
By  (\ref{neworientation}),
\begin{equation}\label{c4}
\begin{split}
i((L^{1}_{2j-1,j_1})^*,l_{j,1})&= i( L^*_{2j-1,1},L_{2j-1,j_1}^*)=sign(h_2(y_{2j}))=-1\\
i((L^{3}_{2j-1,j_1})^*,l_{j,1})&= i( L^*_{2j-1,1},L_{2j-1,j_1}^*)=sign(h_2(y_{2j-1}))=-1\\
i((L^{3}_{2j-2,j_2})^*,l_{j,1})&= i( L^*_{2j-1,1},L_{2j-2,j_2}^*)=sign(h_2(y_{2j-1}))=-1,\\
i((L^{1}_{2j,j_2})^*,l_{j,1})&= i( L^*_{2j-1,1},L_{2j,j_2}^*)=sign(h_2(y_{2j}))=-1,
\end{split}
\end{equation}
$1<j_1\leq m, 1\leq j_2\leq m$.

Case 5: Construct $l_I$ and $l_{I,1}$.

Let $a_1$ and $a_2$ be two arcs in two disks with radius greater than $\e$ centered at $\hat{c}_{2(I-1)}$ and $\hat{c}_{2I+1}$ respectively.
$a_1$ is parallel to $L_{2(I-1)-1,1}^*$ and intersects $(L_{2(I-1)-1,j_1}^1)^*$ and $(L_{2(I-1),j_2}^1)^*$ once, $1< j_1\leq m, 1\leq j_2\leq m$.
$a_2$ is parallel to $L_{2(I+1)-1,1}^*$ and intersects $(L_{2(I+1)-1,j_1}^3)^*$ and $(L_{2I,j_2}^3)^*$ once, $1< j_1\leq m, 1\leq j_2\leq m$.
Let $b_1$ and $b_2$ be two arcs parallel to $L^*_{2I-1,1}$ in the disks with radius greater than $\e$ centered at $\hat{c}_{2I-1}$ and $\hat{c}_{2I}$ respectively. 
$b_1$ intersects $(L_{2I-1,j_1}^1)^*$ and $(L_{2(I-1),j_2}^1)^*$ once, and $b_2$ intersects $(L_{2I-1,j_1}^3)^*$ and $(L_{2I,j_2}^3)^*$ once, $1< j_1\leq m, 1\leq j_2\leq m$.
$a_1,a_2,b_1,b_2$ are shown as in Figure \ref{Lijstar6c} and Figure \ref{610}.

$l_I$ is obtained by connecting $a_1$ and $b_1$ by two arcs. 
One is parallel to $(L_{2(I-1),1}^1)^*$ and in a small neighborhood of $(L_{2(I-1),1}^1)^*$. The other is parallel to 
\begin{equation*}
\begin{cases}
(L_{2(I-1),(m+1)/2}^1)^* \text{ and in a small neighborhood of } (L_{2(I-1),(m+1)/2}^1)^*, &\text{if }h_2(y_{2I-1})=-\bar{1},\\
(L_{2(I-1),(m+3)/2}^1)^* \text{ and in a small neighborhood of } (L_{2(I-1),(m+3)/2}^1)^*, &\text{if }h_2(y_{2I-1})=\bar{1}.
\end{cases}
\end{equation*}
Similarly, $l_{I,1}$ is obtained by connecting $a_2$ and $b_2$ by two arcs. 
One is parallel to $(L_{2I,1}^3)^*$ and in a small neighborhood of $(L_{2I,1}^3)^*$. The other is parallel to
\begin{equation*}
\begin{cases}
(L_{2I,(m+1)/2}^3)^* \text{ and in a small neighborhood of } (L_{2I,(m+1)/2}^3)^*, &\text{if }h_2(y_{2I})=\bar{1},\\
(L_{2I,(m+3)/2}^3)^* \text{ and in a small neighborhood of } (L_{2I,(m+3)/2}^3)^*, &\text{if }h_2(y_{2I})=-\bar{1}.
\end{cases}
\end{equation*}
Note that the intersections of $l_I$ ($l_{I,1}$) and $L^*$ are contained in $a_1\cup b_1$ ($a_2\cup b_2$).
In addition, we can arrange $\{l_i,l_{j,1}:1\leq i\leq n/2, j=I-1,I,I+1\}$ such that they are mutually disjoint. 
Orient $l_{I}$ and $l_{I,1}$ such that $a_1$ and $a_2$ have the different orientations from $l_{I-1,1}$ and $l_{I+1,1}$ respectively. 
By the discussion in Case 3 and Case 4,
\begin{equation}\label{c5}
\begin{split}
i((L_{2(I-1)-1,j_1}^1)^*, a_1)&=-i((L_{2(I-1)-1,j_1}^1)^*,l_{I-1,1})=1;\\
i((L_{2(I-1),j_2}^1)^*, a_1)&=-i((L_{2(I-1),j_2}^1)^*,l_{I-1,1})=1;\\ 
i((L_{2(I+1)-1,j_1}^3)^*, a_2)&=-i((L_{2(I+1)-1,j_1}^3)^*,l_{I+1,1})=1; \\
i((L_{2I,j_2}^3)^*, a_2)&=-i((L_{2I,j_2}^3)^*,l_{I+1,1})=1,
\end{split}
\end{equation} $1<j_1\leq m, 1\leq j_2\leq m$.
By the construction of $l_I$ and $l_{I,1}$,
\begin{equation}\label{c6}
\begin{split}
i((L_{2(I-1),j_2}^1)^*, b_1)&=-i((L_{2(I-1),j_2}^1)^*,a_1)=-1;\\ 
i((L_{2I,j_2}^3)^*, b_2)&=-i((L_{2I,j_2}^3)^*,a_2)=-1,
\end{split}
\end{equation}
$1<j_1\leq m, 1\leq j_2\leq m$. 
 $b_1$ and $b_2$ are parallel to $L_{2I-1,1}^*$. $(L_{2I-1,j_1}^1)^*$ and $(L_{2(I-1),j_2}^1)^*$ intersect $L_{2l-1,1}^*$ in the same direction. $(L_{2I-1,j_1}^3)^*$ and $(L_{2I,j_2}^3)^*$ intersect $L_{2l-1,1}^*$ in the same direction. By (\ref{neworientation})
\begin{equation}\label{c7}
\begin{split}
i((L_{2I-1,j_1}^1)^*, b_1)&=i((L_{2(I-1),j_2}^1)^*,b_1)=-1;\\ 
i((L_{2I-1,j_1}^3)^*, b_2)&=i((L_{2I,j_2}^3)^*,b_2)=-1,
\end{split}
\end{equation} 
$1<j_1\leq m, 1\leq j_2\leq m$.

Next, we check $i((L^{2l-1}_{i,j})^*,{\cal C})<0$, $(i,j)\in A, l=1,2$. Assume that $n\geq 10$ (i.e. $I\geq 3$) and $1<j_1 \leq m, 1\leq j_2\leq m$. (The case $n=6$ is  similar. (c.f. Figure \ref{Lijstar6c}))

By (\ref{c1}) and (\ref{c2}),
\begin{align*}
i((L_{2i-1,j_1}^1)^*,{\cal C})&=i((L_{2i-1,j_1}^3)^*, {\cal C})=-1, 1\leq i\leq I-2, \text{or } I+2\leq i\leq n/2;\\
i((L_{2i,j_2}^1)^*,{\cal C})&=i((L_{2i,j_2}^3)^*, {\cal C})=-1, 1\leq i\leq I-3, \text{or } I+2\leq i\leq n/2.\\
i((L_{2(I-2),j_2}^1)^*,{\cal C})
&=i((L_{2(I-2),j_2}^1)^*,l_{I-1})=-1, \text{by (\ref{c1}) and (\ref{c2})}.\\
i((L_{2(I-2),j_2}^3)^*,{\cal C})
&=i((L_{2(I-2),j_2}^3)^*,l_{I-2})+i((L_{2(I-2),j_2}^3)^*,l_{I-1,1})\\
&=-1+(-1)=-2, \text{by (\ref{c1})-(\ref{c4})}.\\
i((L_{2(I-1)-1,j_1}^1)^*,{\cal C})&=i((L_{2(I-1)-1,j_1}^1)^*,l_{I-1})+i((L_{2(I-1)-1,j_1}^1)^*,l_{I-1,1})+i((L_{2(I-1)-1,j_1}^1)^*,a_1)\\
&=-1+(-1)+1=-1, \text{by (\ref{c1})-(\ref{c5})}.\\
i((L_{2(I-1)-1,j_1}^3)^*,{\cal C})
&=i((L_{2(I-1)-1,j_1}^3)^*,l_{I-1})+i((L_{2(I-1)-1,j_1}^3)^*,l_{I-1,1})\\
&=-1+(-1)=-2, \text{by (\ref{c1})-(\ref{c4})}.\\
i((L_{2(I-1),j_2}^1)^*, {\cal C})
&=i((L_{2(I-1),j_2}^1)^*,l_{I-1,1})+i((L_{2(I-1),j_2}^1)^*,a_1)+i((L_{2(I-1),j_2}^1)^*,b_1)\\
&=-1+1+(-1)=-1, \text{by (\ref{c3})-(\ref{c6})}.\\
i((L_{2(I-1),j_2}^3)^*, {\cal C})
&=i((L_{2(I-1),j_2}^3)^*,l_{I-1})=-1, \text{by (\ref{c1})-(\ref{c2})}.\\
i((L_{2I-1,j_1}^1)^*, {\cal C})
&=i((L_{2I-1,j_1}^1)^*,b_1)=-1;\\  
i((L_{2I-1,j_1}^3)^*, {\cal C})
&=i((L_{2I-1,j_1}^3)^*,b_2)=-1,\text{by (\ref{c7})}.\\
i((L_{2I,j_2}^1)^*,{\cal C})
&=i((L_{2I,j_2}^1)^*,l_{I+1})=-1, \text{by (\ref{c1}) and (\ref{c2})}.\\
i((L_{2I,j_2}^3)^*,{\cal C})
&=i((L_{2I,j_2}^3)^*,l_{I+1,1})+i((L_{2I,j_2}^3)^*,a_2)+i((L_{2I,j_2}^3)^*,b_2)\\
&=-1+1+(-1)=-1, \text{by (\ref{c3})-(\ref{c6})}. 
\end{align*}
\begin{align*}
i((L_{2(I+1)-1,j_1}^1)^*,{\cal C})
&=i((L_{2(I+1)-1,j_1}^1)^*, l_{I+1})+i((L_{2(I+1)-1,j_1}^1)^*, l_{I+1,1})\\
&=-1+(-1)=-2, \text{by (\ref{c1})-(\ref{c4})};\\
i((L_{2(I+1)-1,j_1}^3)^*,{\cal C})
&=i((L_{2(I+1)-1,j_1}^3)^*,l_{I+1})+i((L_{2(I+1)-1,j_1}^3)^*,l_{I+1,1})+i((L_{2(I+1)-1,j_1}^3)^*,a_2)\\
&=-1+(-1)+1=-1, \text{by (\ref{c1})-(\ref{c5})}.\\
i((L_{2(I+1),j_2}^1)^*,{\cal C})
&=i((L_{2(I+1),j_2}^1)^*,l_{I+2})+i((L_{2(I+1),j_2}^1)^*,l_{I+1,1})\\
&=-1+(-1)=-2, \text{by (\ref{c1})-(\ref{c4})}\\
i((L_{2(I+1),j_2}^3)^*,{\cal C})
&=i((L_{2(I+1),j_2}^3)^*,l_{I+1})=-1, \text{by (\ref{c1}) and (\ref{c2})}.
\end{align*}
$1<j_1 \leq m, 1\leq j_2\leq m$.

Then $i((L^{2l-1}_{i,j})^*,{\cal C})<0$, $(i,j)\in A, l=1,2$.
Let $\breve{\cal F}'_{1,s}$ be the new surface bundle of $\breve{M}_{1,s}$ after enough times Denh twist operations along $\G$, $s=1,2$.
By Claim \ref{torus}, there is an isotopy of $\breve{\cal F}'_{1,s}$ such that $\{\breve{L}_{i,j,t,s}^{2l-1}: (i,j)\in A, s,t,l=1,2\}$ are transverse to the new surface bundle.
Then the exterior of $\breve{L}$ in $\breve{Y}$ is a surface bundle. 
$\breve{Y}$ is a $p$-fold cover of $W_K$, which is a 2-fold cover of the exterior of $K$ is $S^3$. 
Thus $K$ is virtually fibered. 
We finish the proof of Theorem \ref{main} in Case II.

\end{document}